\documentclass[twoside]{book}

\newcommand{\arxiv}[2]{#1} 
\arxiv{}{\pdfminorversion=3}

\usepackage{problems}
\usepackage[colorlinks=true,
citecolor=black,
linkcolor=black,
anchorcolor=black,
filecolor=black,
menucolor=black,
urlcolor=black,
pdftitle={PIGTIKAL: puzzles in geometry that I know and love},
pdfsubject={Geometry},
pdfauthor={Anton Petrunin},
]{hyperref}

\makeindex

\arxiv{}{\usepackage[x-1a]{pdfx}}

\usepackage{attachfile}

\begin{document}


\begin{titlepage}
\fontfamily{cmr}\selectfont  
    \vspace*{4.5\baselineskip}
    \centering

\centering
\vspace*{3mm}
{\huge PIGTIKAL} \\[0em]
{\huge (puzzles in geometry\\
that I know and love)}\\[14mm]
{\Large Anton Petrunin} 
\par
\vspace*{55mm}
{\large Association for Mathematical Research\\[0mm]
Monographs Volume 2}
\par
\vspace*{14mm}
\arxiv{\includegraphics[angle=0, scale=.17]{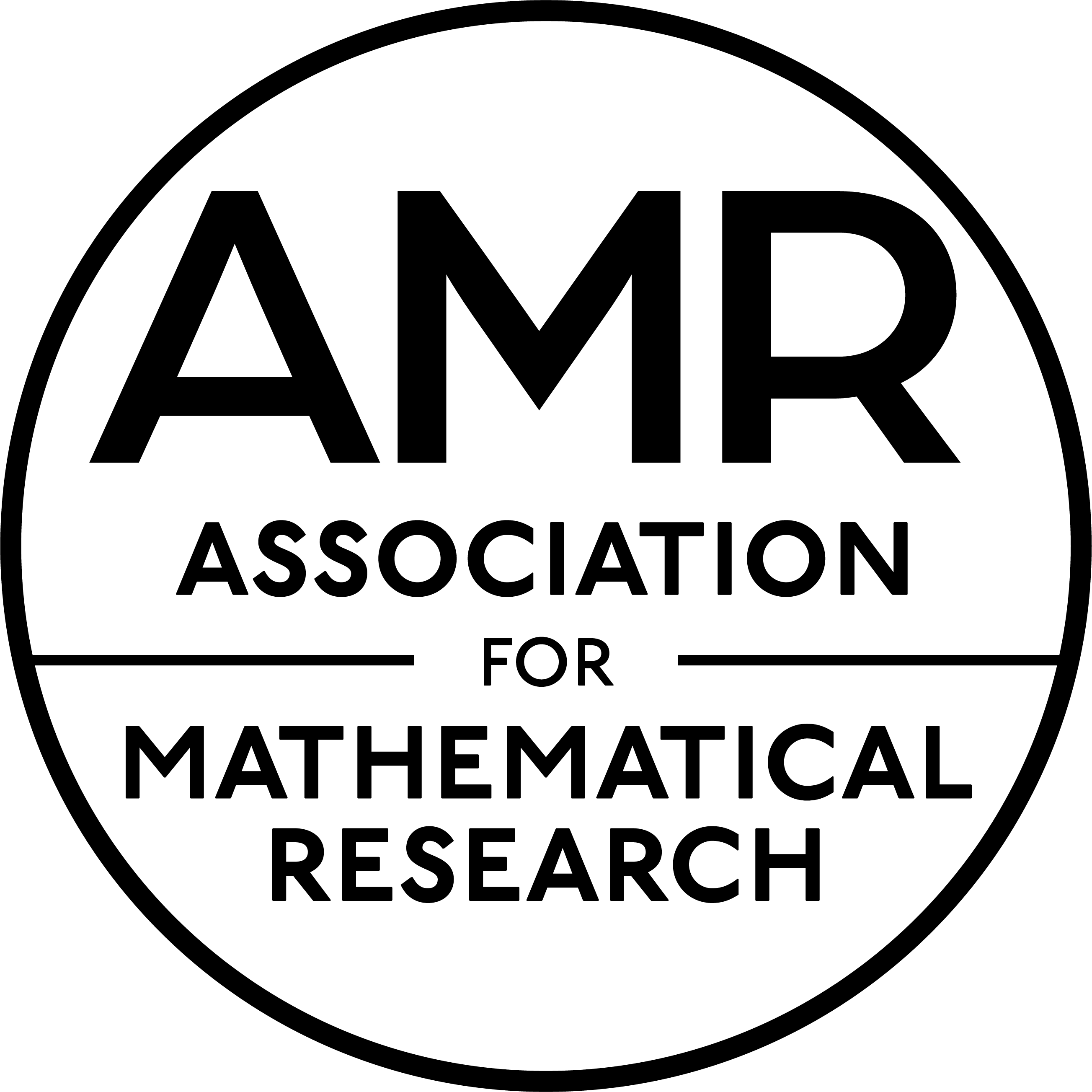}}{\includegraphics[angle=0, scale=.17]{cover/amr-small.jpg}}
\end{titlepage}

\noindent
Anton Petrunin\\
Department of Mathematics\\
The Pennsylvania State University\\
University Park, PA 16802\\
USA\\[20mm]
\noindent
AMR Monographs\\
Editorial Board:

\qquad George Andrews

\qquad Colin Adams

\qquad 
Eric Friedlander

\qquad 
Robert Ghrist

\qquad 
Joel Hass

\qquad 
Robion Kirby

\qquad 
Alex Kontorovich

\qquad 
Sergei Tabachnikov\\
[30mm]
\noindent{\includegraphics[scale=0.5]{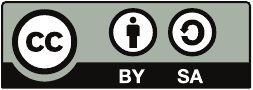}
\vspace*{1mm}
\\
Copyright © 2022 by Anton Petrunin\\
This work is licensed under a CC BY-SA 4.0 license\\
To view a copy of this license, visit:\\
\texttt{https://creativecommons.org/licenses/by-sa/4.0/legalcode}\\
\null\vfill
\noindent Second printing: 2025
\null\vfill
\noindent 
\texttt{arXiv:0906.0290v19}\\[3mm]
Association for Mathematical Research

\qquad Davis, CA; Jenkintown PA\\[3mm]
\noindent Cover design: Robert Ghrist

\thispagestyle{empty}
\tableofcontents
\thispagestyle{empty}

\newpage
\thispagestyle{empty}
\section*{Preface}

This collection is about ideas, and it is not about theory.
An idea might feel more comfortable in a suitable theory,
but it has its own life and history, and it can speak for itself.

I am collecting these problems for fun, 
but they might be used to improve 
the problem-solving skills in geometry.
Every problem has a short elegant solution ---
this gives a hint which was not available
when the problem was discovered.


\parbf{How to read it.}
Open a random chapter; make sure you like the practice problem --- if yes try to solve a random problem in the chapter.
A semisolution is given at the end of the chapter,
but think before reading,
otherwise, it will not help. 

Some problems are marked by $\circ$, $*$, $+$ or $\sharp$.
\begin{itemize}
\item[$\circ$] --- easy problem;
\item[$*$] --- the solution requires at least two ideas;
\item[$+$] --- the solution requires knowledge of a theorem;
\item[$\sharp$] --- there are interesting solutions based on different ideas.
\end{itemize}

\parbf{Acknowledgments.} 
I want to thank everyone who helped me;
here is an incomplete list:
Stephanie Alexander,
Ilya Alexeev,
Miroslav Ba\v{c}\'{a}k, 
Christopher Croke,
Bogdan Georgiev,
Sergei Gelfand,
Mohammad Ghomi,
Jouni Luukkainen,
Alexander Lytchak,
Andrei Malyutin, 
Rostislav Matveyev, 
Dmitri Panov, 
Peter Petersen, 
Idzhad Sabitov,
Thomas Sharpe,
Serge Tabachnikov, and
Sergio Zamora Barrera.

This collection is partly inspired by \emph{connoisseur's collection} of puzzles of Peter Winkler \cite[][]{winkler}.
Many problems were suggested on MathOverflow~\cite[][]{One-step}.

This work was partially supported by the following grants:
NSF grants DMS 
0103957,
0406482,
0905138,
1309340,
2005279,
Simons Foundation grants 
245094 and 584781,
and
Minobrnauki grant 075-15-2022-289.


\csname @openrightfalse\endcsname
\chapter{Curves}

Recall that a \index{curve}\emph{curve} is a continuous map 
from a real interval into a space (for example, Euclidean plane)
and 
a {}\emph{closed curve} is a continuous map defined on a circle.
If the map is injective then the curve is called {}\emph{simple}.

We assume that the reader is familiar with related definitions including 
length of curve 
and its curvature.
The necessary material is covered in the first couple of lectures 
of a standard introduction to differential geometry, [see Part I in \ncite{petrunin-zamora}, Chapter 1 in \ncite{toponogov-curves-and-surfaces}, or \S26--27 in \ncite{hilbert-cohn-vossen}].

\medskip

We give a practice problem with a solution --- after that, you are on your own.

\subsection*{Spiral}
\label{spiral}

The following problem states that 
if you drive on the plane and turn the steering wheel to the right all the time,
then you will not be able to come back to the same place.

\begin{pr}
Let $\gamma$ be a smooth regular plane curve with strictly monotonic curvature. 
Show that $\gamma$ has no self-intersections.
\end{pr}

\begin{wrapfigure}{o}{35 mm}
\centering
\includegraphics{mppics/pic-102}
\end{wrapfigure}

\parit{Semisolution.}
The trick is to show that the osculating circles of $\gamma$ are nested.

\medskip

Without loss of generality, we may assume that the curve is parametrized by its length and its
curvature decreases.

Let $z(t)$ be the center of the osculating circle at $\gamma(t)$
and $r(t)$ its radius.
Note that 
\begin{align*}
z(t)&=\gamma(t)+\tfrac{\gamma''(t)}{|\gamma''(t)|^2},
&
r(t)&=\tfrac{1}{|\gamma''(t)|}.
\end{align*}

Straightforward calculations show that
\[|z'(t)|= r'(t).\]
Note that the curve $z(t)$ has no straight arcs;
therefore 
\[|z(t_1)-z(t_0)|<r(t_1)-r(t_0).\leqno({*})\]
if $t_1>t_0$.

Denote by $D_t$ the osculating disk of $\gamma$ at $\gamma(t)$;
it has a center at $z(t)$ and radius $r(t)$.
By $({*})$, $D_{t_1}$ lies in the interior of $D_{t_0}$ for any $t_1>t_0$.
Hence the result follows.\qeds

This problem was considered by Peter Tait \cite{tait}
and later rediscovered by Adolf Kneser \cite{kneser}.
The osculating circles of the curve give a peculiar decomposition of an annulus into circles; it has the following property: if a smooth function is constant on each osculating circle it must be constant in the annulus \cite[see Lecture 10 in][]{fuchs-tabachnikov}.
The same idea leads to a solution of the following problem: 

\begin{pr}
Let $\gamma$ be a smooth regular plane curve with strictly monotonic curvature. 
Show that no line can be tangent to $\gamma$ at two distinct points.
\end{pr} 

It is instructive to check that the 3-dimensional analog does not hold;
that is, there are self-intersecting smooth regular space curves with strictly monotonic curvature. 

Note that if the curve $\gamma(t)$ is defined for $t\in[0,\infty)$ and its curvature tends to $\infty$ as $t\to \infty$, 
then the problem implies the existence of the limit of $\gamma(t)$ as $t\to\infty$.
The latter result could be considered as a continuous analog of the Leibniz test for alternating series.

{

\begin{wrapfigure}{r}{35 mm}
\vskip2mm
\centering
\includegraphics{mppics/pic-104}
\end{wrapfigure}

\subsection*{Moon in a puddle}
\label{moon-in-puddle}

\begin{pr}
A smooth closed simple plane curve with curvature less than $1$ at every point bounds figure~$F$. 
Prove that $F$ contains a unit disk.
\end{pr}

}

\subsection*{Wire in a tin}
\label{A spring in a tin} 

\begin{pr}
Let $\alpha$ be a closed smooth curve immersed
in a unit disk. 
Prove that the average absolute curvature of $\alpha$ is at least $1$, with
equality if and only if $\alpha$ is the unit circle possibly traversed more than once.
\end{pr}

\subsection*{Curve on a sphere}
\label{A curve in a sphere}

\begin{pr}
Show that if a closed curve on the unit sphere intersects every equator then its length is at least $2\cdot\pi$.
\end{pr}


{

\begin{wrapfigure}{r}{25 mm}
\vskip2mm
\centering
\includegraphics{mppics/pic-106}
\end{wrapfigure}

\subsection*{Oval in an oval}
\label{Oval in oval}

\begin{pr}
Consider two closed smooth strictly convex planar curves, one inside the other. 
Show that there is a chord of the outer curve that is tangent to the inner curve at the midpoint of the chord.
\end{pr}

}

\subsection*{Capture a sphere in a knot\hard}
\label{Capture a sphere in a knot}

The following formulation uses the notion of smooth isotopy of knots,
that is, a one-parameter family of embeddings 
\[f_t\:\mathbb{S}^1\z\to \RR^3,\ \ t\in[0,1]\] 
such that the map $[0,1]\times \mathbb{S}^1\to\RR^3$ is smooth.

\begin{pr}
Show that one cannot capture a sphere in a knot.

More precisely, let $B$ be the closed unit ball in $\RR^3$
and $f\:\mathbb{S}^1\z\to \RR^3\setminus B$ a knot.
Show that there is a smooth isotopy 
$$f_t\:\mathbb{S}^1\to \RR^3\setminus B,\ \ \ t\in [0,1]$$ 
such that $f_0=f$,
the length of $f_t$ is non-increasing with respect to $t$
and $f_1(\mathbb{S}^1)$ can be separated from $B$ by a plane.
\end{pr}

{

\begin{wrapfigure}{r}{20 mm}
\vskip2mm
\centering
\includegraphics{mppics/pic-108}
\end{wrapfigure}

\subsection*{Linked circles}
\label{linked-circles}

\begin{pr}
Suppose that two linked simple closed curves in $\RR^3$
lie at a distance at least $1$ from each other.
Show that the length of each curve is at least $2\cdot\pi$.
\end{pr}

}

\subsection*{Surrounded area}
\label{Surrounded area}

\begin{pr}
Consider two simple closed plane curves 
$\gamma_1,\gamma_2\:\mathbb S^1\to\RR^2$.
Assume 
\[|\gamma_1(v)-\gamma_1(w)|\le|\gamma_2(v)-\gamma_2(w)|\]
for any $v,w\in \mathbb S^1$.
Show that the area surrounded by $\gamma_1$ does not exceed the area surrounded by $\gamma_2$. 
\end{pr}

\subsection*{Crooked circle}

\label{Crooked circle}

\begin{pr}
Construct 
a bounded set in $\RR^2$
homeomorphic to an open disk
such that 
its boundary contains no simple curves.
\end{pr}

\subsection*{Rectifiable curve}
\label{Rectifiable curve}

For the following problem we need the notion of 
\index{Hausdorff measure}\emph{Hausdorff measure}.
Choose a compact set $X\subset\RR^2$ and $\alpha>0$.
Given $\delta>0$, set
\[h(\delta)=\inf\left\{\,\sum_i(\diam X_i)^\alpha\,\right\},\]
where the greatest lower bound is taken over all finite coverings $\{X_i\}$ of $X$ 
such that $\diam X_i<\delta$ for each $i$.

Note that the function $\delta\mapsto h(\delta)$ is not decreasing in $\delta$.
In particular, $h(\delta)\to \mathcal{H}_\alpha(X)$ as $\delta\to 0$ for some (possibly infinite) value $\mathcal{H}_\alpha(X)$.
This value $\mathcal{H}_\alpha(X)$ is called the $\alpha$-dimensional Hausdorff measure of $X$.

\begin{pr}
Let $X\subset \RR^2$ be a compact connected set
with finite 1-dimensional Hausdorff measure. 
Show that there is a rectifiable curve passing thru all the points in $X$.
\end{pr}

\subsection*{Shortcut\hard}

\begin{pr}
Let $X\subset \RR^2$ be a compact connected set. 
Show that any two points $x,y\in X$ can be connected by a path $\alpha$ such that the complement $\alpha\setminus X$ has arbitrarily small length.
\end{pr}

Note that it might be impossible to connect $x$ and $y$ by a path in~$X$.
In fact, there are connected sets in the plane (for example, the pseudo-arc) that contain no curves.

\subsection*{Straight set}

A set $X$ in the plane is called $\delta$-straight if, for any disc $\DD(x,r)$ with radius $r>0$ and center at a point $x\in X$, the intersection 
$X\cap \DD(x,r)$ is $\delta\cdot r$-close in the sense of Hausdorff to a diameter of $\DD(x,r)$. 

\begin{pr}
Show that for any $\eps>0$ there is $\delta>0$ such that any $\delta$-straight closed plane set $X$ is a curve that admits a locally $(1\pm\eps)$-bi-Hölder parametrization.
That is, $X$ is connected and for any $x\in X$ there is a curve $\alpha\:[0,1]\to \RR^2$ that covers a neighborhood of $x$ in $X$, and the inequality
\[c\cdot |t_1-t_0|^{1+\eps}\le|\alpha(t_0)-\alpha(t_1)| \le C\cdot |t_1-t_0|^{1-\eps}\]
holds for any $t_0$ and $t_1$ and positive constants $c$ and $C$,
\end{pr}

Before trying to solve the problem, it might be useful to look at the following example.
Consider the set $X_k$ formed by the origin and the two logarithmic spirals $\rho = \pm e^{k\varphi}$ in the polar coordinates $(\rho,\phi)$.
Observe that if $k$ is large, then $X_k$ is $\delta$-straight.

\noindent
\hskip-.9in
\arxiv{\includegraphics[scale=1]{mppics/pic-111}}{\includegraphics[scale=1]{mppics/pic-110}}

\vskip-.3in
Another example of a $\delta$-straight set is shown; it is constructed in a way similar to the Koch snowflake curve, where we used a very obtuse isosceles triangle instead of equilateral.
The Hausdorff dimension of this example is larger than 1;
in particular, it shows that one can not expect to have a Lipschitz parametrization of $X$.

\subsection*{Typical convex curves}

Formally we do not need it in the problem, 
but it is worth noting that the curvature of a convex curve is defined almost everywhere;
it follows from the fact that monotonic functions are differentiable almost everywhere.

\begin{pr}
Show that \emph{most} of the convex closed curves in the plane
have vanishing curvature at every point where it is defined.
\end{pr}

We need to explain the meaning of the word ``most'' in the formulation;
it uses \index{Hausdorff distance}\emph{Hausdorff distance} and \index{G-delta set}\emph{G-delta sets}.

The Hausdorff distance $|A-B|_H$ between two closed bounded sets $A$ and $B$ in the plane is defined by 
\[|A-B|_H=\sup_{x\in\RR^2}\{|\dist_A(x)-\dist_B(x)|\},\]
where $\dist_A(x)$ denotes the smallest distance from $A$ to $x$.
Equivalently, $|A-B|_H$ can be defined as the greatest lower bound of the positive numbers $r$ such that the $r$-neighborhood of $A$ contains $B$ and the $r$-neighborhood of $B$ contains $A$.

It is straightforward to show that the Hausdorff distance defines a metric on the space of all closed plane curves.
The obtained metric space is locally compact.
The latter follows from the \index{selection theorem}\emph{selection theorem} \cite[see \S18 in][]{blaschke},
which states that closed subsets of a fixed closed bounded set in the plane form a compact set with respect to the Hausdorff metric. 

A G-delta set in a metric space $X$ is defined as a countable intersection of open sets.
According to the \index{Baire category theorem}\emph{Baire category theorem}, 
in locally compact metric spaces $X$,
the intersection of a countable collection of open dense sets 
has to be dense.
(The same holds if $X$ is complete, but we will not need it.)

In particular, in $X$, 
the intersection of a finite or countable collection of G-delta dense sets is also a G-delta dense set. 
It means that each G-delta dense set contains {}\emph{most} of $X$.
This is the meaning of the word {}\emph{most} used in the problem.

\section*{Semisolutions}


\parbf{Moon in a puddle.}
In the proof we will use the \emph{cut locus}
of $F$ with respect to its boundary (also known as \emph{medial axis}.);
it will be further denoted by $T$.
The cut locus can be defined as the closure
of the set of points $x\in F$ 
for which there exist two or more points in $\partial F$ minimizing the distance to $x$.

\begin{wrapfigure}{r}{36 mm}
\vskip-2mm
\centering
\includegraphics{mppics/pic-112}
\end{wrapfigure}

For each point $x\in T$, consider the subset $X\subset\partial F$ where the minimal distance to $x$ is attained.
If $X$ is not connected then we say that $x$ is a \emph{cut point};
equivalently it means that for any sufficiently small neighborhood $U\ni x$, 
the complement $U\setminus T$ is disconnected.
If $X$ is connected 
then we say that $x$ is a \emph{focal point};
equivalently it means that the osculating circle to $\partial F$ at any point of $X$ is centered at $x$.

The trick is to show that $T$ contains a focal point, say $z$.
Since $\partial F$ has curvature of at most $1$, the radius of any osculating circle is at least $1$.
Hence the distance from $\partial F$ to $z$ is at least 1,
and the statement will follow.

\medskip

After a small perturbation of $\partial F$, we may assume that $T$ is a graph embedded in $F$ with a finite number of edges.

Note that $T$ is a
deformation retract of $F$.
The retraction $F\to T$ can be obtained the following way:
(1) given a point $x\in F\setminus T$,
consider the (necessarily unique) point $\hat x\in \partial F$ that minimizes the distance $|x-\hat x|$ and
(2) move $x$ along the extension of the line segment $[\hat x x]$ behind $x$ until it hits $T$.

In particular, $T$ is a tree.
Therefore $T$ has
an end vertex, say $z$.
The point $z$ is focal since there are arbitrarily small neighborhoods $U$ of $z$ such that the complements $U\setminus T$ are connected.
\qeds

We proved a slightly stronger statement;
namely, \textit{there are at least two points on $\partial F$ at which osculating circles lie in $F$}.
Note that these points are \index{vertex of curve}\emph{vertices} of $\partial F$;
that is, they are critical points of its curvature.

Note further that inversion respects osculating circles.
That is, if $\gamma$ is an osculating circle of curve $\alpha$ at $t_0$,
$\gamma'$ is the inversion of $\gamma$, and 
$\alpha'$ is the inversion of $\alpha$,
then $\gamma'$ is an osculating circle of curve $\alpha'$ at $t_0$.
Therefore applying an inversion about a circle with the center in $F$, we also get a pair of osculating circles of $\partial F$ which surround $F$.
This way we obtain 4 osculating circles that lie on one side of $\partial F$.
The latter statement is a generalization of the four-vertex theorem \cite{petrunin-zamora:moon}.

The case of convex curves of this problem appears in the book of Wilhelm Blaschke \cite[see \S 24 in][]{blaschke}.
In full generality, the problem was discussed by Vladimir Ionin and German Pestov \cite{pestov-ionin}. 
A solution via curve shortening flow of a weaker statement 
was given by Konstantin Pankrashkin \cite{pankrashkin}.
The statement still holds if the curve fails to be smooth at one point \cite{petrunin-zamora:moon}.
A spherical version of the later statement 
was used by Dmitri Panov and me \cite{panov-petrunin-ramification}.

The statement admits a generalization to curves that bound a disc $F$ in a surface with nonpositive curvature.
The latter can be used to prove the following problem which was suggested by Dmitri Burago.

\begin{pr}
Let $\gamma$ be a closed space curve with curvature at most 1.
Show that $\gamma$ cannot be filled by a disc with area less than $\pi$.
\end{pr}

As you can see from the following problem, the 3-dimensional analog of this statement does not hold.

\begin{pr}
Construct a smooth embedding $\mathbb{S}^2\hookrightarrow \RR^3$ 
with all the principal curvatures between $-1$ and $1$
such that it does not surround a ball of radius~1.
\end{pr}

Such an example can be obtained by fattening a nontrivial contractible 2-complex in $\RR^3$ 
[Bing's house constructed in \ncite{bing} will do the job].
This problem is discussed by Abram Fet and Vladimir Lagunov \cite{lagunov-2,lagunov-fet-1,lagunov-fet-2} 
and it was generalized to Riemannian manifolds with boundary by Stephanie Alexander and Richard Bishop \cite{alexander-bishop}.

A similar argument shows that for any Riemannian metric $g$ on the 2-sphere $\mathbb S^2$ 
and any point $p\in(\mathbb S^2,g)$ there is a minimizing geodesic $[pq]$ with conjugate ends.
On the other hand, for $(\mathbb S^3,g)$ this is not true.
Moreover, there is a metric $g$ on $\mathbb{S}^3$ 
with sectional curvature bounded above by arbitrarily small $\eps>0$ and $\diam(\mathbb{S}^3,g)\le 1$.
In particular, $(\mathbb S^3,g)$ has no minimizing geodesic with conjugate ends.
An example was originally constructed by Mikhael Gromov \cite{gromov-almost-flat}; 
a simplification was given by 
Peter Buser
and Detlef Gromoll \cite{buser-gromoll}.

\parbf{Wire in a tin.}
To solve this problem,
you should imagine that you travel on a train along the curve $\alpha(t)$
and watch the position of the center of the disk in the frame of your train car.

\medskip

Denote by $\ell$ the length of $\alpha$.
Equip the plane with complex coordinates so that $0$ is the center of the unit disk.
We can assume that $\alpha$ is equipped with an $\ell$-periodic parametrization by arc length.

Consider the curve $\beta(t)=t-\tfrac{\alpha(t)}{\alpha'(t)}$.
Observe that 
\[\beta(t+\ell)=\beta(t)+\ell\] 
for any $t$.
In particular 
\[\length (\beta|_{[0,\ell]}) 
\ge 
|\beta(\ell)-\beta(0)|
=
\ell.\leqno({*})\]

Also 
\begin{align*}
|\beta'(t)|&=|\tfrac{\alpha(t)\cdot\alpha''(t)}{\alpha'(t)^2}|\le
\\
&\le|\alpha''(t)|.
\end{align*}
Note that $|\alpha''(t)|$ is the absolute curvature of $\alpha$ at $t$.
Therefore, the result follows from~$({*})$.
\qeds

The statement was originally proved 
by Istv\'an F\'ary in \cite{fary};
several different proofs are discussed by Serge Tabachnikov [see \ncite{tabachnikov} and also 19.5 in \ncite{fuchs-tabachnikov}].

If instead of the disk we have a region bounded by a closed convex curve $\gamma$, 
then it is still true that the average curvature of $\alpha$ is at least as big as the average curvature of $\gamma$. 
The proof was given by Jeffrey Lagarias
and Thomas Richardson [see \ncite{lagarias-richardson} and also \ncite{nazarov-petrov}].

Our solution can be adapted to the unit ball of arbitrary dimension.
Further, the same argument together with Liouville's theorem (geodesic flow preserves the phase volume) implies that a closed smooth submanifold in a unit ball has average normal curvatures at least~1.
Let us formulate two more related problems.

\begin{pr}
Let $M$ be a closed smooth $n$-dimensional submanifold in the unit ball in $\RR^q$.
Denote by $H_p$ the mean curvature vector of $M$ at $p$.
Show that the average value of $|H_p|$ is at least $n$.
\end{pr}

\begin{pr}
Suppose $T$ is a smoothly embedded 2-torus in the unit ball in $\RR^q$.
Show that the average value of squared normal curvatures of $T$ is at least $\tfrac32$.
\end{pr}

See also the problem ``Small-twist embedding'' on page \pageref{Small-twist}.

\parbf{Curve on a sphere.} 
Let us present two solutions.
We assume that $\alpha$ is a closed curve in $\mathbb{S}^2$ of length $2\cdot\ell$ that intersects each equator.

\parit{A solution with Crofton's formula.}
Given a unit vector $u$, denote by $e_u$ the equator with the pole at $u$.
Let $k(u)$ be the number of intersections
of $\alpha$ and $e_u$.

Note that for almost all $u\in \mathbb{S}^2$, the value $k(u)$ is even or infinite.
Since each equator intersects $\alpha$, we get $k(u)\ge 2$ for almost all $u$.

Then we get
\begin{align*}
2\cdot\ell&=\tfrac14\cdot\int\limits_{u\in \mathbb{S}^2}k(u)\ge 
\\
&\ge\tfrac12\cdot\area\mathbb{S}^2=
\\
&=2\cdot\pi.
\end{align*}

The first identity above is called \index{Crofton's formula}\emph{Crofton's formula}.
To prove this formula, start with the case when the curve is formed by one geodesic segment, summing up we get it for broken lines, and by approximation it holds for all curves.
\qeds

\parit{A solution by symmetry.}
Let $\check\alpha$ be a sub-arc of $\alpha$ of length $\ell$, with endpoints $p$ and $q$. 
Let $z$ be the midpoint of a minimizing geodesic $[pq]$ in~$\mathbb{S}^2$. 

Let $r$ be a point of intersection of $\alpha$ with the equator with pole at~$z$. 
Without loss of generality, we may assume that $r\in\check\alpha$. 

The arc $\check\alpha$ together with its reflection with respect to the point $z$ forms a closed curve of length $2\cdot \ell$ passing thru both $r$ and its antipodal point $r^{*}$.
Therefore 
\[\ell=\length \check\alpha\ge |r-r^{*}|_{\mathbb S^2}=\pi.\]
Here $|r-r^{*}|_{\mathbb S^2}$ 
denotes the angle metric in the sphere $\mathbb S^2$.\qeds

Different solutions of this problem are discussed in a short note by Robert Foote \cite{foote};
the second proof is due to Stephanie Alexander.
The problem was suggested by Nikolai Nadirashvili.
It is nearly equivalent to Fenchel's theorem: 

\begin{pr}
Show that total curvature of any closed smooth regular space curve is at least $2\cdot\pi$.
\end{pr}

Let us also mention the problem of Karol Borsuk that was solved by John Milnor and István Fáry;
it states that any embedded circle of total curvature less than $4\cdot\pi$ is unknotted.
Six proofs of this statement are surveyed by Stephan Stadler and the author \cite{petrunin-stadler-Fary-Milnor}.

\parbf{Oval in an oval.}
Choose a chord that divides the area of the bigger oval in the minimal (or maximal) ratio.

If the chord is not divided into equal parts, then you can rotate it slightly
to decrease the ratio.
Hence the problem follows.
\qeds

\begin{wrapfigure}{r}{25 mm}
\vskip-4mm
\centering
\includegraphics{mppics/pic-107}
\end{wrapfigure}

\parit{Alternative solution.}
Given a unit vector $u$, denote by $x_u$ the point on the inner curve
with the outer normal vector $u$.
Draw a chord of the outer curve that is tangent to the inner curve at $x_u$;
denote by $r\z=r(u)$ and $l=l(u)$ the lengths of the segments of this chord to the right and to the left of $x_u$, respectively.

Arguing by contradiction, assume that $r(u)\z\ne l(u)$ for all $u\in\mathbb{S}^1$.
Since the functions $r$ and $l$ are continuous,
we can assume that 
$$r(u)<l(u)\ \ \text{for all}\ \ u\in\mathbb{S}^1.\leqno{({*})}$$

Prove that
each of the following two integrals 
\begin{align*}
\tfrac12\cdot\int\limits_{u\in\mathbb{S}^1}r^2(u)
\quad\text{and}\quad
\tfrac12\cdot\int\limits_{u\in\mathbb{S}^1}l^2(u)
\end{align*}
give the area between the curves.
In particular, 
the integrals are equal. 
The latter contradicts $({*})$.\qeds

This is a problem of Serge Tabachnikov \cite{tabachnikov-mi}.
A closely related {}\emph{equal tangents problem} is discussed by the same author in \cite{tabachnikov-tan}.

\parbf{Capture a sphere in a knot.}
We can assume that the knot lies on the sphere $\partial B$.

Choose a Möbius transformation 
$m\:\mathbb{S}^2\to\mathbb{S}^2$ close to the identity and not a rotation.

Note that $m$ is a conformal map;
that is, there is a function $u$ defined on $\mathbb{S}^2$ 
as 
\[u(x)=\lim_{y,z\to x}\frac{|m(y)-m(z)|}{|y-z|}.\]
(The function $u$ is called the \emph{conformal factor} of $m$.)

Applying the area formula for $m$,
we get 
$$\frac1{\area \mathbb{S}^2}\cdot\int\limits_{u\in\mathbb{S}^2} u^2=1.$$ 
By Bunyakovsky inequality, 
$$\frac1{\area \mathbb{S}^2}\cdot\int\limits_{u\in\mathbb{S}^2} u<1.$$ 

It follows that after a suitable rotation of $\mathbb{S}^2$, 
the map $m$ decreases the length of the knot.

Iterating this construction we get a sequence of knots $f_n\:\mathbb{S}^1\z\to\mathbb{S}^2$ with length decreasing  and tending to zero.
Passing to the limit as $m\to\id$, we get a continuous one-parameter family of Möbius transformations which shorten the length of the knot.
Therefore it drifts the knot to a single hemisphere and allows the ball to escape. 
\qeds

This is a problem by Zarathustra Brady, 
the given solution is based on the idea of David Eppstein \cite{zeb}.
A solution to the following problem is based on the same idea.

\begin{pr}
Show that a sphere cannot be captured in a link with 3 components, but can be captured in a link with 4 components.
\end{pr}

It seems unknown which convex bodies can be captured by a knot.
The following problem is quite tricky already.

\begin{pr}
Show that one can capture some convex bodies in a knot.
\end{pr}

The following problem of Abram Besicovitch is closely related \cite{besicovitch-sphere}; it can be solved using spherical Crofton's formula.

\begin{pr}
Show the total length of strings in a net that can hold a unit sphere has to be larger than $3\cdot\pi$.
\end{pr}


\parbf{Linked circles.} 
Denote the linked circles by $\alpha$ and $\beta$. 

Choose a point $x\in\alpha$. 
Note that there is a point $y\in\alpha$ such that the line segment 
$[xy]$ intersects $\beta$, say at $z$. 
Indeed, if this is not the case, 
then a homothety with center $x$ to $\alpha$ would shrink it to $x$ without crossing $\beta$.
The latter contradicts that $\alpha$ and $\beta$ are linked. 

\begin{wrapfigure}{r}{26 mm}
\vskip-2mm
\centering
\includegraphics{mppics/pic-109}
\end{wrapfigure}

Let $\alpha^{*}$ be the image of $\alpha$ under the central projection onto the unit sphere around $z$.
Since $|\alpha(t)-z|\ge1$ for any $t$, we have that
$$\length \alpha\ge \length\alpha^{*}.$$

Observe that $\alpha^{*}$ is a closed spherical curve that contains two antipodal points,
one corresponds to $x$ and the other to $y$.
It follows that
$$\length \alpha^{*}\ge 2\cdot\pi.$$
Hence the result follows.\qeds

This problem was proposed by Frederick Gehring \cite[see 7.22 in][]{gehring};
solutions and generalizations are surveyed in \cite{mateljevic}. 
The presented solution is attributed to Marvin Ortel in \cite{CJKSW} and it is close to the solution of Michael Edelstein and Binyamin Schwarz \cite{edelstein-schwatz}.

\parbf{Surrounded area.} The trick is to use the Kirszbraun theorem:
\emph{Any $L$-Lipschitz  map $f\:Q\to \RR^n$ defined on a subset $Q\subset \RR^m$ can be extended to an $L$-Lipschitz  map $\bar f\:\RR^m\to \RR^n$.}

This theorem  appears in the thesis of Mojżesz Kirszbraun \cite[][]{kirszbraun};
it was rediscovered later by Frederick Valentine \cite[][]{valentine}.
An interesting survey is given by 
Ludwig Danzer, Branko Grünbaum and Victor Klee \cite[][]{danzer-grunbaum-klee}.

\medskip

Let $C_1$ and $C_2$ be the compact regions bounded by $\gamma_1$ and $\gamma_2$ respectively.

By the Kirszbraun theorem, there is a 1-Lipschitz map $f\:\RR^2\to\RR^2$ 
such that $f(\gamma_2(v))\z=f(\gamma_1(v))$ for any $v\in\mathbb S^1$.

Note that $f(C_2)\supset C_1$.
Hence the statement follows.\qeds

\parbf{Crooked circle.}
A continuous function $f\:[0,1]\to [0,1]$
will be called $\eps$-crooked 
if $f(0)=0$, $f(1)=1$ 
and for any segment $[a,b]\subset [0,1]$ 
one can choose $a\le x\le y\le b$ 
such that
\[|f(y)-f(a)|\le\eps\ \ \text{and}\ \ |f(x)-f(b)|\le\eps.\]

A sequence of $\tfrac1n$-crooked maps can be constructed recursively. 
Figure out the construction by looking at the following diagram.

\begin{figure}[ht!]
\centering
\includegraphics{mppics/pic-113}
\end{figure}

Now, start with the unit circle, 
$\gamma_0(t)=(\cos 2\pi t,\sin 2 \pi t)$.
Choose a sequence of positive numbers $\eps_n$ converging to zero very fast. 
Construct recursively a sequence of simple closed curves $\gamma_n\:[0,1]\z\to\RR^2$ such that $\gamma_{n+1}$ runs outside of the disk bounded by $\gamma_n$
and 
\[|\gamma_{n+1}(t)-\gamma_n\circ f_n(t)|<\eps_n,\]
for an $\eps_n$-crooked function $f_n$.
(It is hard to draw $\gamma_{n+1}$; it runs like crazy back-and-forth almost along $\gamma_n$.)

Denote by $D$ the union of all disks bounded by the curves $\gamma_n$.
Clearly $D$ is homeomorphic to an open disk.
For the right choice of the sequence $\eps_n$, 
the set $D$ is bounded.
By construction, the boundary of $D$ contains no simple curves. \qeds

In fact, the only curves in the boundary of the constructed set are constant. Compare to the problem \emph{Simple path} on page~\pageref{Simple path}.

The proof uses the so-called \emph{pseudo-arc} 
constructed by Bronis\l{}aw Knaster \cite{knaster}.
The proof resembles the construction of Cantor's set.
Here are a few similar problems:

\begin{pr}
 Construct three distinct open sets in $\RR$ with the same boundary.
\end{pr}

\begin{pr}
 Construct three open disks in $\RR^2$ having the same boundary.
\end{pr}

These disks are called \index{lakes of Wada}\emph{lakes of Wada}; it is described by Kuniz\^{o} Yo\-ne\-ya\-ma \cite{yoneyama}.

\begin{pr}
 Construct a Cantor set in $\RR^3$ with a non-simply-connected complement.
\end{pr}

This example is called \index{Antoine's necklace}\emph{Antoine's necklace} \cite{antoine}.

\begin{pr}
 Construct an open set in $\RR^3$ with a fundamental group isomorphic to the additive group of rational numbers.
\end{pr}

More advanced examples include
the \emph{Whitehead manifold}, 
\emph{Dogbone space}, 
and \emph{Casson handle};
see also the problem ``Conic neighborhood'' on page \pageref{Conic neighborhood}.

\parbf{Rectifiable curve.}
The 1-dimensional Hausdorff measure will be denoted by $\mathcal{H}_1$. 

Set $L=\mathcal{H}_1(K)$.
Without loss of generality, we may assume that $K$ has diameter $1$.

Since $K$ is connected, we get 
\[\mathcal{H}_1(B(x,\eps)\cap K)\ge\eps\leqno(*)\]
for any $x\in K$ and $0<\eps<\tfrac12$.

Let $x_1,\dots, x_n$ be a maximal set of points in $K$ with 
\[|x_i-x_j|\z\ge\eps\] for all $i\ne j$. 
From $(*)$ we have $n\le2\cdot L/\eps$.

Note that there is a tree $T_\eps$ with vertices $x_1,\dots, x_n$ and straight edges with lengths at most $2\cdot\eps$ each.
Therefore the total length of $T_\eps$ is below $2\cdot n\cdot\eps\le 4\cdot L$.
By construction, 
$T_\eps$ is $\eps$-close to $K$ in the Hausdorff metric.

Clearly, there is a closed curve $\gamma_\eps$ whose image is $T_\eps$ and its length is twice the total length of $T_\eps$;
that is, 
\[\length\gamma_\eps\le 8\cdot L.\]

Passing to a subsequential limit of $\gamma_\eps$ as $\eps\to 0$,
we get the needed curve. \qeds

In terms of measure, the optimal bound is $2\cdot L$;
if in addition the diameter $D$ is known then it is $2\cdot L-D$.
The problem is due to 
Samuel Eilenberg 
and Orville Harrold 
\cite{eilenberg-harrold};
it also appears in the book of Kenneth Falconer \cite[see Exercise 3.5 in][]{falconer}.


\parbf{Shortcut.} 
Choose $\eps>0$.
Let us construct a set $X'\subset X$ 
 and a collection of paths $\alpha_0,\dots,\alpha_n$ such that 
\begin{enumerate}[(i)]
\item the total length of $\alpha_i\setminus X$ is at most $\eps$,
\item the set $X'$ is a union of a finite collection of closed connected sets $X_0,\dots,X_n$,
\item diameter of each $X_i$ is at most $\eps$,
\item $x\in\alpha_0$, $y\in\alpha_n$, and
\item the union $X'\cup \alpha_0\cup\dots\cup\alpha_n$ is connected
\end{enumerate}

Imagine that the construction is given.
Let us show that the statement can be proved by applying this construction recursively for a sequence of $\eps_n$ that converges to zero very fast.

Indeed we can apply the construction to each of the subsets $X_n$ and take as $X''$ the union of all closed subsets provided by the construction.
This way we obtain a nested sequence of closed sets $X\supset X'\supset X''\supset\dots$ which break into a finite union of closed connected subsets of arbitrary small diameter and a countable collection of arcs with total length at most $\eps_1+\eps_2+\dots$ 
outside of $X$.
Consider the Cantor set 
\[Y= X\cap  X'\cap X''\cap\dots\]
Note that there is a simple curve from $x$ to $y$ that runs in the constructed arcs and $Y$.
The total length of the constructed curves outside of $X$ can not exceed the sum
$\eps_1+\eps_2+\dots$;
whence the result.

It remains to do the construction.

Let us cover $X$ by a grid of $\tfrac\eps2$-squares $Q_1,\dots,Q_k$.
Denote by $\Delta$ the union of all the sides of the squares.

By the regularity of length, we may cover $\Delta\cap X$ by a finite collection of arcs with total length arbitrarily close to the length of $\Delta\cap X$.
Denote these arcs by $\alpha_0,\dots,\alpha_n$.
Without loss of generality, we may assume that $x\in\alpha_0$ and $y\in\alpha_n$.

Consider a finite graph with vertices labeled by $\alpha_0,\dots,\alpha_n$;
two vertices $\alpha_i$ and $\alpha_j$ are adjacent if there is a connected set $\Theta\subset  X\cap Q_k$ 
for some $k$ such that $\Theta$ intersects $\alpha_i$ and $\alpha_j$.
Note that the graph is connected. Therefore we may choose a path from $\alpha_0$ to $\alpha_n$ in the graph.

The path corresponds to a sequence of arcs $\alpha_i$ and a sequence of $\Theta$-sets.
The $\Theta$-sets that correspond to the edges in the path can be taken as $ X_i$ in the construction.
\qeds

This solution was found by Taras Banakh \cite{banakh};
it was used by Stephan Stadler and the author \cite{petrunin-stadler:revisited}.
The proof works only in dimension two and we are not aware of a generalization to higher dimensions. 
Namely, the following question is open:

\begin{pr}
Let $x,y$ be two points in a compact connected subset $X\subset \RR^3$. 
Is it always possible to connect $x$ and $y$ by a path $\alpha$ such that the complement $\alpha\setminus X$ has arbitrarily small length?
\end{pr}

\parbf{Straight set.}
A set $L$ in the plane will be called $r$-reachable
if it is formed by a collection of disjoint proper smooth curves 
with $r$-tubular neighborhood;
that is, the closest-point projection to $L$ is uniquely defined in the $r$-neighborhood of $L$.
Note that in this case, $L$ has curvature at most $\tfrac1r$ at any point.

Further, $\delta_i$ will denote a positive value that depends on $\delta$ such that $\delta_i\to0$ as $\delta\to 0$.
In fact, every statement below holds for $\delta_i=239\cdot\delta$ for any $i$.
We assume that $\delta$ and, therefore, each $\delta_i$ are small.

Fix an $\delta$-straight set $X$.
Note that for any $r>0$, there is a $r$-reachable set $L$ 
that is $\delta_1\cdot r$-close to $X$ in the sense of Hausdorff.
Indeed, using $\delta$-straightness one may approximate $X$ by a polygonal line that has sides about $r$, very obtuse angles, and that does not come close to a fixed point twice.
Smoothing its corners produces $L$.
\begin{figure}[ht!]
\centering
\includegraphics{mppics/pic-114}
\end{figure}

Let $L_n$ be a sequence of such approximations for $r_n=\tfrac1{2^n}$.
Note that $L_{n-1}$ lies in a $\delta_2\cdot r_n$-neighborhood of $L_{n}$.
In particular, the closest point projection $f_n\:L_{n-1}\to L_{n}$ is uniquely defined.
Moreover, it is straightforward to check that $f_n$ is $(1\mp \delta_3)$-bi-Lipschitz.

Given $x_1\in L_1$, consider the sequence $x_n\in L_n$ defined by $x_{n+1}\z=f_n(x_n)$.
From above, we have that 
\[
\begin{aligned}|x_{n+1}-x_n|&<\delta_4\cdot r_n,
&
|x_{n+1}-y_{n+1}|&\gtrless(1\mp\delta_4)\cdot |x_{n}-y_{n}|.
\end{aligned}
\eqno({*})
\]
In particular, for any $x_1\in L_1$, the sequence $(x_n)$ converges in itself.
Denote its limit by $F(x_1)$;
evidently, $F\:L_1\to X$ is onto.
It follows in particular, that any pair of points in $X$ on distance at most 1 lie in one connected component.
Since the $\delta$-straightness does not depend on the scale, it implies that $X$ is connected.

Note that $({*})$ implies the following two pairs of inequalities
\begin{align*}
|x_{n+1}-y_{n+1}|&\le (1+\delta_5)\cdot |x_{n}-y_{n}|,
\\
|x_{n+1}-y_{n+1}|&\le|x_{n}-y_{n}|+ \delta_5\cdot r_n,
\intertext{and}
|x_{n+1}-y_{n+1}|&\ge (1-\delta_5)\cdot |x_{n}-y_{n}|,
\\
|x_{n+1}-y_{n+1}|&\ge |x_{n}-y_{n}|- \delta_5\cdot r_n.
\end{align*}
To prove that $F$ is bi-Hölder,
apply recursively the first inequality in each pair until $|x_{n}-y_{n}|\z>r_n$ and after the second one. 
\qeds

This puzzle is a baby case of the so-called \emph{Reifenberg's lemma}.
It was introduced by Ernst Reifenberg \cite{reifenberg} and became a useful tool in metric geometry; in particular, it is used to study the limit spaces with lower Ricci curvature bound \cite{cheeger-colding-1, naber}.


\parbf{Typical convex curves.}
Denote by $\mathfrak{C}$ the space of all closed convex curves in the plane equipped with the Hausdorff metric.
Recall that $\mathfrak{C}$ is locally compact.
In particular, by the Baire theorem, a countable intersection of everywhere dense open sets is everywhere dense.

Note that if a curve $\gamma\in\mathfrak{C}$ has non-zero second derivative at a point~$p$,
then $\gamma$ lies between two nested circles tangent to each other at $p$.

Fix these two circles.
It is straightforward to check that there is $\eps>0$ such that 
the Hausdorff distance from any convex curve $\gamma$ squeezed between the circles 
to any convex $n$-gon is at least $\frac{\eps}{n^{100}}$.

\begin{wrapfigure}{o}{41 mm}
\vskip-0mm
\centering
\includegraphics{mppics/pic-116}
\end{wrapfigure}

Choose a countable dense set of convex polygons $\mathfrak{p}_1,\mathfrak{p}_2,\dots$ in $\mathfrak{C}$.
Denote by $n_i$ the number of sides in $\mathfrak{p}_i$.
For any positive integer $k$,
consider the set $\Omega_k\subset\mathfrak{C}$ defined by 
\[\Omega_k
=
\set{\xi\in \mathfrak{C}}{\min_i\{|\xi-\mathfrak{p}_i|_{H}\}<\tfrac1{k\cdot n_i^{100}}},\]
where $|{*}-{*}|_H$ denotes the Hausdorff distance 

From the above, we get that $\gamma\notin\Omega_k$ for some $k$. 

Note that $\Omega_k$ is open and everywhere dense in $\mathfrak{C}$.
Therefore 
\[\Omega=\bigcap_k\Omega_k\]
is a G-delta dense set.
Hence the statement follows.\qeds

This problem states that typical convex curves have an unexpected property.
This is common --- it is hard to see the typical objects and these objects often have surprising properties.

For example, it was proved by
Bernd Kirchheim, 
Emanuele Spadaro, 
and 
L{\'a}szl{\'o} Sz{\'e}kelyhidi,
that \textit{typical 1-Lipschitz maps from the plane to itself preserve the length of all curves} \cite{KSS}.
The same way, one could show that \textit{the boundaries of typical open sets in the plane contain no nontrivial curves}, 
but the construction of a concrete example is not trivial
[see ``Crooked circle'', page \pageref{Crooked circle}].
More problems of that type are surveyed by Tudor Zamfirescu \cite{zamfirescu}.

\csname @openrightfalse\endcsname
\chapter{Surfaces}

We assume that the reader is familiar with smooth surfaces and the related definitions
including intrinsic metric, 
geodesics,
convex and saddle surfaces
as well as different types of curvature.
An introductory course in differential geometry should cover all necessary background material 
[see \ncite{toponogov-curves-and-surfaces}, \ncite{petrunin-zamora}, or \S28--29 in \ncite{hilbert-cohn-vossen}].

{

\begin{wrapfigure}{r}{27 mm}
\vskip4mm
\centering
\includegraphics{mppics/pic-202}
\end{wrapfigure}

\subsection*{Convex hat}
\label{Convex hat}

\begin{pr}
Suppose that a plane $\Pi$ cuts from a smooth closed convex surface $\Sigma$ a disk $\Delta$.
Assume that the reflection of $\Delta$ with respect to $\Pi$ lies inside of $\Sigma$.
Show that $\Delta$ is \index{convex set}\emph{convex} with respect to the intrinsic metric  of $\Sigma$;
that is, 
if both ends of a minimizing geodesic in $\Sigma$ 
lie in $\Delta$,
then the entire geodesic lies in $\Delta$.
\end{pr}

}

\parit{Semisolution.}
Assume the contrary,
then there is a minimizing geodesic $\gamma\not\subset\Delta$ with ends $p$ and $q$ in $\Delta$.

Without loss of generality, we may assume that only one arc of $\gamma$ lies outside of $\Delta$.
Reflection of this arc  with respect to $\Pi$ together with the remaining part of $\gamma$ forms another curve $\hat\gamma$ from $p$ to $q$;
it runs partly along $\Sigma$ 
and partly outside $\Sigma$,
but does not get inside $\Sigma$.
Note that
\[\length\hat\gamma=\length\gamma.\]

Denote by $\bar\gamma$ the closest point projection of $\hat\gamma$ on $\Sigma$.
Since $\Sigma$ is convex, the closest point projection decreases the length.
Therefore 
the curve $\bar\gamma$ lies in $\Sigma$, 
it has the same ends as $\gamma$,
and
\[\length\bar\gamma<\length\gamma.\]
This means that $\gamma$ is not length-minimizing 
--- a contradiction.\qeds


{

\begin{wrapfigure}{r}{39 mm}
\vskip0mm
\centering
\includegraphics{mppics/pic-204}
\end{wrapfigure}

\subsection*{Involute of geodesic}
\label{Involute of  geodesic}

\begin{pr}
Let $\Sigma$ be a smooth closed strictly convex surface 
in $\RR^3$ 
and $\gamma\:[0,\ell]\z\to \Sigma$ a unit-speed minimizing geodesic.
Set $p\z=\gamma(0)$, $q=\gamma(\ell)$, and 
$$p_t=\gamma(t)-t\cdot\gamma'(t),$$ 
where $\gamma'(t)$ denotes the velocity vector of $\gamma$ at~$t$.
\end{pr}

}
\vskip-\medskipamount
\textit{Show that for any $t\in (0,\ell)$,
one {}\emph{cannot see}  $q$ from $p_t$;
that is, the line segment $[p_tq]$ intersects $\Sigma$ at a point distinct from $q$.}

\subsection*{Simple geodesic}
\label{Simple geodesic}

\begin{pr}
Let $\Sigma$ be a complete unbounded convex surface in $\mathbb R^3$.
Show that there is a two-sided infinite geodesic in $\Sigma$ with no self-intersections.
\end{pr}

Let us review a couple of statements 
about Gauss curvature which might help to solve the problem \cite[see \S28 in][for more details]{hilbert-cohn-vossen}.

If $\Sigma$ is a convex surface in $\RR^3$ then its Gauss curvature is nonnegative.

Assume that a simply-connected region $\Omega$ in the surface $\Sigma$ is bounded by a closed broken geodesic $\gamma$.
Denote by $\kappa(\Omega)$ the integral of the Gauss curvature over $\Omega$.

For any point $p\in\Sigma$ consider the outer unit normal vector $n(p)\z\in\mathbb{S}^2$.
Then 
\begin{align*}
\kappa(\Omega)&=\area[n(\Omega)]
\intertext{and by the Gauss--Bonnet formula}
\kappa(\Omega)&=2\cdot\pi-\sigma(\gamma),
\end{align*}
where $\sigma(\gamma)$ denotes the sum of the signed exterior angles of $\gamma$.
In particular,  $|\sigma(\gamma)|\le2\cdot\pi$.

\subsection*{Geodesics for birds}
\label{liberman}

The \index{total curvature}\emph{total curvature} of a space curve $\gamma$ is defined as the integral of its curvature.
That is, if a curve $\gamma\:[a,b]\to\RR^3$ has unit speed parametrization, 
then its total curvature equals 
\[\int\limits_a^b|\gamma''(t)|\cdot dt,\]
the vector $\gamma''(t)$ is called the \index{curvature vector}\emph{curvature vector} and its length $|\gamma''(t)|$ is the \index{curvature}\emph{curvature} of $\gamma$ at time $t$.
The above definition makes sense for $C^{1,1}$ smooth curves,
that is, in the case when $\gamma'(t)$ is locally Lipschitz;
in this case the curvature $|\gamma''(t)|$ is defined almost everywhere.

The \index{geodesic}\emph{geodesics} in the following problem are defined as the curves locally minimizing the length;
that is, any sufficiently short arc of the curve containing a given value of the parameter is length-minimizing.

\begin{pr}
Let $f\:\RR^2\to\RR$ be a smooth $\ell$-Lipschitz function.
Let $W\subset \RR^3$ be the epigraph of $f$;
that is,
$$W=\set{(x,y,z)\in\RR^3}{z\ge f(x,y)}.$$
Equip $W$ with the induced intrinsic metric.

Show that any geodesic in $W$ 
 has  total curvature at most $2\cdot\ell$. 
\end{pr}

Actually, geodesics in $W$ are $C^{1,1}$-smooth;
in particular, the formula for the total curvature mentioned above makes sense.
This is an easy exercise in real analysis which can be also taken for granted.

\subsection*{Immersed surface}
\label{Immersed surface}

\begin{pr}
Let $\Sigma$ be a smooth connected immersed surface in $\RR^3$ with strictly positive Gauss curvature and nonempty boundary $\partial\Sigma$.
Assume that the boundary $\partial\Sigma$ lies in a plane $\Pi$
and $\Sigma$ lies entirely on one side of~$\Pi$.
Prove that $\Sigma$ is an embedded disk.
\end{pr}

\subsection*{Periodic asymptote}
\label{Asymptotic geodesic}

\begin{pr}
Let $\Sigma$ be a closed smooth surface with non-positive curvature at every point
and $\gamma$ a geodesic in $\Sigma$.
Assume that $\gamma$ is not periodic
and the curvature of $\Sigma$ vanishes at every point of $\gamma$.
Show that $\gamma$ does not have a periodic asymptote;
that is, there is no periodic geodesic $\delta$ such that the distance from $\gamma(t)$ to $\delta$  converges to $0$ as $t\to\infty$. 
\end{pr}

\subsection*{Saddle surface}
\label{Saddle surface}

Recall that a smooth surface $\Sigma$ in $\RR^3$
is \index{saddle surface}\emph{saddle} at a point $p$ if its principal curvatures at $p$ have opposite signs. 
We say that $\Sigma$ is {}\emph{saddle} if it is saddle at all points.

\begin{pr}
Let $\Sigma$ be a saddle surface in $\RR^3$
homeomorphic to a disk.
Assume that the orthogonal projection to the $(x,y)$-plane
maps the boundary of $\Sigma$
injectively to a convex closed curve.
Show that this projection is injective on~$\Sigma$.

In particular, $\Sigma$ is the graph $z=f(x,y)$ of a function $f$ defined on a convex domain in the $(x,y)$-plane.
\end{pr}

\subsection*{Asymptotic line}
\label{asymptotic-line}

The saddle surfaces are defined in the previous problem.

Recall that a smooth curve $\gamma$ on a smooth surface $\Sigma\subset \RR^3$ is called \index{asymptotic line}\emph{asymptotic} if $\gamma''(t)$ is tangent to the surface at $\gamma(t)$ for any $t$.

\begin{pr}
Let $\Sigma\subset \RR^3$ be the graph $z\z=f(x,y)$
of a smooth function $f$ 
and $\gamma$ a closed smooth asymptotic line in $\Sigma$.
Assume that $\Sigma$ is saddle in a neighborhood of $\gamma$.
Show that the projection of $\gamma$ to the $(x, y)$-plane cannot be star-shaped;
that is, there is no point $p$ in the plane such that each half-line from $p$ intersects the projection at exactly one point.
\end{pr}

\subsection*{Minimal surface}
\label{min-surf}

Recall that a smooth surface in $\RR^3$ is called \index{minimal surface}\emph{minimal} if its mean curvature vanishes at all points.
The \index{mean curvature}\emph{mean curvature} at each point is defined as the sum of the principal curvatures at that point.

\begin{pr}
Let $\Sigma$ be a minimal surface in $\RR^3$ having its boundary on a unit sphere.
Assume that $\Sigma$ passes thru the center of the sphere.
Show that the area of $\Sigma$ is at least $\pi$.
\end{pr}

{

\begin{wrapfigure}[3]{r}{43 mm}
\vskip-4mm
\centering
\includegraphics{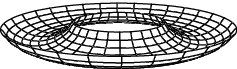}
\end{wrapfigure}

\subsection*{Round gutter\hard}
\label{half-torus}

A round gutter is the surface shown on the picture.

}

More precisely, consider the torus $T$; a surface generated by revolving a circle in $\RR^3$ around an axis coplanar with the circle.
Let $\gamma\z\subset T$ be one of the circles in $T$ that locally separates positive and negative curvature on $T$;
a plane containing $\gamma$ is tangent to $T$ at all points of $\gamma$.
Then a neighborhood of $\gamma$ in $T$ is called 
a \index{gutter}\emph{round gutter}
and the circle $\gamma$ is called its {}\emph{main latitude}.

\begin{pr}
Let $\Omega\subset \RR^3$ be a round gutter with the main latitude $\gamma$. 
Assume that $\iota\:\Omega\z\to\RR^3$ 
is a smooth length-preserving embedding that is sufficiently close to the identity.
Show that $\gamma$ and $\iota(\gamma)$ are congruent;
that is, there is an isometric motion of $\RR^3$ sending $\gamma$ to $\iota(\gamma)$
\end{pr}

\subsection*{Non-contractible geodesics}
\label{torus}

\begin{pr}
Give an example of a non-flat metric 
on the $2$-torus such that no closed geodesic is contractible.
\end{pr}

\subsection*{Two disks}
\label{Two disks}

\begin{pr}
Let $\Sigma_1$ and $\Sigma_2$ be two smoothly embedded open disks in $\mathbb R^3$ 
that have a common closed smooth curve $\gamma$.
Show that there is a pair of points  $p_1\in \Sigma_1$ and $p_2\in \Sigma_2$ with parallel tangent planes.
\end{pr}

\subsection*{Second derivative bounds first\easy}
\label{Second derivative bounds first}

\begin{pr}
Let $f$ be a smooth function on the Lobachevsky plane with bounded Hessian.
Show that $f$ is Lipschitz.
\end{pr}

\section*{Semisolutions}
\parbf{Involute of  geodesic.}
Let $W$ be the closed unbounded set formed by $\Sigma$ and its exterior points.
Choose $t\in (0,\ell)$;
denote by $\gamma_t$ the concatenation of the line segment $[p_t\gamma(t)]$ and the arc $\gamma|_{[t,\ell]}$.
The key step is to show the following:

\begin{cl}{$({*})$} The curve $\gamma_t$ is a minimizing geodesic in the intrinsic metric induced on $W$.
\end{cl}

Try to prove it before reading further.

\medskip

Let $\Pi_t$ be the tangent plane to $\Sigma$ at $\gamma(t)$.
Consider the curve $\alpha(t)=p_t$.
Note that  
$\alpha(t)\z\in\Pi_t$,
$\alpha'(t)\perp\Pi_t$,
and $\alpha'(t)$ points to the side of $\Pi_t$ opposite to $\Sigma$.

It follows that for any $x\in\Sigma$ the function  
\[t\mapsto |x - p_t|
\quad\text{and, therefore,}\quad
t\mapsto |x - p_t|_W\] are non-decreasing;
here $|x - p_t|_W$ stands for the intrinsic distance from $x$ to $p_t$ in $W$.

\begin{wrapfigure}{r}{43 mm}
\vskip-6mm
\centering
\includegraphics{mppics/pic-205}
\end{wrapfigure}

On the other hand, by construction 
\[|q - p_t|_W\le |q - p|_\Sigma;\] 
therefore, from the above 
\[|q - p_t|_W= |q - p|_\Sigma\]
for any $t$.
Hence $(*)$ follows.

Now assume that $q$ is visible from $p_t$ for some $t$;
that is, the line segment $[qp_t]$ intersects $\Sigma$ only at $q$.
From the above, 
$\gamma_t$  coincides with the line segment $[qp_t]$.
On the other hand, $\gamma_t$ contains $\gamma(t)\in\Sigma$ --- a contradiction.\qeds

This problem is based on an observation used by Anatoliy Milka in the proof of the following generalization of the comparison theorem for convex surfaces \cite[Theorem 2]{milka-geod}.

\begin{pr}
Let $W$ be the closed unbounded set formed by a (not necessarily smooth) closed convex surface $\Sigma$ and its exterior points.
Suppose $\gamma_1\:[0,\ell_1]\to \Sigma$ and $\gamma_2\:[0,\ell_2]\to \Sigma$ are unit-speed minimizing geodesics that start at the same point $p$.
Set $x_i=\gamma_i(\ell_i)$ and $y_i=p+\ell_i\cdot\gamma_i'(0)$.
Show that 
\[|x_1-x_2|_W\le |y_1-y_2|_W,\]
where $|\  - \ |_W$ stands for the intrinsic distance in $W$.
\end{pr}


\begin{wrapfigure}{r}{43 mm}
\vskip-0mm
\centering
\includegraphics{mppics/pic-206}
\end{wrapfigure}

\parbf{Simple geodesic.} 
Look at two combinatorial types of a self-intersection shown in the diagram.
One of them can and the other cannot appear as self-intersections of a geodesic on an unbounded convex surface.
Try to determine which is which before reading further.

\medskip

Let $\gamma$ be a two-sided infinite geodesic in $\Sigma$.
The following is the key statement in the proof.

\begin{cl}{$({*})$}
The geodesic $\gamma$ contains at most one simple loop.
\end{cl}

To prove $({*})$, we use the following observation.

\begin{cl}{$({*}{*})$}
The integral curvature $\omega$ of $\Sigma$ cannot exceed $2\cdot\pi$.
\end{cl}

Indeed, since $\Sigma$ is unbounded and convex,
it surrounds a half-line.
Consider a coordinate system with this half-line as the positive half of its $z$-axis. 
In these coordinates, the surface $\Sigma$ is described as a graph $z=f(x,y)$ of a convex function $f$.
In particular, all outer normal vectors to $\Sigma$ have a negative $z$-coordinate;
in other words, they point to the southern hemisphere.
Therefore the area of the spherical image of $\Sigma$ is at most $2\cdot\pi$.
The area of this image is the integral of the Gauss curvature over $\Sigma$. 
Hence $({*}{*})$ follows.

From the Gauss--Bonnet formula, we get the following conclusion.
If $\phi$ is the angle at the base of a simple geodesic loop then the integral curvature surrounded by the loop equals $\pi+\phi$. 
In particular, $({*}{*})$ implies that $\phi\le\pi$; in other words, there are no \emph{concave loops}.

Now assume that $({*})$ does not hold, so that a geodesic has two simple loops.
Note that the disks bounded by the loops  have to overlap,
otherwise the curvature of $\Sigma$ would exceed $2\cdot\pi$.
But if they overlap, then it is easy to show that the curve also contains a concave loop, 
which contradicts the above observation.%
\footnote{This observation implies that the right picture on the above diagram cannot be realized by a geodesic.}

If a geodesic $\gamma$ has a self-intersection,
then it contains a simple loop.
From $({*})$, there is only one such loop;
it cuts a disk from $\Sigma$ 
and goes around it either clockwise or counterclockwise.
This way we divide all the self-intersecting geodesics 
into two sets which we will call {}\emph{clockwise} and {}\emph{counterclockwise}.

Note that the geodesic $t\mapsto \gamma(t)$ is clockwise 
if and only if 
$t\z\mapsto \gamma(-t)$
is counterclockwise.
The sets of clockwise and counterclockwise are open and the space of all geodesics is connected. 
It follows that there are geodesics that aren't clockwise nor counterclockwise.
Those geodesics have no self-intersections.\qeds

Note that the proof implies that a two-sided infinite geodesic can be found among geodesics containing a given point in $\Sigma$.

The problem is due to Stephan Cohn-Vossen \cite[Satz 9 in][]{convossen};
generalizations were obtained  by 
Vladimir Streltsov and Alexandr Alexandrov 
\cite{streltsov-alexandrov} 
and 
by Victor Bangert \cite{bangert}.

The following problem is of the same style \cite{petrunin-self-crossing-geodesics}.

\begin{pr}
Let $\gamma$ be a closed geodesic on a closed surface of positive curvature $\Sigma$.
\begin{figure}[ht!]
\vskip-0mm
\centering
\includegraphics{mppics/pic-207}
\end{figure}
Show that $\gamma$ cannot look like one of the curves on the diagram.
In other words, $\gamma$ cannot subdivide $\Sigma$ into (1) one hexagon, one triangle, and three monogons, and (2) one pentagon, one quadrangle, and three monogons.
\end{pr}

\parbf{Geodesics for birds.}
Choose a unit-speed geodesic in $W$, say
\[\gamma\:t\mapsto(x(t),y(t),z(t)).\]
We can assume that $\gamma$ is defined on a closed interval $[a,b]$.
The key step is to show the following:

\begin{cl}{$({*})$} 
The function $t\mapsto z$ is concave.
\end{cl}

Parametrize the plane curve $t\mapsto (x(t),y(t))$ by the arc length $s$
and reparametrize $\gamma$ by $s$.

\begin{wrapfigure}{o}{35 mm}
\vskip-0mm
\centering
\includegraphics{mppics/pic-208}
\end{wrapfigure}

Note that the function $s\mapsto z$ is concave.
Indeed, suppose $s\mapsto z$ is not concave around $s_0$.
Then one could shorten $\gamma$ by increasing its $z$ component in a small interval around $s_0$ while keeping its endpoints fixed.
After the deformation, the curve still lies in $W$.
The latter contradicts that $\gamma$ is locally length-minimizing.

Finally, note that the concavity of $s\mapsto z$ is equivalent to the concavity of $t\mapsto z$.
Hence $({*})$ follows.

Since $f$ is smooth, 
the curve $\gamma(t)$ is $C^{1,1}$; 
that is, its first derivative $\gamma'(t)$ is a well-defined Lipschitz function.
It follows that its second derivative $\gamma''(t)$ is defined almost everywhere.

Since $z(t)$ is concave, we have $z''(t)\le 0$.
Since $f$ is $\ell$-Lipschitz, $z(t)$ is $\tfrac{\ell}{\sqrt{1+\ell^2}}$-Lipschitz.
It follows that 
\[\int\limits_a^b |z''(t)|\cdot dt\le 2\cdot\tfrac{\ell}{\sqrt{1+\ell^2}}.\]

The curvature vector $\gamma''(t)$ is perpendicular to the surface.
Since the surface has slope at most $\ell$,
we get 
\[|\gamma''(t)|\le |z''(t)|\cdot\sqrt{1+\ell^2}.\]
Hence 
\[\int\limits_a^b |\gamma''(t)|\cdot dt\le 2\cdot\ell.\qedsin\]
\medskip

The statement holds for general $\ell$-Lischitz functions,
not necessarily smooth.
The given bound is optimal, the equality is attained by a two-side infinite geodesic on the graph of  
\[f(x,y)=-\ell\cdot\sqrt{x^2+y^2}.\]

The problem is due to David Berg \cite{berg},
the same bound for convex $\ell$-Lipschitz surfaces was proved earlier by Vladimir Usov \cite{usov}.
The observation $({*})$
is called \index{Liberman’s lemma}\emph{Liberman’s lemma}; 
it was used earlier 
to bound the total curvature
of a geodesic on a convex surface \cite{liberman}.\footnote{It is a part of the thesis of Joseph Liberman, defended a couple of months before he died in WWII.}
This lemma is often useful when working with geodesics on general convex surfaces.

\parbf{Immersed surface.}
Let $\ell$ be a linear function that vanishes on $\Pi$ 
and is positive on $\Sigma$. 
We will apply a Morse-type argument for the restriction of $\ell$ to $\Sigma$.

\medskip

Let $z_0$ be a maximum of $\ell$ on $\Sigma$;
set $s_0=\ell(z_0)$.
Given $s<s_0$, denote by $\Sigma_s$ the connected component of $z_0$ in $\Sigma\cap\ell^{-1}([s,s_0])$.
Note that 
\begin{itemize}
\item $\Sigma_s$ is an embedded disk, and
\item $\partial\Sigma_s$ is a convex plane curve
\end{itemize}
for all $s$ sufficiently close to $s_0$.

Consider the set $A\subset [0,s_0)$ such that for any $a\in A$ these two conditions hold for any $s\in [a,s_0)$.
Observe that $A$ is open and closed in $[0,s_0)$.
Whence $A=[0,s_0)$; in particular, these conditions hold for $s=0$.

Since $\Sigma$ is connected, $\Sigma_0=\Sigma$.
Hence the result follows.\qeds

This problem is discussed in the lectures of Mikhael Gromov \cite[see \S$\tfrac12$~in][]{gromov-SGMC}.

\parbf{Periodic asymptote.}
Arguing by contradiction, assume that there is a geodesic $\gamma$ on the surface $\Sigma$ with a periodic asymptote $\delta$. 

Passing to a finite cover of $\Sigma$, we can ensure that the asymptote has no self-intersections.
In this case, 
the restriction $\gamma|_{[a,\infty)}$  
has no self-intersections 
if $a$ is sufficiently large.

Cut $\Sigma$ along $\gamma([a,\infty))$ and then cut from the obtained surface an infinite triangle $\triangle$. 
The triangle has two sides formed by both sides of cuts along $\gamma$;
let us denote these sides of $\triangle$ by $\gamma_-$ and $\gamma_+$.
Note that 
\[\area\triangle<\area \Sigma<\infty,
\leqno(*)\]
and both sides $\gamma_\pm$ 
are infinite minimizing geodesics in $\triangle$.

Consider the Busemann function $f$ for $\gamma_+$ [defined on page~\pageref{page:Busemann function}];
denote by $\ell(t)$ the length of the level curve $f^{-1}(t)$.
Let $-\kappa(t)$  be the total curvature of the sup-level set $f^{-1}([t,\infty))$.  
From the Gauss--Bonnet formula,
\[\ell'(t)=\kappa(t).\leqno({*}{*})\]

The level curve $f^{-1}(t)$ can be parametrized by a unit-speed curve, say $\theta_t\:[0,\ell(t)]\to \triangle$.
By the coarea formula we have
\[\kappa'(t)
=
-\int\limits_0^{\ell(t)} K_{\theta_t(\tau)}\cdot d\tau,
\]
where $K_x$ denotes the Gauss curvature of $\Sigma$ at the point $x$.
Since $K_{\theta_t(0)}\z=K_{\theta_t(\ell_t)}=0$ and the surface is smooth,
there is a constant $C$ such that $|K_{\theta_t(\tau)}|\le C\cdot \ell(t)^2$ for all $t$, $\tau$.
Therefore
\[\kappa'(t)\le C\cdot \ell(t)^3 \leqno(\asterism)\]

Together, $({*}{*})$ and $(\asterism)$ imply that there is $\eps>0$ such that
\[\ell(t)\ge \frac\eps{t-a}\]
for sufficiently large $t$.
By the coarea formula we get 
\[\area\triangle=\int\limits_a^\infty\ell(t)\cdot dt=\infty;\]
the latter contradicts $(*)$.\qeds

I learned the problem from 
Dmitri Burago 
and Sergei Ivanov, 
it originated from a discussion with
Keith Burns, 
Michael Brin, 
and Yakov Pesin.

Here is a motivation.
Let $\Sigma$ be a closed surface with non-positive curvature that is not flat.
The space $\Gamma$ of all unit-speed geodesics $\gamma\:\RR\z\to\Sigma$ can be identified with the unit tangent bundle $\UU\Sigma$. 
In particular $\Gamma$ comes with a natural choice of measure.
Denote by $\Gamma_0\subset \Gamma$ the set of geodesics that run in the set of zero curvature all the time.
It is expected that $\Gamma_0$ has a vanishing measure.
In all known examples $\Gamma_0$ contains only periodic geodesics in only finitely many homotopy classes \cite[see also][]{hertz}.

\parbf{Saddle surface.}
Denote by $\Sigma^\circ$ the interior of $\Sigma$.
Choose a plane $\Pi$. 
Note that the intersection $\Pi\cap \Sigma^\circ$ 
locally looks either like a curve or like two curves intersecting transversally;
in the latter case, $\Pi$ is tangent to $\Sigma^\circ$ at the intersection point.

Further, note that $\Pi\cap \Sigma^\circ$ has no cycle.
Otherwise, $\Sigma$ would fail to be saddle at the point of the disk surrounded by that cycle maximizing the distance to $\Pi$.

If $\Sigma$ is not a graph then there is a point $p\in\Sigma$ with a vertical tangent plane;
denote this plane by $\Pi$.
The intersection $\Pi\cap\Sigma$ has a cross-point at~$p$.

Since the boundary of $\Sigma$ projects injectively to a closed convex curve in $(x,y)$-plane,
the intersection of $\Pi\cap\partial \Sigma$ has at most 2 points --- these are the only endpoints of $\Pi\cap\Sigma$.

It follows that the connected component of $p$ in $\Pi\cap\Sigma$ is a tree 
with a vertex of degree 4 at $p$ and at most two end-points --- a contradiction.\qeds

The described idea can be used to prove the result of Richard Schoen and Shing-Tung  Yau \cite{schoen-yau-2D} which gives a sufficient condition for a harmonic map between surfaces to be a diffeomorphism.
Unlike the original proof, it requires no calculations.

The proof above is based on the observation 
that for any saddle surface $\Sigma$ and plane $\Pi$,
each connected component of $\Sigma\setminus \Pi$ is either unbounded or intersects the boundary curve.
This observation plays a central role in the proof of Sergei Bernstein \cite{bernshtein}
of the following problem:

\begin{pr}
Show that a smooth bounded function $f\:\RR^2\to\RR$ cannot have a strictly saddle graph.
\end{pr}

One could go further and define a \emph{generalized saddle surface} as an arbitrary (non-necessarily smooth) surface satisfying the observation above.
The geometry of these surfaces is far from being understood,
Samuil Shefel has a number of beautiful results about them, 
\cite[see][and the references therein]{shefel, AKP-invitation}.
The statement of the problem holds for these generalized saddle surfaces, but
the proof is tricky \cite{petrunin-stadler}.

\parbf{Asymptotic line.}
Denote by $\Pi_t$ the tangent plane to $\Sigma$ at $\gamma(t)$ and by $\ell_t$ the tangent line of $\gamma$ at time $t$.

Since $\gamma$ is asymptotic, the plane $\Pi_t$ rotates around $\ell_t$ as $t$ changes.
Since $\Sigma$ is saddle, the speed of rotation cannot vanish.%
\footnote{By the Beltrami--Enneper theorem, if $\gamma$ has unit speed, then the speed of rotation is $\pm\sqrt{-K}$, where $K$ is the Gauss curvature which cannot vanish on a saddle surface.}

Note that $\Pi_t$ is a graph of a linear function, say $h_t$, defined on the $(x, y)$-plane.
Denote by $\bar\ell_t$ the projection of $\ell_t$ to the $(x, y)$-plane.
The described rotation of $\Pi_t$ can be expressed algebraically:
the derivative $\tfrac{d}{dt}h_t(w)$ vanishes at the point $w$ if and only if $w\in \bar\ell_t$ 
and the derivative changes sign if $w$ changes the side of $\bar\ell_t$.

Denote by $\bar\gamma$ the projection of $\gamma$ to the $(x, y)$-plane.
If $\bar\gamma$ is star-shaped with respect to a point $w$, then $w$ cannot cross $\bar\ell_t$.
Therefore the function $t\mapsto h_t(w)$ is monotone on $\SSS^1$.
Observing that this function cannot be constant, we arrive at a contradiction.\qeds

This is a stripped version of the result of Galina Kovaleva \cite{kovaleva}, which was rediscovered by Dmitri Panov \cite{panov-curves,arnold}.



\parbf{Minimal surface.}
Without loss of generality, we may assume that the sphere is centered at the origin of $\RR^3$.

Let $h$ be the restriction of the function $x\mapsto \tfrac12\cdot|x|^2$ to the surface $\Sigma$.
Direct calculations show that $\Delta_\Sigma h =  2$;
here $\Delta_\Sigma$ denoted Laplacian on $\Sigma$.
Applying the divergence theorem for the gradient $\nabla_\Sigma h$
in $\Sigma_r\z=\Sigma\cap B(0,r)$, we get
\[2\cdot \area\Sigma_r\le r\cdot\length [\partial\Sigma_r].\]

Set $a(r)=\area\Sigma_r$.
By the coarea formula, $a'(r)\ge\length [\partial\Sigma_r]$ for almost all $r$.
Therefore the function
\[f\:r\mapsto \frac{\area\Sigma_r}{r^2}
\]
is non-decreasing in the interval $(0,1)$.

Since $f(r)\to \pi$ as $r\to0$, the result follows.\qeds

We described a partial case of the so-called \index{monotonicity formula}\emph{monotonicity formula} for minimal surfaces.

The same argument shows that if $0$ is a double point
of $\Sigma$ then $\area\Sigma\z\ge 2\cdot \pi$.
This observation was used to prove 
that the minimal disk bounded by a simple closed curve with total curvature $\le 4\cdot\pi$ 
is necessarily embedded.
It was proved by 
Tobias Ekholm, 
Brian White,
and Daniel Wienholtz
\cite{EWW};
an amusing variation of this proof
was obtained by 
Stephan Stadler \cite{stadler-FM}.
This result also implies that any embedded circle of total curvature at most $4\cdot\pi$ is unknot.
The latter was proved independently by Istv{\'a}n F{\'a}ry \cite{fary-knot} and  John Milnor \cite{milnor}.

Note that if we assume in addition that the surface is a disk,
then the statement holds for any saddle surface. 
Indeed, denote by $S_r$ the sphere of radius $r$ concentrical with the unit sphere. 
Then according to the problem ``A curve on a sphere'' [page \pageref{A curve in a sphere}], 
\[\length[\partial\Sigma_r]\ge 2\cdot\pi\cdot r.\]
Then by the coarea formula, we get $\area\Sigma\ge \pi$.

On the other hand, there are saddle surfaces homeomorphic to the cylinder
having an arbitrarily small area in the ball. 

If $\Sigma$ does not pass thru the center 
and we only know the distance, say $r$, 
from the center to $\Sigma$,
then the optimal bound is $\pi\cdot(1-r^2)$.
This question was open for about 40 years and proved by Simon Brendle and Pei-Ken Hung \cite{brende-hung};
their proof is based on a similar idea and is quite elementary.
Earlier Herbert Alexander, 
David Hoffman,
and Robert Osserman 
proved it for two cases: (1) for disks and (2) for arbitrary area minimizing surfaces, any dimension and codimension
 \cite{alexander-osserman,alexander-hoffman-osserman}.

\parbf{Round gutter.}
Without loss of generality, we can assume that the length of $\gamma$ is $2{\cdot}\pi$ and its intrinsic curvature is $1$ at all points.

Let $K$ be the convex hull of $\hat\Omega=\iota(\Omega)$.
Part of $\hat\Omega$ touches the boundary of $K$ and the rest lies in the interior of $K$. 
Denote by $\hat\gamma$ the curve in $\hat\Omega$ dividing these two parts.

First note that the Gauss curvature of $\hat\Omega$ should vanish at the points of $\hat\gamma$;
in other words, $\hat\gamma=\iota(\gamma)$.
Indeed, since $\hat\gamma$ lies on the convex part, 
the Gauss curvature at the points of $\hat\gamma$ should be non-negative. 
On the other hand, $\hat\gamma$ bounds a flat disk in $\partial K$;
therefore its integral intrinsic curvature should be $2{\cdot}\pi$.
If the Gauss curvature is positive at a point of $\hat\gamma$, 
then by the Gauss--Bonnet formula, the total intrinsic curvature of $\hat\gamma$ should be smaller than $2{\cdot}\pi$ --- a contradiction.

On the other hand, $\hat\gamma$ is an asymptotic line.
Indeed, if the direction of $\hat\gamma$ is not asymptotic at $t_0$,
then $\hat\gamma(t_0 \pm\eps)$ lies the interior of $K$ for some small $\eps>0$ --- a contradiction.

Therefore, as the space curve,
$\hat\gamma$ has to be a closed curve with constant curvature $1$.
Any such curve is congruent to a unit circle.\qeds

It is not known whether $\hat\Omega$ is congruent to $\Omega$ or not.

The solution presented above is based on my answer 
to the question of Joseph O'Rourke \cite{rourke}.
Here are some related statements.
\begin{itemize}
\item A gutter is second-order rigid;
this was proved by Eduard Rembs
[see \ncite{rembs} and also page 135 in \ncite{efimov}].
\item Any second-order rigid surface does not admit analytic deformation 
[proved by Nikolay Efimov; page 121 in \ncite{efimov}]
and for the surfaces of revolution, the assumption of analyticity can be removed 
\cite[proved by Idzhad Sabitov, see][]{sabitov}.
\end{itemize}

\parbf{Non-contractible geodesics.}
A torus of revolution is an example.

Indeed, it has a family of {}\emph{meridians} --- a family of circles that form closed geodesics.
Note that a geodesic on the torus is either a meridian
or it intersects meridians transversally.
In the latter case, all the meridians are crossed by the geodesic in the same direction.

Note that a contractible curve has to cross each meridian an equal number of times in both directions.
Therefore no geodesic of the torus is contractible.\qeds

I learned this problem 
from the book of Mikhael Gromov \cite{gromov-MetStr},
where it is attributed to Y. Colin de Verdi\`ere.
The same argument can be used to show that a torus with a geodesic foliation has no contractible closed geodesics.
I do not know other examples of that type \cite{petrunin-torus}; namely, the following question is open:

\begin{pr}
Are there examples of Riemannian metrics on the 2-torus without closed null-homotopic geodesics and without a geodesic foliation?
\end{pr}

\parbf{Two disks.}
Choose a continuous map $h\:\Sigma_1\to \Sigma_2$
that is the identity on~$\gamma$.
Let us prove that for some $p_1\in \Sigma_1$ and $p_2=h(p_1)\in \Sigma_2$,
the tangent planes $\T_{p_1} \Sigma_1$ and  $\T_{p_2} \Sigma_2$ are parallel;
this fact is stronger than the one required.

\medskip

Arguing by contradiction,
assume that such a point does not exist.
Then for each $p\in\Sigma_1$
there is a unique line $\ell_p\ni p$ 
 parallel to both $\T_{p} \Sigma_1$ and $\T_{h(p)} \Sigma_2$.

Note that the lines $\ell_p$ form a tangent line distribution over $\Sigma_1$
and $\ell_p$ is tangent to $\gamma$ at all $p\in\gamma$.

Let $\Delta$ be the disk in $\Sigma_1$ bounded by $\gamma$.
Consider the doubling of $\Delta$ along  $\gamma$;
it is diffeomorphic to $\mathbb S^2$.
The line distribution $\ell$ lifts to a line distribution on the doubling.
The latter contradicts the hairy ball theorem.\qeds

This proof was suggested nearly simultaneously 
by Steven Sivek 
and Damiano Testa \cite{two-disks}.

Note that the same proof works when $\Sigma_i$ are oriented open surfaces and $\gamma$ cuts a compact domain in each $\Sigma_i$.

There are examples of three disks $\Sigma_1$, $\Sigma_2$, and $\Sigma_3$
with a common closed curve $\gamma$ such that there is
no triple of points $p_i\in\Sigma_i$ with parallel tangent planes.
Such examples can be found among ruled surfaces \cite{three-disks}.

\parbf{Second derivative bounds first.}
Observe that the gradient $v=\nabla f$ is almost parallel;
that is, there is a constant $C$ such that $|\nabla_u v|\le C$ for any unit tangent vector $u$.
In particular, the parallel translation of $v(p)$ around a circle has to be close to $v(p)$.
If $|v(p)|$ is large, the latter contradicts the Gauss--Bonnet formula.

\medskip

I learned this problem from Christopher Criscitiello \cite{ccriscitiello}.

\csname @openrightfalse\endcsname
\chapter{Comparison geometry}

In this chapter, we consider Riemannian manifolds with curvature bounds.

This chapter is very demanding;
we assume that the reader is familiar with   
the shape operator and second fundamental form, 
equations of Riccati and Jacobi,
comparison theorems,
and Morse theory.
The classical book \cite{cheeger-ebin} covers all the  necessary  material.

\subsection*{Geodesic immersion\hard}
\label{Geodesic immersion}

An isometric immersion $\iota\:N\looparrowright M$ from one Riemannian manifold to another is called \index{totally geodesic}\emph{totally geodesic} if it maps any geodesic in $N$ to a geodesic in $M$.

\begin{pr}
Let $M$ and $N$ be simply-connected positively-curved Riemannian manifolds and $\iota\:N\looparrowright M$ a totally geodesic immersion.
Assume that 
\[\dim N>\tfrac 12\cdot \dim M.\]
Prove that $\iota$ is an embedding.
\end{pr}

\parit{Semisolution.}
Set $n=\dim N$, $m=\dim M$.

Choose a smooth increasing strictly concave function $\phi$.
Consider the function $f=\phi\circ\dist_N$,
where $\dist_N(x)$ denotes the distance from $x\in M$ to $N$.

Note that if $f$ is smooth at a point $x\in M$, then the Hessian of $f$ at $x$ (briefly $\Hess_xf$) has at least $n+1$ negative eigenvalues.

Moreover, at any point $x\notin \iota(N)$ the same holds in the barrier sense\label{page:barrier sense}.
That is, there is a smooth function $h$ defined on $M$
such that $h(x)=f(x)$, $h(y)\ge f(y)$ for any $y$
and $\Hess_xh$ has at least $n+1$ negative eigenvalues.

Use that $m< 2\cdot n$ and the described property to prove the following
analog of Morse lemma for $f$.

\begin{cl}{$({*})$}
 Given $x\notin \iota(N)$, there is a neighborhood $U\ni x$ such that the set 
\[U_-=\set{z\in U}{f(z)<f(x)}\] is simply-connected.
\end{cl}

Since $M$ is simply-connected,
any closed curve in $\iota(N)$
can be contracted by a disc, say $s_0\:\mathbb D\to M$.

Applying the claim $({*})$, one can construct an $f$-decreasing homotopy $s_t$ that starts at $s_0$ and ends in $\iota(N)$;
that is, there is 
a homotopy $s_t\:\mathbb D\z\to M$, $t\in [0,1]$ 
such that $s_t(\partial \mathbb D)\subset \iota(N)$ for any $t$ 
and $s_1(\mathbb D)\subset \iota(N)$.
It follows that $\iota(N)$ is simply-connected.

Finally, assume that $a$ and $b$ are distinct points in $N$ such that $\iota(a)\z=\iota(b)$.
If $\gamma$ is a path from $a$ to $b$ in $N$ then the loop $\iota\circ\gamma$ is not contractible in $\iota(N)$.
Therefore if $\iota\:N\to M$ has a self-intersection,
then the image
$\iota(N)$ is not simply-connected.
Hence the result follows.\qeds

The statement was proved by 
Fuquan Fang, 
S\'ergio Mendon\c{c}a,
and Xiaochun Rong \cite{FMR}.
The main idea was discovered by 
Burkhard Wilking \cite{wilking-2003}.

\subsection*{Geodesic hypersurface}
\label{Geodesic hypersurface}

The totally geodesic embedding is defined before the previous problem.

\begin{pr}
Assume a compact connected positively-curved manifold $M$ has a totally geodesic embedded hypersurface.
Show that either $M$ or its double covering is homeomorphic to the sphere.
\end{pr}

\subsection*{If convex, then embedded}
\label{If convex then embedded} 

\begin{pr}
Let $M$ be a complete simply-connected Riemannian manifold 
with non-positive curvature 
and dimension at least $3$.
Prove that any immersed locally convex
compact hypersurface $\Sigma$ in $M$ is embedded.
\end{pr}

Let us summarize some statements about complete simply-connected Riemannian manifolds 
with non-positive curvature.

By the Cartan--Hadamard theorem, for any point $p\in M$
the exponential map $\exp_p\:\T_p\to M$ is a diffeomorphism.
In particular, $M$ is diffeomorphic to the Euclidean space of the same dimension.
Moreover, any geodesic in $M$ is minimizing,
and any two points in $M$ are connected by a unique minimizing geodesic,

Further, $M$ is a $\CAT(0)$ space; that is, it satisfies a global angle comparison which we are about to describe.
Let $[xyz]$ be a triangle in $M$;
that is, $[xyz]$ is formed by three distinct points $x,y,z$ pairwise connected by geodesics.
Consider its model triangle $[\tilde x\tilde y\tilde z]$ in the Euclidean plane;
that is, $[\tilde x\tilde y\tilde z]$ has the same side lengths as $[xyz]$.
Then each angle in $[xyz]$ cannot exceed the corresponding angle in $[\tilde x\tilde y\tilde z]$.
This inequality can be written as
\[\tilde\measuredangle(y\,^x_z)\ge\measuredangle\hinge yxz,\]
where $\measuredangle\hinge yxz$ denotes the angle of the hinge $\hinge yxz$ formed by two geodesics $[yx]$ and $[yz]$ 
and $\tilde\measuredangle(y\,^x_z)$ denotes the corresponding angle in the model triangle $[\tilde x\tilde y\tilde z]$.

From this comparison, it follows that any connected closed locally convex sets in $M$ is globally convex.
In particular, if $\Sigma$ is embedded then it bounds a convex set.

\subsection*{Immersed ball\hard}
\label{Immersed ball}

\begin{pr}
Prove that any immersed locally convex
hypersurface $\iota\:\Sigma\looparrowright M$
in a compact positively-curved manifold $M$ of dimension $m\ge 3$ is the boundary of an immersed ball. 
That is, there is an immersion of a closed ball $f\:\bar B^m\looparrowright M$ and a diffeomorphism $h\:\Sigma\to\partial \bar B^m$
such that $\iota=f\circ h$.
\end{pr}

\subsection*{Minimal surface in the sphere}
\label{minimal surface}\label{almgren} 

A  smooth $n$-dimensional surface $\Sigma$ in
an $m$-dimensional Riemannian manifold $M$ is called \index{minimal surface}\emph{minimal}
if it locally minimizes the $n$-dimensional area;
that is, sufficiently small regions of $\Sigma$ do not admit area-decreasing deformations with a fixed boundary.

The minimal surfaces can be also defined via the mean curvature vector as follows.
Let $\T=\T\,\Sigma$ and $\mathrm{N}=\mathrm{N}\,\Sigma$ denote the tangent and the normal bundle respectively.
Let $s$ denote the second fundamental form of $\Sigma$;
it is a quadratic from on $\T$ with values in $\mathrm{N}$,
see the remark after problem ``Hypercurve'' below. 
Given an orthonormal basis $(e_i)$ in $\T_x$, set 
$$H_x=\sum_i s(e_i,e_i).$$
The vector $H_x$ lies in the normal space $\mathrm{N}_x$
and does not depend on the choice of orthonormal basis $(e_i)$.
This vector $H_x$ is called the mean curvature vector at $x\in \Sigma$. 

We say that $\Sigma$ is \index{minimal surface}\emph{minimal} if $H\equiv 0$.

\begin{pr}
Let $\Sigma$ be a closed $n$-dimensional 
minimal surface
in the unit $m$-dimensional sphere $\mathbb{S}^m$.
Prove that
$\vol_n \Sigma\ge \vol_n \mathbb{S}^n$.
\end{pr}

\subsection*{Hypercurve}
\label{codim=2}

The Riemann curvature tensor $R$
can be viewed as an operator $\text{\bf R}$ on the space of tangent bi-vectors $\bigwedge^2 \T$;
it is uniquely defined by the identity
$$\langle\mathbf{R}(X\wedge Y),V\wedge W\rangle
=
\langle R(X,Y)V,W\rangle.$$
The operator $\mathbf{R}\:\bigwedge^2 \T\to \bigwedge^2 \T$ is called the \index{curvature operator}\emph{curvature operator} and it is said to be {}\emph{positive definite} if
$\langle\mathbf{R}(\phi),\phi\rangle>0$ for all non-zero
bi-vector $\phi\in\bigwedge^2 \T$.

\begin{pr}
Let $M^m\hookrightarrow \RR^{m+2}$ be a closed smooth $m$-dimensional
submanifold and let  $g$ be the  induced Riemannian metric on $M^m$.
Assume that sectional curvature of $g$ is positive.
Prove that the curvature operator of $g$ is positive definite.
\end{pr}

The second fundamental form for manifolds of arbitrary codimension which we are about to describe might help to solve this problem.

Let $M$ be a smooth submanifold in $\RR^m$.
Given a point $p\in M$, denote by $\T_p$ and $\mathrm{N}_p=\T_p^\bot$
the tangent and normal space of $M$ at $p$.
The \index{second fundamental form}\emph{second fundamental form}\label{page:second fundamental form} of $M$ at $p$ is defined by 
\[s(X,Y)=(\nabla_X Y)^\bot,\]
where $(\nabla_X Y)^\bot$ a denotes the orthogonal projection of covariant derivative $\nabla_X Y$ onto the normal bundle.

The curvature tensor of $M$ can be found from the second fundamental form using the following  formula
\[\langle R(X,Y)V,W\rangle=\langle s(X,W),s(Y,V)\rangle-\langle s(X,V),s(Y,W)\rangle,\]
which is a direct generalization of the formula for Gauss curvature of a surface.

\subsection*{Horo-sphere}
\label{Horosphere}

We say that a Riemannian manifold has negatively pinched sectional curvature if its sectional curvatures at all points in all sectional directions lie in an interval $[-a^2, -b^2]$, for fixed constants $a>b>0$.

Let $M$ be a complete Riemannian manifold
and $\gamma$ a ray in $M$; 
that is, $\gamma\:[0, \infty)\to M$ is a minimizing unit-speed geodesic.

The \label{page:Busemann function}\index{Busemann function}\emph{Busemann function} $\bus_\gamma\:M\to\RR$ is defined by
$$\bus_\gamma(p)=\lim_{t\to\infty}\left(|p-\gamma(t)|_M-t\right).$$
By the triangle inequality, 
the expression under the limit is non-increasing in $t$; 
therefore  the limit above is defined for any $p$.

A \index{horo-sphere}\emph{horo-sphere} in $M$ is defined as a level set of a Busemann function
on $M$.

We say that a complete Riemannian manifold $M$ has \index{polynomial volume growth}\emph{polynomial volume growth} if, for some (and therefore any) $p\in M$, we have 
$$\vol B(p,r)_M\z\le C\cdot (r^k+1),$$ 
where $B(p,r)_M$ denotes the ball in $M$ and  $C$, $k$ are constants.

\begin{pr} Let $M$ be a complete simply-connected manifold with negatively pinched sectional curvature
and $\Sigma\subset M$ an horo-sphere in $M$.
Show that
$\Sigma$ with the induced intrinsic metric 
has polynomial volume growth.
\end{pr}

\subsection*{Number of conjugate points}
\label{Number of conjugate points}

Recall that points $p$ and $q$ on a geodesic $\gamma$ are called \index{conjugate points}\emph{conjugate} if there exists a non-zero Jacobi field along $\gamma$ that vanishes at $p$ and $q$. 

\begin{pr}
Let $s\:N\to M$ be a Riemannian submersion.
Suppose $N$ has nonpositive sectional curvature.
Show that any point $p$ in $M$ has at most $k=\dim N-\dim M$ conjugate points on any geodesic $\gamma\ni p$.
\end{pr}

\subsection*{Minimal spheres}
\label{Minimal spheres}

Recall that two subsets $A$ and $B$ in a metric space $X$ are called \index{equidistant sets}\emph{equidistant} if the distance function $\dist_A\:X\to\RR$ is constant on $B$ and $\dist_B$ is constant on $A$.

The minimal surfaces are defined on page \pageref{minimal surface}.

\begin{pr}
Show that a 
$4$-dimensional
compact 
positively-curved 
Riemannian manifold 
cannot contain an infinite number of  mutually
 equidistant minimal 2-spheres.
\end{pr}

\subsection*{Positive curvature and symmetry\thm}
\label{kleiner-hopf} 

\begin{pr}
Assume that $\mathbb S^1$ acts isometrically on a closed $4$-dimensional Riemannian manifold with positive sectional curvature.
Show that the action 
has at most $3$ isolated fixed points.
\end{pr}

The following statement might be useful.
If $(M,g)$ is a Riemannian manifold with sectional curvature $\ge \kappa$ that admits a continuous isometric action of a compact group $G$, 
then the quotient space $A=(M,g)/G$ is an Alexandrov space with curvature $\ge \kappa$;
that is, the conclusion of the Toponogov comparison theorem holds in $A$. 

For more on Alexandrov geometry see \cite{akp}.

\subsection*{Energy minimizer}
\label{Energy minimizer}

Let $F$ be a smooth map from a closed Riemannian manifold $M$ to a Riemannian manifold $N$.
The energy functional of $F$ is defined by
\[E(F)=\int\limits_{x\in M} |d_xF|^2.\]
We assume that  
\[|d_xF|^2=\sum_{i,j}a_{i,j}^2,\]
where $(a_{i,j})$ denote the components 
of the differential $d_xF$ 
written in the orthonormal bases of the tangent spaces $\T_xM$ and $\T_{F(x)}N$.

\begin{pr}
Show that the identity map on $\RP^m$ is 
energy-minimizing in its homotopy class.
Here we assume that $\RP^m$ is equipped with the canonical metric.
\end{pr}

\subsection*{Curvature against injectivity radius\thm}
\label{scalar-curv}

\begin{pr} 
Let $(M,g)$ be a closed 
Riemannian $m$-dimensional manifold.
Assume average of sectional curvatures over all sectional directions of $(M,g)$ is $1$. 
Show that the injectivity radius of $(M,g)$ is at most $\pi$.
\end{pr}

Solutions to this and the previous problem use the fact that geodesic flow on the tangent bundle to a Riemannian manifold preserves the volume form; this is a corollary of Liouville's theorem about phase volume.

\subsection*{Approximation of a quotient}

\begin{pr}
Let $(M,g)$ be a compact Riemannian manifold 
and $G$ a compact Lie group acting by isometries on $(M,g)$.
Construct a sequence of metrics $g_n$ on a fixed manifold $N$ such that $(N,g_n)$ converges to the quotient space $(M,g)/G$ in the sense of Gromov--Hausdorff.
\end{pr}

\subsection*{Polar points\many}
\label{milka-polar} 

\begin{pr}
Let $M$ be a compact Riemannian manifold with sectional curvature at least $1$ 
and dimension at least $2$. 
Prove that for any point $p\in M$ there is a point $p^*\in M$ such that 
\[|p-x|_M+|x-p^*|_M\le \pi\]
for any $x\in M$.
\end{pr}

\subsection*{Isometric section\hard}
\label{Isometric section}

\begin{pr}
Let $M$ and $W$ be compact Riemannian manifolds,
$\dim W>\dim M$,
and $s\:W\to M$ a Riemannian submersion.
Assume that $W$ has positive sectional curvature.
Show that $s$ does not admit an isometric section;
that is, there is no isometric embedding $\iota\:M\hookrightarrow W$ such that $s\circ\iota(p)=p$ for any $p\in M$.
\end{pr}

\subsection*{Warped product}
\label{Warped product}
\label{page:warped product}

Let $(M,g)$ and $(N,h)$ be Riemannian manifolds 
and $f$ a smooth positive function defined on $M$.
Consider the product manifold $W\z=M\times N$.
Given a tangent vector 
$X\z\in \T_{(p,q)} W
\z=\T_p M\times \T_p N$, denote by 
$X_M\z\in \T M$ and $X_N\z\in \T N$ its projections.
Let us equip $W$ with the Riemannian metric defined by
\[s(X,Y)=g(X_M,Y_M)+f^2\cdot h(X_N,Y_N).\]
The obtained Riemannian manifold $(W,s)$ is called a \index{warped product}\emph{warped product} of $M$ and $N$ with respect to $f\:M\to \RR$;
it can be written as  
\[(W,g)\z=(N,h)\times_f(M,g).\]

\begin{pr}
Let $M$ be an oriented 3-dimensional Riemannian manifold with positive scalar curvature 
and $\Sigma\subset M$ an oriented smooth hypersurface that is area minimizing in its homology class.

Show that there is a positive smooth function $f\:\Sigma\to \RR$
such that the warped product $\mathbb S^1\times_f \Sigma$
has positive scalar curvature;
here $\Sigma$ is equipped with the Riemannian metric
induced from $M$.
\end{pr}

\subsection*{No approximation\many}
\label{No approximation}

\begin{pr}
Prove that 
if $p\not=2$,
then $\RR^m$ 
equipped with the metric induced by the $\ell^p$-norm 
cannot be a
Gromov--Hausdorff limit of
$m$-dimensional Riemannian manifolds $(M_n,g_n)$ with $\Ric_{g_n}\z\ge C$ for a constant $C$.
\end{pr}

\subsection*{Area of spheres}
\label{Area of spheres}

\begin{pr}
Let $M$ be a complete non-compact Riemannian manifold 
with non-negative Ricci curvature and $p\in M$.
Show that there is $\eps>0$ such that 
$$\area\left[\partial B(p,r)\right]>\eps$$
for all sufficiently large $r$.
\end{pr}

\subsection*{Flat coordinate planes}
\label{Flat coordinate planes}

\begin{pr}
Let $g$ be a complete Riemannian metric on $\RR^3$ such that the coordinate planes $x=0$, $y=0$, and $z=0$ are flat and totally geodesic.
Assume the sectional curvature of $g$ is either non-negative or non-positive.
Show that in both cases $g$ is flat. 
\end{pr}

\subsection*{Two-convexity\many}
\label{Two-convexity}

An open subset $V$ with smooth boundary in the Euclidean space  
is called \index{two-convex set}\emph{two-convex} if at most one principal curvature in the outward direction to $V$ is negative.

The two-convexity of $V$ is equivalent to the following property.
For any plane $\Pi$ and any closed curve $\gamma$ in the intersection  $V\cap \Pi$,
if $\gamma$ is contactable in $V$ then it is contactable in $\Pi\cap V$.

\begin{pr}
Let $K$ be a closed set bounded by a smooth surface
in $\RR^4$.
Assume that $K$ contains two coordinate planes 
$$\{(x,y,0,0)\in\RR^4\}
\quad
\text{and}
\quad
\{(0,0,z,t)\in\RR^4\}$$
in its interior 
and also belongs to the closed $1$-neighborhood of these two planes.

Show that the complement to $K$ is not two-convex.
\end{pr}

\subsection*{Convex lens\thm}
\label{Convex lens}

\begin{pr} Let $D$ and $D'$ be two smooth discs with a common boundary that bound a convex set (a lens) $L$ in a positively-curved 3-dimensional Riemannian manifold $M$.
Assume that the discs meet at a small angle.

\begin{wrapfigure}{o}{25 mm}
\vskip-2mm
\centering
\includegraphics{mppics/pic-301}
\vskip0mm
\end{wrapfigure}

Show that the integral 
\[\int\limits_{D}k_1\cdot k_2\]
is small; here $k_1$ and $k_2$ denote the principal curvatures of $D$.
\end{pr}

The expected solution uses the following relative version of the Bochner formula.
Let $M$ be a Riemannian manifold with boundary $\partial M$.
If $f\:M\to \RR$ is a smooth function that vanishes on $\partial M$,
then 
\[\int\limits_M \left(|\Delta f|^2
-|\Hess f|^2
-\langle\mathrm{Ric}(\nabla f),\nabla f\rangle\right)
=\int\limits_{\partial M}
H\cdot|\nabla f|^2.\]
Here we denote by $\Ric$ the Ricci curvature of $M$ 
and by $H$ the mean curvature of $\partial M$, we assume that $H\ge 0$ is the boundary is convex.

\subsection*{Small-twist}
\label{Small-twist}

\begin{pr}
Show that any smooth closed manifold admits an immersion into the unit ball in a Euclidean space of sufficiently large dimension
with all its normal curvatures less than $\sqrt{3}$.
\end{pr}

\section*{Semisolutions}

\parbf{Geodesic hypersurface.}
Let $\Sigma$ be the totally geodesic embedded hypersurface in the positively-curved manifold $M$.
Without loss of generality, we can assume that $\Sigma$ is connected.%
\footnote{In fact, by Frankel's theorem [see page \pageref{page:frankel}] $\Sigma$ is connected.}

The complement $M\setminus\Sigma$ has one or two connected components.
First let us show that if the number of connected components is two, 
then $M$ is homeomorphic to a sphere.

By cutting $M$ along $\Sigma$ 
we get two manifolds
with geodesic boundaries.
It is sufficient to show that each of them is homeomorphic to a Euclidean ball.

Choose one of these manifolds; denote it by $N$.
Denote by $f\:N\z\to\RR$ the distance functions to the boundary $\partial N$.
By the Riccati equation $\Hess f\le 0$ at any smooth point,
and for any point the same holds in the barrier sense [defined on page \pageref{page:barrier sense}].
It follows that $f$ is concave.

Choose an increasing strictly concave function $\phi\:\RR\to\RR$.
Note that $\phi\circ f$ is strictly concave in the interior of $N$.

Choose a compact subset $K$ in the interior of $N$ and
smooth $\phi\circ f$ in a neighborhood of $K$ keeping it concave. 
This can be done by applying the smoothing theorem of Robert Greene and Hung-Hsi Wu \cite[Theorem~2]{greene-wu}.

After the smoothing, the obtained strictly concave function, say $h$, has a single critical point which is its maximum.
In particular, by Morse lemma, we get that if the set  
\[N'_s=\set{x\in N}{h(x)\ge s}\]
is not empty and lies in $K$ then it is diffeomorphic to a Euclidean ball.

For appropriately chosen set $K$ and the smoothing $h$, the set $N'_s$ can be made arbitrarily close to $N$;
moreover, its boundary $\partial N'_s$ can be made $C^\infty$-close to $\partial N$.
It follows that $N$ is diffeomorphic to a Euclidean ball.
This finishes the proof of the first case.

Now assume $M\setminus\Sigma$ is connected.
In this case, there is a double covering $\tilde M$ of $M$ that induces a double covering $\tilde\Sigma$ of $\Sigma$,
so $\tilde M$ contains a geodesic hypersurface $\tilde\Sigma$ that divides $\tilde M$ into two connected components. 
From the case which has already been considered, $\tilde M$ is homeomorphic to a sphere;
hence the second case follows.
\qeds

The problem was suggested by Peter Petersen.

\parbf{If convex, then embedded.}
Set 
\[m=\dim \Sigma=\dim M-1.\]

Given a point $p$ on $\Sigma$, denote by $p_r$ the point at distance $r$ from $p$
that lies on the geodesic starting at $p$ in the outer normal direction to $\Sigma$.
Note that for fixed $r\ge 0$,
the points $p_r$ sweep an immersed locally convex hypersurface which we denote by $\Sigma_r$.

\begin{wrapfigure}{o}{61 mm}
\vskip-2mm
\centering
\includegraphics{mppics/pic-302}
\end{wrapfigure}

Choose $z\in M$. 
Denote by $d$ the maximal distance from $z$ to the points in $\Sigma$.
Note that 
any point on $\Sigma_r$
lies at a distance at least $r-d$ from $z$.

By comparison, 
\[\measuredangle\hinge{p_r}zp\le \arcsin\tfrac dr.\]
In particular, for large $r$, 
each infinite geodesic starting at $z$ intersects $\Sigma_r$ transversally.

The space of geodesics starting at $z$ can be identified with the sphere $\mathbb{S}^m$.
Therefore the map that sends a point $x\in \Sigma_r$ to a geodesic from $z$ thru $x$ induces a local diffeomorphism $\phi_r\:\Sigma\z\to\mathbb{S}^m$.

Since $m\ge 2$, the sphere $\mathbb{S}^m$ is simply-connected.
Since $\Sigma$ is connected, the map $\phi_r$ is a diffeomorphism.
Therefore $\Sigma_r$ is star-shaped with a center at $z$.
In particular, $\Sigma_r$ is embedded.
Since $\Sigma_r$ is locally convex, it bounds a convex region and is embedded.

Consider the set $W$ of reals $r\ge 0$ such that $\Sigma_r$ is embedded and bounds a convex region.
Note that $W$ is open and closed in $[0,\infty)$.
Therefore $W=[0,\infty)$, and, in particular, $\Sigma_0=\Sigma$ is embedded.\qeds

The problem is due to Stephanie Alexander \cite{alexander}.

\parbf{Immersed ball.}
Equip $\Sigma$ with the induced intrinsic metric.
Denote by $\kappa$ the lower bound for principal curvatures of $\Sigma$.
Note that we can assume that $\kappa>0$.

Choose a sufficiently small $\eps=\eps(M,\kappa)>0$.
Given $p\in \Sigma$, denote by $\Delta(p)$ the $\eps$-ball in $\Sigma$ centered at $p$.
Consider the lift $\tilde h_p\:\Delta(p)\to \T_{h(p)}$ along the exponential map $\exp_{h(p)}\:\T_{h(p)}\to M$.
More precisely:
\begin{enumerate}
\item Connect each point $q\in \Delta(p)\subset \Sigma$ to $p$
by a minimizing geodesic  path $\gamma_q\:[0,1]\to \Sigma$
\item Consider the lifting $\tilde\gamma_q$ in $\T_{h(p)}$; 
that is, $\tilde\gamma_q$ is the curve such that $\tilde\gamma_q(0)=0$ 
and $\exp_{h(p)}\circ\,\tilde\gamma_q(t)=\gamma_q(t)$ for each $t\in[0,1]$.
 \item Set $\tilde h(q)=\tilde\gamma_q(1)$.
\end{enumerate}

Show that all the hypersurfaces $\tilde h_p(\Delta(p))\subset \T_{h(p)}$ have principal curvatures at least $\tfrac\kappa2$.

Use the same idea as in the solution of ``Immersed surface'' [page ~\pageref{Immersed surface}] to show that 
one can fix $\eps\z=\eps(M,\kappa)>0$ such that the restriction of $\tilde h_p|_{\Delta(p)}$ is injective.
Conclude that the restriction $h|_{\Delta(p)}$ is injective for any $p\in\Sigma$.
(Here we use that $m\ge 3$.)

Now consider locally equidistant surfaces $\Sigma_t$ in the inward direction for small $t$. 
The principal curvatures of $\Sigma_t$ remain at least $\kappa$ in the barrier sense;
that is, at any point $p$, the surface $\Sigma_t$ can be supported by a smooth surface with principal curvatures at least $\kappa$ at $p$.
By the same argument as above, any $\eps$-ball in $\Sigma_t$ is embedded.

We get a one-parameter family of locally convex locally equidistant surfaces $\Sigma_t$
for $t$ in the maximal interval $[0,a]$,
where the surface $\Sigma_a$ degenerates to a point, say $p$. 

To construct the immersion $\partial \bar B^m\looparrowright M$,
take the point $p$ as the image of the center $\bar B^m$ 
and take the surfaces $\Sigma_t$ as the restrictions of the  embedding to the spheres;
the existence of the immersion follows from the Morse lemma.\qeds

\begin{wrapfigure}{r}{30 mm}
\vskip0mm
\centering
\includegraphics{mppics/pic-304}
\end{wrapfigure}

As you see from the picture, 
the analogous statement does not hold in the two-dimensional case.

The proof presented above was indicated in the lectures of Mikhael Gromov \cite{gromov-SGMC} and written rigorously by Jost Eschenburg \cite{eschenburg}.

A variation of Gromov's proof 
was obtained independently by Ben Andrews \cite{andrews}.
Instead of equidistant deformation, 
he uses the so-called \index{inverse mean curvature flow}\emph{inverse mean curvature flow};
this way one has to perform some calculations to show that convexity survives in the flow, 
but one does not have to worry about non-smoothness of the hypersurfaces~$\Sigma_t$.

\parbf{Minimal surface in the sphere.}
Choose a  geodesic $n$-dimensional sphere $\tilde\Sigma=\mathbb{S}^n\subset \mathbb{S}^m$.

Denote by $U_r$ and $\tilde U_r$ the closed tubular $r$-neighborhood 
of $\Sigma$ and $\tilde\Sigma$ in $\mathbb{S}^m$ respectively.

Note that 
\[U_{\frac\pi2}=\tilde U_{\frac\pi2}=\mathbb{S}^m.
\leqno({*})\]
Indeed, clearly $\tilde U_{\frac\pi2}=\mathbb{S}^m$.
If $U_{\frac\pi2}\ne\mathbb{S}^m$, fix $x\in \mathbb{S}^m\setminus U_r$.
Choose a closest point $y\in \Sigma$ to $x$.
Since $r=|x-y|_{\mathbb{S}^m}>\tfrac\pi2$ the $r$-sphere $\mathrm{S}_r\subset \mathbb{S}^m$ with center $x$ is concave.
Note that $\mathrm{S}_r$ supports $\Sigma$ at $y$;
in particular, the mean curvature vector of $\Sigma$ at $y$ cannot vanish --- a contradiction.

By the Riccati equation, 
\[H_r(x)\ge \tilde H_r\] 
for any $x\in \partial U_r$,
where $H_r(x)$ denotes the mean curvature of $\partial U_r$  at a point $x$
and $\tilde H_r$ is the mean curvature of $\partial\tilde U_r$,
the latter is the same at all points.

Set 
\begin{align*}
a(r)&=\vol_{m-1} \partial U_r,
&
\tilde a(r)&=\vol_{m-1} \partial\tilde U_r,
\\
v(r)&=\vol_m U_r,
&
\tilde v(r)&=\vol_m \tilde U_r.
\intertext{By the coarea formula,}
\tfrac d{dr} v(r)&\aall a(r),
&
\tfrac d{dr}\tilde v(r)&=\tilde a(r).
\end{align*}
Note that
\begin{align*}\tfrac d{dr}a(r)&\le \int\limits_{x\in\partial U_r} H_r(x)\le
\\
&\le a(r)\cdot \tilde H_r,
\end{align*}
and
\begin{align*}
\tfrac d{dr}\tilde a(r)
&= \tilde a(r)\cdot \tilde H_r.
\intertext{It follows that}
\frac {v''(r)}{v(r)}&\le \frac {\tilde v''(r)}{\tilde v(r)}
\end{align*}
for almost all $r$. 
Therefore
\[v(r)\le\frac{\area\Sigma}{\area \tilde\Sigma}\cdot \tilde v(r)\]
for any $r>0$.

According to $({*})$,
\[v(\tfrac\pi2)=\tilde v(\tfrac\pi2)=\vol\mathbb{S}^m.\]
Hence the result follows.\qeds

This problem is the geometric lemma in the proof given by Frederick Almgren of his isoperimetric inequality \cite{almgren}.
The argument is similar to 
the proof of isoperimetric inequality for manifolds with positive Ricci curvature
given by Mikhael Gromov \cite{gromov-apendix}.

\parbf{Hypercurve.}
Choose $p\in M$.
Denote by $s$ 
the second fundamental form of $M$ at $p$.
Recall that $s$ is a symmetric bilinear form on the tangent space $\T_pM$ of $M$ with values in the normal space $\mathrm{N}_pM$ to $M$, see page~\pageref{page:second fundamental form}.

By the Gauss formula
\[\langle R(X,Y)Y,X\rangle=\langle s(X,X),s(Y,Y)\rangle-\langle s(X,Y),s(X,Y)\rangle.\]
Since the sectional curvature of $M$ is positive, 
we get
\[\<s(X,X),s(Y,Y)\> > 0\leqno({*})\]
for any pair of non-zero vectors $X,Y\in\T_pM$.

The normal space $\mathrm{N}_pM$ is two-dimensional.
By $({*})$ there is an orthonormal basis $e_1,e_2$ in $\mathrm{N}_pM$ 
such that the real-valued quadratic forms 
\begin{align*}
s_1(X,X)&=\<s(X,X),e_1\>,
&
s_2(X,X)&=\<s(X,X),e_2\>
\end{align*}
are positive definite.

Note that the curvature operators $\mathbf{R}_1$ and $\mathbf{R}_2$ 
are defined by the formula
\[\mathbf{R}_{i}(X\wedge Y), V\wedge W\rangle 
=s_i(X,W)\cdot s_i(Y,V)-s_i(X,V)\cdot s_i(Y,W)\]
are positive.
Finally, note that $\mathbf{R}=\mathbf{R}_{1}+\mathbf{R}_{2}$ is the curvature operator of $M$ at $p$.\qeds

The problem is due to Alan Weinstein \cite{weinstein}.
Note that from \cite{micallef-moore}/\cite{boehm-wilking} it follows that the universal covering of $M$ is homeomorphic/\hskip0mm diffeomorphic to a standard sphere.

\parbf{Horo-sphere.}
Set 
$m=\dim \Sigma=\dim M-1$.

Let $\bus\:M\to\RR$ be the Busemann function such that 
\[\Sigma=\bus^{-1}\{0\}.\]
Set  $\Sigma_r=\bus^{-1}\{r\}$, so $\Sigma_0=\Sigma$.

Let us equip each $\Sigma_r$ with the induced Riemannian metric.
Note that all $\Sigma_r$ have bounded curvature.
In particular, the unit balls in $\Sigma_r$ have volume bounded above by a universal constant, say $v_0$.
 
Given $x\in \Sigma$, denote by $\gamma_x$ 
the unit-speed geodesic
such that $\gamma_x(0)=x$ and $\bus(\gamma_x(t))=t$ for any $t$.
Consider the map $\phi_{r}\:\Sigma\to\Sigma_r$ defined by
$\phi_r\:x\mapsto \gamma_x(r)$.
In other words, $\phi_{r}$ is the closest point projection from $\Sigma$ to $\Sigma_r$.

Notice that $\phi_r$ is a bi-Lipschitz map with the Lipschitz constants $e^{a\cdot r}$ and $e^{b\cdot r}$.
In particular, the ball of radius $R$ in $\Sigma$ is mapped by $\phi_r$
to a ball of radius $e^{a\cdot r}\cdot R$ in $\Sigma_r$.
Therefore
\[\vol_m B(x,R)_\Sigma
\le 
e^{m\cdot b\cdot r}\cdot \vol_m B(\phi_r(x),e^{a\cdot r}\cdot R)_{\Sigma_r}\]
for all $R,r>0$.
Taking $R=e^{-a\cdot r}$, we get
\[\vol_m B(x,R)_\Sigma\le v_0\cdot R^{m\cdot \frac ba}\]
for any $R\ge1$. 
Hence the statement follows.
\qeds

The problem was suggested by Vitali Kapovitch.

There are examples of horo-spheres as above with a degree of polynomial growth higher than $m$.
For example, consider the horo-sphere $\Sigma$ in the complex hyperbolic space 
of real dimension $4$.
Clearly, $m=\dim \Sigma=3$, but the degree of its volume growth is $4$.

In this case, $\Sigma$ is isometric to the Heisenberg group.%
\footnote{\index{Heisenberg group}\emph{Heisenberg group}
is the group of $3\times3$ upper triangular matrices of the form
\[\begin{pmatrix}
 1 & a & c\\
 0 & 1 & b\\
 0 & 0 & 1\\
\end{pmatrix}\]
under the operation of matrix multiplication.} 
It is instructive to show that any such metric has volume  growth of degree $4$.

\parbf{Number of conjugate points.}
Choose a geodesic $\gamma$ in $M$ and a point $p\in \gamma$.
Note that $\gamma$ can be lifted to a horizontal geodesic $\bar\gamma$ in $N$.
That is, $\gamma=s\circ\bar\gamma$ and $\bar\gamma$ is perpendicular to the fibers of~$s$ (in particular, $\gamma$ and $\bar\gamma$ have equal speeds).

Observe that each conjugate point of $p$ on $\gamma$ corresponds to a \index{focal point}\emph{focal point} on $\bar\gamma$ to the fiber $F$ over $p$ in $N$;
that is, $\bar\gamma$ lies in a family of geodesics $\bar\gamma_t$ that are perpendicular to $N$ 
such that the corresponding Jacobi field along $\bar\gamma$ vanish at $q$.

Note that $F$ has dimension $k=\dim N-\dim M$.
It remains to prove that any smooth $k$-dimensional submanifold $F$ in a complete nonpositively-curved manifold $N$ has at most $k$ focal points on any geodesic $\bar \gamma$ that is perpendicular to $F$.\qeds

The problem is inspired by the paper of Alexander Lytchak~\cite{lytchak:conjugate}.
Applying it together with the Poincaré recurrence theorem
leads to a solution of the following problem.

\begin{pr}
Let $s\:N\to M$ be a Riemannian submersion.
Suppose $N$ has nonpositive sectional curvature and $M$ is compact.
Show that $M$ has no conjugate points.
\end{pr}

In fact, no compact negatively curved manifold $N$ admits a nontrivial Riemannian submersion $s\:N\to M$~\cite[see Theorem F in][]{zeghib}.

\parbf{Minimal spheres.}
Assuming the contrary,
we can choose a pair of sufficiently close minimal spheres $\Sigma$ and $\Sigma'$ in the 4-dimensional manifold $M$;
we can assume that the distance $a$ between $\Sigma$ and $\Sigma'$ is strictly smaller than the injectivity radius of the manifold.
Note that in this case, there is a unique bijection $\Sigma\to \Sigma'$, denoted by $p\mapsto p'$ such that the distance $|p-p'|_M=a$ for any $p\in\Sigma$.

Let $\iota_p\:\T_p\to\T_{p'}$ be the parallel translation along the (necessarily unique) minimizing geodesic $[pp']$.
Note that there is a pair $(p,p')$ such that $\iota_p(\T_p\Sigma)=\T_{p'}\Sigma'$.
Indeed, if this is not the case, then $\iota_p(\T_p\Sigma)\z\cap\T_{p'}\Sigma'$ forms a continuous line distribution over $\Sigma'$.
Since $\Sigma'$ is a two-sphere, the latter contradicts the hairy ball theorem.

Consider pairs of unit-speed geodesics $\alpha$ and $\alpha'$ 
in $\Sigma$ and $\Sigma'$  
that start at $p$ and $p'$ respectively
and go in parallel directions, say $\nu$ and $\nu'$. 
Set $\ell_\nu(t)=|\alpha(t)-\alpha'(t)|$.

Use the second variation formula together with the lower bound on Ricci curvature
to show that $\ell_\nu''(0)$ has a negative average for all tangent directions $\nu$ to $\Sigma$ at $p$. 
In particular, $\ell_\nu''(0)<0$ for a vector $\nu$ as above.
For the corresponding pair $\alpha$ and $\alpha'$,
it follows that there are points $v=\alpha(\eps)\in\Sigma$ near $p$ 
and $v'=\alpha'(\eps)\in\Sigma'$ near $p'$
such that 
\[|v-v'|<|p-p'|\]
--- a contradiction.\qeds

Likely, any compact positively-curved 
4-dimensional manifold
cannot contain a pair of equidistant spheres.
The argument above implies that the distance between such a pair has to exceed the injectivity radius of the manifold.

The problem was suggested by Dmitri Burago.
Here is a short list of classical problems which use the second variation formula in a similar fashion:

\begin{pr}
Any compact even-dimensional orientable manifold with positive sectional curvature is
simply-connected.
\end{pr}

This is called Synge's lemma \cite{synge}.

\begin{pr}
Any two compact minimal hypersurfaces in a Riemannian manifold with positive Ricci curvature must intersect.
\end{pr}

\begin{pr}
Let $\Sigma_1$ and $\Sigma_2$ be two compact geodesic submanifolds in a manifold with positive sectional curvature $M$ and 
\[\dim \Sigma_1+\dim \Sigma_2\ge \dim M.\] 
Then $\Sigma_1\cap\Sigma_2\ne\emptyset$.
\end{pr}

These two statements have been proved by Theodore Frankel \cite{frankel}.\label{page:frankel}

\begin{pr}
Let $(M,g)$ be a closed Riemannian manifold with negative Ricci curvature.
Prove that $(M,g)$ does not admit an isometric $\mathbb{S}^1$-action.
\end{pr}

This is a theorem of Salomon Bochner \cite{bochner}.

The problem ``Geodesic immersion'' [page \pageref{Geodesic immersion}] can be considered as further development of the idea.

\parbf{Positive curvature and symmetry.}
Let $M$ be a 4-dimensional Riemannian manifold with isometric $\mathbb{S}^1$-action.
Consider the quotient space $X=M/\mathbb{S}^1$.
Note that $X$ is a positively-curved 3-dimensional Alexandrov space.
In particular, the angle $\measuredangle\hinge xyz$ between any two geodesics $[xy]$ and $[xz]$ is defined.
Further, for any non-degenerate triangle $[xyz]$ 
formed by the minimizing geodesics $[xy]$, $[yz]$, and $[zx]$  in $X$ we have
\[\measuredangle\hinge xyz+\measuredangle\hinge yzx+\measuredangle\hinge zxy> \pi.
\leqno({*})\]

Assume that $p\in X$ corresponds to a fixed point $\bar p\in M$ of the $\mathbb{S}^1$-action.
Each direction of a geodesic starting at $p$ in $X$ corresponds to an $\mathbb{S}^1$-orbit of the induced isometric action $\mathbb{S}^1\z\acts\mathbb{S}^3$ on the sphere of unit vectors at $\bar p$.
Any such action is conjugate to the action $\mathbb{S}^1_{p,q}\z\acts\mathbb{S}^3\subset\CC^2$ induced by complex matrices 
$
\left(
\begin{smallmatrix}
z^p&0
\\
0&z^q
\end{smallmatrix}
\right)
$
with $|z|=1$ and relatively prime positive integers $p,q$.
The possible quotient spaces $\Sigma_{p,q}=\mathbb{S}^3/\mathbb{S}^1_{p,q}$ 
have diameter $\tfrac\pi2$ and perimeter of any triangle in $\Sigma_{p,q}$ is at most $\pi$;
this is straightforward to check but requires some work.

Therefore for any three geodesics $[px]$, $[py]$, and $[pz]$ in $X$ we have
\[\measuredangle\hinge pxy+\measuredangle\hinge pyz+\measuredangle\hinge pzx\le \pi.\leqno({*}{*})\]
and
\[\measuredangle\hinge pxy,\ \measuredangle\hinge pyz,\ \measuredangle\hinge pzx\le \tfrac\pi2.\leqno(\asterism)\]

Arguing by contradiction,
assume that there are 4 fixed points $q_1$, $q_2$, $q_3$, and $q_4$.
Connect each pair by a minimizing geodesic $[q_iq_j]$.

Denote by $\omega$ the sum of all 12 angles of the type  $\measuredangle\hinge{q_i}{q_j}{q_k}$.
By $(\asterism)$, each triangle $[q_iq_jq_k]$ is non-degenerate.
Therefore by $({*})$, we have
\[\omega>4\cdot\pi.\]
On the other hand, applying $({*}{*})$ at each vertex $q_i$, we have 
\[\omega\le 4\cdot\pi\]
--- a contradiction.\qeds

The problem is due to 
Wu-Yi Hsiang 
and Bruce Kleiner 
\cite{hsiang-kleiner}.
The connection of this proof to Alexandrov geometry was noticed by Karsten Grove \cite{grove}.
An interesting new twist of the idea 
is given by 
Karsten Grove 
and Burkhard Wilking 
\cite{grove-wilking}.

\parbf{Energy minimizer.}
Denote by $\mathcal{U}$ the unit tangent bundle over $\RP^m$
and by $\mathcal{L}$ the space of projective lines in $\ell\:\RP^1\to \RP^m$.
The spaces $\mathcal{U}$ and $\mathcal{L}$ 
have dimension $2\cdot m-1$ 
and $2\cdot(m-1)$
respectively.

According to Liouville's theorem about phase volume, the identity
\[\int\limits_{v\in \mathcal{U}}f(v)
=
\int\limits_{\ell\in\mathcal{L}}\ \int\limits_{t\in\RP^1} f(\ell'(t))\]
holds for any integrable function $f\:\mathcal{U}\to\RR$.

Let $F\:\RP^m\to\RP^m$ be a smooth map.
Note that up to a multiplicative constant,
the energy of $F$ can be expressed the following way
\[\int\limits_{v\in\mathcal{U}} |dF(v)|^2
=
\int\limits_{\ell\in\mathcal{L}}\ \int\limits_{t\in\RP^1} |[d(F\circ \ell)](t)|^2.\]

Notice that any noncontractible curve in $\RP^m$ has length at least $\pi$.
Therefore, by Bunyakovsky inequality, we get
\begin{align*}
\int\limits_{t\in\RP^1} \left|[d(F\circ \ell)](t)\right|^2
&\ge 
\tfrac1\pi\cdot \left(\,\int\limits_{t\in\RP^1} \left|[d(F\circ \ell)](t)\right|\right)^2=
\\
&=\tfrac1\pi\cdot (\length F\circ\ell)^2\ge
\\
&\ge \pi.
\end{align*}
for any line $\ell\:\RP^1\to \RP^m$.
Hence the result follows.\qeds

\label{page:liouville}
The problem is due to Christopher Croke \cite{croke-energy}. 
He uses the same idea to show that the identity map on $\CP^m$ is energy-minimizing in its homotopy class.
For $\mathbb S^m$, an analogous statement does not hold if $m\ge 3$.
In fact, 
if a closed Riemannian manifold $M$ 
has dimension at least 3 
and $\pi_1M=\pi_2M=0$,
then the identity map on $M$ is homotopic 
to a map with arbitrarily small energy;
the latter was shown by Brian White \cite{white}.

The same idea is used to prove the so-called Loewner's inequality \cite{gromov-filling}.
\begin{pr}
Let $g$ be a Riemannian metric on $\RP^m$ that is conformally equivalent to the canonical metric $g_0$.
Assume that any noncontractible curve in $(\RP^m,g)$ has length at least $\pi$.
Then
\[\vol(\RP^m,g)\ge\vol(\RP^m,g_0).\]

\end{pr}

A more advanced application is the sharp isoperimetric inequality for 
4-dimensional Hadamard manifolds proved by Christopher Croke [see \ncite{croke-4d} and also \ncite{croke-eigenvalue}].

\parbf{Curvature against  injectivity radius.}
We will show that 
if the injectivity radius of the manifold $(M,g)$ is at least $\pi$,
then the average of sectional curvatures on $(M,g)$ is at most $1$.
This is equivalent to the problem.

Choose a point $p\in M$ and two orthonormal vectors $U,V\in\T_p M$.
Consider the geodesic $\gamma$ in $M$ such that $\gamma'(0)=U$.

Set $U_t=\gamma'(t)\in \T_{\gamma(t)}$, and let $V_t\in \T_{\gamma(t)}$ be the parallel translation of $V=V_0$ along $\gamma$.

Consider the field $W_t=\sin t\cdot V_t$ on $\gamma$.
Set 
\begin{align*}
\gamma_\tau(t)&=\exp_{\gamma(t)} (\tau\cdot W_t),
\\
\ell(\tau)&=\length(\gamma_\tau|_{[0,\pi]}),
\\
q(U,V)&=\ell''(0).
\end{align*}
Note that
$$q(U,V)
=
\int\limits_{0}^\pi [(\cos t)^2-K(U_t,V_t)\cdot (\sin t)^2]\cdot dt,
\leqno({*})$$
where $K(U,V)$ is the sectional curvature 
in the direction spanned by $U$ and $V$. 

Since any geodesics of length $\pi$ is minimizing,
we get $q(U,V)\ge0$ for any pair of orthonormal vectors $U$ and $V$.
It follows that the average value of the right-hand side in $({*})$ is non-negative.

By Liouville's theorem about phase volume, while taking the average of $({*})$, we can switch the order of integrals;
therefore  
\[0\le \tfrac\pi2\cdot(1-\bar{K}),\]
where $\bar{K}$ denotes the average of sectional curvatures on $(M,g)$.
Hence the result follows.\qeds

The problem illustrates the idea of Eberhard Hopf \cite{hopf-conjugate}
which was developed further by Leon Green \cite{green}.
Hopf used it to show that a metric on 2-dimensional torus without conjugate points must be flat
and Green showed that the average of sectional curvature on a closed manifold without conjugate points cannot be positive.
For more on the subject see the paper of Mikhael Gromov \cite{gromov2021}.

\parbf{Approximation of a quotient.} The proof will use that for any Riemannian submersion $s\:M\to N$
the lower bound on sectional curvature of $M$ can non exceed the lower bound on sectional curvature of~$N$.

This statement follows from O'Nail's formula \cite[Theorem 3.20]{cheeger-ebin} 
which gives the following relation between sectional curvatures of $M$ ad $N$
\[K_M(X,Y)=K_N(\bar X, \bar Y)+\tfrac34|[\bar X,\bar Y]^V|^2,\]
where $X,Y$ are orthonormal vector fields on $N$, $\bar X, \bar Y$ their horizontal lifts to $M$, $[{*},{*}]$ is the Lie bracket and ${*}^V$ is the projection to the vertical distribution of the submersion.
Indeed, since $\tfrac34|[\bar X,\bar Y]^V|^2\ge 0$, we have $K_M(X,Y)\ge K_N(\bar X, \bar Y)$.

\medskip

Note that $G$ admits an embedding into a compact connected Lie group $H$;
in fact, we can assume that $H=\SO(n)$, for sufficiently large~$n$.

Suppose that the curvature of $(M,g)$ is bounded below by~$\kappa$.

The bi-invariant metric $h$ on $H$ is non-negatively curved.
Therefore for any positive integer $n$ the product $(H,\tfrac1n\cdot h)\times (M,g)$ is a Riemannian manifold with  curvature bounded below by~$\kappa$.

The diagonal action of $G$ on $(H,\tfrac1n\cdot h)\times (M,g)$ is isometric and free. 
Therefore 
the quotient $(H,\tfrac1n\cdot h)\times (M,g)/G$
is a Riemannian manifold, say $(N,g_n)$.
Note that the quotient map $(H,\tfrac1n\cdot h)\times (M,g)\to (N,g_n)$ is a Riemannian submersion.
Therefore $(N,g_n)$ has sectional curvature bounded below by $\kappa$.

It remains to observe that the spaces $(N,g_n)$ converge to $(M,g)/G$ as $n\z\to \infty$.\qeds

The used construction is called \index{Cheeger's trick}\emph{Cheeger's trick}.
The earliest use of this trick I found in \cite{GKM}; 
it was used there to show that Berger's spheres have positive curvature.
This trick is used in the construction of most of the known examples of positively and non-negatively curved manifolds
 \cite{cheeger,aloff-wallach,gromoll-meyer,eschenburg-spaces,bazajkin}.
 
The quotient space $(M,g)/G$ has a finite dimension and its curvature is bounded below in the sense of Alexandrov. 
It is expected that not all finite-dimensional Alexandrov spaces admit approximation by Riemannian manifolds with curvature bounded below
[some partial results are discussed in \ncite{pwz,kapovitch}].

\parbf{Polar points.}
Choose a unit-speed geodesic $\gamma$ that starts at $p$;
that is, $\gamma(0)=p$.
Apply the Toponogov comparison to show that $p^*=\gamma(\pi)$ is a solution. 
\qeds

\parit{Alternative proof.} 
Assume the contrary;
that is, for any $x\in M$ there is a point $x'$ such that 
\[|x-x'|_M+|p-x'|_M>\pi.\]

Given $x\in M$, denote by $f(x)$ a point that maximizes the following sum:
\[|x-f(x)|_M+|p-f(x)|_M.\]
Show that $f$ is uniquely defined and continuous.

Choose sufficiently small $\eps>0$.
Prove that the set $W_\eps=M\setminus B(p,\eps)$ 
is homeomorphic to a ball 
and the map $f$ sends $W_\eps$ into itself.

By Brouwer's fixed-point theorem, $x=f(x)$ for some $x$.
In this case, 
\begin{align*}
|x-f(x)|_M+|p-f(x)|_M&=|p-x|_M\le
\\
&\le\pi
\end{align*}
--- a contradiction.\qeds
 
The problem is due to Anatoliy Milka \cite{milka-poly}.

\parbf{Isometric section.}
Arguing by contradiction, 
assume there is an isometric section $\iota\: M\z\to W$.
It makes it possible to treat $M$ as a submanifold in $W$.

Given $p\in M$, denote by $\mathrm{N}^1_p$ the unit normal space to $M$ at $p$.
Given $v\in \mathrm{N}^1_p$ and a real number $k$,
set 
\[p^{k\cdot v}=s\circ\exp_{p} (k\cdot v).\]
Note that 
\[p^{0\cdot v}=p\ \ \text{for any}\ \  p\in M\ \ \text{and}\ \ v\in \mathrm{N}^1_p.\leqno({*})\]

Choose sufficiently small $\delta>0$.
By Rauch comparison \cite[Corollary 1.36]{cheeger-ebin}, 
if $w\in \mathrm{N}^1_q$ 
is the parallel translation of $v\in \mathrm{N}^1_q$ 
along a minimizing geodesic from $p$ to $q$ in $M$,
then 
\[|p^{k\cdot v}-q^{k\cdot w}|_M<|p-q|_M
\leqno({*}{*})\]
assuming that $|k|\le \delta$.
The same comparison implies that 
\[|p^{k\cdot v}-q^{k'\cdot w}|_M^2<|p-q|_M^2+ (k-k')^2
\leqno(\asterism)\]
assuming that $|k|,|k'|\le \delta$.

Choose $p$ and $v \in \mathrm{N}^1_p$ so that $r=|p-p^{\delta\cdot v}|$ 
takes the maximal possible value.
From $({*}{*})$ it follows that $r>0$.

Let $\gamma$ be the extension of the unit-speed minimizing geodesic from $p^{\delta\cdot v}$ to $p$;
denote by $v_t$ the parallel translation of $v$ to $\gamma(t)$ along $\gamma$. 

We can choose the parameter of $\gamma$ so that $p=\gamma(0)$, $p^{\delta\cdot v}=\gamma(-r)$.
Set $p_n=\gamma(n\cdot r)$, so $p=p_0$ and $p^{\delta\cdot v}=p_{-1}$. 
Choose a large integer $N$ and set $w_n=\delta\cdot(1-\tfrac nN)\cdot v_{n\cdot r}$, $q_n=p_n^{w_n}$, $x_n=\exp_{p_n} (w_n)$, and $q_n=p_n^{w_n}=s(x_n)$.

\begin{figure}[ht!]
\vskip0mm
\centering
\includegraphics{mppics/pic-306}
\end{figure}

By $(\asterism)$, there is a constant $C$ independent of $N$ such that
\[|q_k-q_{k+1}|<r+\tfrac C{N^2}\cdot\delta^2.\]
Therefore 
\[|q_{k+1}-p_{k+1}|>|q_k-p_k|-\tfrac C{N^2}\cdot\delta^2.\]
By induction, we get 
\[|q_N-p_N|>r-\tfrac C{N}\cdot\delta^2.\]
Since $N$ is large we get
\[|q_N-p_N|>0.\]
Note that $w_N=0$;
therefore by $({*})$, we get $q_N=p_N^0=p_N$ --- a contradiction.\qeds

This is the core of Perelman's solution of the Soul conjecture \cite{perelman}.

\parbf{Warped product.}
Given $x\in \Sigma$, denote by $\nu_x$ the normal vector to $\Sigma$ at $x$ that agrees with the orientations of $\Sigma$ and $M$. 
Denote by $\kappa_x$ the non-negative principal curvature of $\Sigma$ at $x$;
since $\Sigma$ is minimal the other principal curvature has to be $-\kappa_x$.

Consider the warped product $W=\mathbb S^1\times_f\Sigma$ for a positive smooth function $f\:\Sigma\to \RR$.
Assume that a point $y\in W$ projects to a point $x\in\Sigma$.
Straightforward computations show that
\begin{align*}
\Sc_W(y)
&=\Sc_\Sigma(x)-2\cdot\frac{\Delta f(x)}{f(x)}=
\\
&=\Sc_M(x)-2\cdot\Ric(\nu_x)-2\cdot\kappa_x^2-2\cdot\frac{\Delta f(x)}{f(x)},
\end{align*}
where $\Sc$ and $\Ric$ denote the scalar and Ricci curvature respectively. 

Consider linear operator $L$ on the space of smooth functions on $\Sigma$ defined by 
\[(Lf)(x)= -[\Ric(\nu_x)+\kappa_x^2]\cdot f(x)-(\Delta f)(x).\]
It is sufficient to find a smooth function $f$ on $\Sigma$ such that
\[f(x)>0 \ \ \text{and}\ \ (Lf)(x)\ge 0\leqno({*})\]
for any $x\in \Sigma$.

Given a smooth function $f\:\Sigma\to \RR$,
extend the field $f(x)\cdot\nu_x$
on $\Sigma$ to a smooth field, say $v$, on whole $M$.
Denote by $\iota_t$ the flow along $v$ for time $t$ and set $\Sigma_t=\iota_t(\Sigma)$.

Denote by $H_t(x)$ the mean curvature of $\Sigma_t$ at $\iota_t(x)$.
Note that the value $(Lf)(x)$ is the derivative of
the function $t\mapsto H_t(x)$  at $t=0$.

Therefore the condition $({*})$
means that we can push $\Sigma$ into one of its sides 
so that its mean curvature does not increase in the first order.
Since $\Sigma$ is area-minimizing,
such push can be obtained by increasing the pressure on one side of $\Sigma$.
(Read further if you are not convinced.)
\qeds

\parit{Formal end of proof.}
Denote by $\delta(f)$ the second variation of area of $\Sigma_t$;
that is, consider the area function $a(t)=\area\Sigma_t$ 
and set $\delta(f)=a''(0)$.
Direct calculations show that
\begin{align*}
\delta(f)
&=
\int\limits_{x\in\Sigma} 
\left(-[\Ric(\nu_x)+\kappa_x^2]\cdot f^2(x)+|\nabla f(x)|^2\right)=
\\
&=\int\limits_{x\in\Sigma} 
(Lf)(x)\cdot f(x).\end{align*}
Since $\Sigma$ is area-minimizing we get 
\[\delta(f)\ge 0\leqno({*}{*})\] for any $f$.

Choose a function $f$ that minimize $\delta(f)$ for all functions such that $\int_{x\in\Sigma} f^2(x)=1$.
Note that $f$ is an eigenfunction 
for the linear operator $L$;
in particular, $f$ is smooth.
Denote by $\lambda$ the eigenvalue of $f$;
by $({*}{*})$,
$\lambda\ge 0$.

Show that $f(x)>0$ at any $x$.
Since $Lf=\lambda\cdot f$, the inequalities $({*})$ follow.\qeds

The problem is due to Mikhael Gromov and Blaine Lawson \cite{gromov-lawson}.
Earlier, in \cite{schoen-yau}, Shing-Tung  Yau and Richard Schoen showed that the same assumptions 
imply the existence of a conformal factor on $\Sigma$ that makes it positively-curved.
Both statements are used the same way
to prove that $\TT^3$ does not admit a metric with positive scalar curvature.

Both statements admit straightforward generalization to higher dimensions
and they can be used to show the non-existence of a metric with positive scalar curvature on $\TT^m$ with $m\le 7$.
For $m=8$, the proof stops working 
since in this dimension area-minimizing hypersurfaces might have singularities.
For example, 
any domain in the cone in $\RR^8$
defined by the identity
\[x^2_1+x^2_2+x^2_3+x^2_4=x^2_5+x^2_6+x^2_7+x^2_8\]
is area-minimizing among the hypersurfaces with the same boundary.

\parbf{No approximation.}
Choose an increasing function $\phi\:(0,r)\to \RR$
such that 
\[\phi''+(n-1)\cdot(\phi')^2+C=0.\]

If $\Ric_{g_n}\ge C$, 
then the function 
$x\mapsto\phi(|q-x|_{g_n})$ is subharmonic.
Therefore for an arbitrary array of points $q_i$ 
and positive reals $\lambda_i$ the function $f_n\:M_n\to \RR$
defined by the formula
$$f(x)=\sum_i\lambda_i\cdot\phi(|q_i-x|_M)$$
is subharmonic.
In particular, $f_n$ does not have a local minimum in $M_n$.

Passing to the limit as $n\to \infty$, we get that any function $f\:\mathbb{R}^m\z\to\mathbb{R}$
of the form 
$$f(x)=\sum_i\lambda_i\cdot\phi(|q_i-x|_{\ell_p})$$
does not have a local minimum in $\mathbb{R}^m$.

Let $e_i$ be the standard basis in $\RR^m$. 
If $p<2$, consider the sum 
$$f(x)=\sum\phi(|q-x|_{\ell_p}),$$
where $q=\pm\eps\cdot e_i$ for all signs and $i$'s.
Straightforward calculations show that if $\eps>0$ is small, then $f$
has a strict local minimum at $0$.

If $p>2$, one has to take the same sum for  $p=\sum_i\pm\eps\cdot e_i$ for all choices of signs.
In both cases, we arrive at a contradiction.
\qeds

The argument given here is close to the proof of Abresch--Gromoll inequality \cite{abresch-gromoll}.
The solution admits a straightforward generalization which implies that if an $m$-dimensional  Finsler manifold $F$ is a Gromov--Hausdorff limit of $m$-dimensional Riemannian manifolds with uniform lower bound on Ricci curvature, then $F$ has to be Riemannian.

An alternative solution to this problem can be built on the almost splitting theorem proved by  Jeff Cheeger and Tobias Colding \cite{cheeger-colding}.

\parbf{Area of spheres.}
Fix $r_0>0$.
Given $r>r_0$, choose a point $q$ on the distance $2\cdot r$ from $p$.

\begin{wrapfigure}{r}{35 mm}
\vskip-5mm
\centering
\includegraphics{mppics/pic-308}
\end{wrapfigure}

Note that any minimizing geodesic from $q$ to a point in $B=B(p,r_0)$
has to cross $S\z=\partial B(p,r)$.
The statement follows since  
\[\vol B\le C_m\cdot r_0\cdot \area S,\]
where $C_m$ is a constant depending only on the dimension $m=\dim M$.
This volume comparison inequality can be proved along the same lines as the Bishop--Gromov inequality.
\qeds

Applying the coarea formula, 
we see that the volume growth of $M$ is at least linear; 
in particular, $M$ has infinite volume.
The latter was proved independently 
by Eugenio Calabi 
and Shing-Tung Yau \cite{calabi,yau-ricci}.

\parbf{Flat coordinate planes.}
Choose $\eps>0$ such that there is a unique geodesic between any two points at distance $<\eps$ from the origin of $\RR^3$.

Consider three points $a$, $b$, and $c$ on the coordinate lines that are $\eps$-close 
to the origin.
The following observation is the key to the proof.

\begin{cl}{$({*})$}
There is a solid flat geodesic triangle in $(\RR^3,g)$ with vertices at $a$, $b$, and $c$.
\end{cl}

Since the coordinate planes are totally geodesic, 
the parallel translation along a coordinate line preserves the directions tangent to a coordinate plane.
Since the parallel translation preserves the angles between vectors, the angles between coordinate planes in $(\RR^3,g)$ are constant.

It follows that the angles of the triangle $[abc]$ coincide with its \emph{model angles},
that is, the angles in the plane triangle with the same sides.

Both curvature conditions imply that the triangle $[abc]$ bounds a solid flat geodesic triangle in   $(\RR^3,g)$.

Use the family of constructed flat triangles to show that at any $x$ point in the $\tfrac\eps{10}$-neighborhood of the origin the sectional curvature vanishes in an open set of sectional directions.
The latter implies that the curvature is identically zero 
in this neighborhood.

Move the origin and apply the same argument locally.
This way we get that the curvature is identically zero everywhere.
\qeds

This problem is based on a lemma discovered by Sergei Buyalo [Lemma 5.8 in \ncite{buyalo}; see also \ncite{andersson-howard} and \ncite{panov-petrunin}].

\parbf{Two-convexity.}
\textit{Morse-style solution.}
Choose $(x,y,z,t)$-coordinates in~$\RR^4$.

Consider a generic linear function $\ell\:\RR^4\to\RR$ that is close to the sum of coordinates $x+y+z+t$.
Note that $\ell$
has non-degenerate critical points on $\partial K$ and all its critical values are different.

For each $s$ consider the set 
$$W_s=\set{w\in \RR^4\setminus K}{\ell(w)<s}.$$
Note that $W_{-1000}$ contains a closed curve, say $\alpha$, that is contactable in $\RR^4\setminus K$, but not constructible in $W_{-1000}$.

Set $s_0$ to be the infimum of the values $s$ such that
the $\alpha$ is contactable in $W_s$.

Note that $s_0$ is a critical value of $\ell$ on $\partial K$;
denote by $p_0$ the corresponding critical point.
By 2-convexity of $\RR^4\setminus K$,
the index of $p_0$ has to be at most $1$.
On the other hand, a disc that contracts $\alpha$ cannot be moved lower $s_0$.
Therefore the index of $p_0$ has to be at least $2$ --- a contradiction.
\qeds

\parit{Alexandrov-style proof.}
Assume that the complement to $K$ is two-convex.

Note that two-convexity is preserved under linear transformation.
Apply a linear transformation of $\RR^4$ that makes the coordinate planes $\Pi_1$ and $\Pi_2$ not orthogonal.

According to the main result in \cite{ABB}, $W\z=\RR^4\setminus (\Int K)$ has non-positive curvature in the sense of Alexandrov.
In particular, the universal metric covering $\tilde W$ of $W$ is a $\CAT(0)$ space.

By rescaling $\tilde W$ and passing to the limit we obtain that the universal Riemannian covering $Z$ of $\RR^4$ that branches in the planes $\Pi_1$ and $\Pi_2$ is a $\CAT(0)$ space.

Note that $Z$ is isometric to the Euclidean cone over universal covering $\Sigma$ of $\mathbb{S}^3$ branching in two great circles $\Gamma_i=\mathbb{S}^3\cap \Pi_i$ that are not orthogonal.
The shortest path in $\mathbb{S}^3$ between $\Gamma_1$ and $\Gamma_2$ traveled 4 times back and forth is shorter than $2\cdot\pi$ and it lifts to a closed geodesic in $\Sigma$.
It follows that $\Sigma$ is not $\CAT(1)$ and therefore $Z$ is not $\CAT(0)$ --- a contradiction.\qeds

The Morse-style proof is based on an idea of Mikhael Gromov \cite[see \S\textonehalf{} in][]{gromov-SGMC}, where two-convexity was introduced.

Note that the $1$-neighborhood of these two planes has two-convex complement $W$ in the sense of the second definition;
that is, if a closed curve $\gamma$ lies in the plane $\Pi$
and is contactable in $W$, then it is contactable in $\Pi\cap W$.
Clearly, the boundary of this neighborhood is not smooth
and as it follows from the problem, it cannot be smoothed in the class of two-convex sets. 

Two-convexity also shows up in comparison geometry --- the maximal open flat sets in the manifolds of nonnegative or nonpositive curvature are two-convex \cite{panov-petrunin}.


\parbf{Convex lens.}
Before going into the proof, let us describe a straightforward idea that does not work.

By the Gauss formula, we get that 
\[\int\limits_{D}k_1\cdot k_2\le\int\limits_{D}K,\] 
where $K$ denotes the intrinsic curvature of $D$.
Therefore it would be sufficient to show that the right-hand side is small;
however, the integral $\int_{D}K$ might be large for an arbitrarily small angle between the discs; for example, if $M=\mathbb{S}^3$ it might be arbitrarily close to $2\cdot \pi$.

\medskip

Denote by $\eps$ the maximal angle between the discs, we can assume that $\eps<\tfrac\pi2$.

Note that the function $h=\dist_{D'}$ is convex in $L$.
Moreover, the gradient $\nabla_xh$ points outside of $L$ for any $x\in D$. 

Consider the restriction $f=h|_D$.
Note that $f$ is a concave function that vanishes on $\partial D$.

Assume that $f$ is smooth.
Since the discs are meeting at angle at most $\eps<\tfrac\pi2$,
we have that $|\nabla f|\le \sin\eps$ and 
\[(\Hess f)(v,v)+ \cos\eps \cdot s(v,v)\le 0,\]
where $s$ denotes the second fundamental form of $D$ in $M$.
It follows that
\begin{align*}
k_1\cdot k_2&=\det s\le
\\
&\le \tfrac1{\cos^2\eps}\cdot \det(\Hess f)=
\\
&=\tfrac1{2\cdot\cos^2\eps}\cdot\left(|\Delta f|^2
-|\Hess f|^2\right).
\end{align*}

Applying the Bochner formula for $f$, we get that
\[\int\limits_D \left(|\Delta f|^2
-|\Hess f|^2
-K\cdot|\nabla f|^2\right)
=\int\limits_{\partial D}
\kappa\cdot|\nabla f|^2,\]
where $K$ and $\kappa$ denotes the curvature of $D$ and geodesic curvature of $\partial D$ in $D$ respectively.
By the Gauss--Bonnet formula, we get that
\[\int\limits_D 
K+\int\limits_{\partial D}\kappa=2\cdot\pi.\]
Therefore
\[\int\limits_Dk_1\cdot k_2\le \tfrac{\sin \eps}{\cos^2\eps}\cdot\pi.\]

If $f$ is not smooth, then one can smooth it using Greene--Wu construction \cite[Theorem~2]{greene-wu} and repeat the above argument for the obtained function.
\qedsf

This estimate was used by Nina Lebedeva and the author \cite{lebedeva-petrunin-curvature}.
For classical applications of Bochner's formula including the vanishing theorems and estimates for eigenvalues of Laplacian see \cite[][II \S 8]{lawson-michelsohn}.

\parbf{Small-twist.}
Given a positive integer $s$,
consider the Clifford torus
\[\TT^s_{\text{Cl}}=\tfrac1{\sqrt{s}}\cdot\underbrace{\SSS^1\times\dots\times\SSS^1}_\text{$s$ times}\subset\SSS^{2\cdot s-1}\subset\RR^{2\cdot s}.\]
Note that the normal curvatures of $\TT^s_{\text{Cl}}$ lie in the range from $1$ to $\sqrt{s}$, and $\TT^s_{\text{Cl}}$ comes with a flat metric.

Show that given a positive integer $n$, one can choose large $s=s(n)$
so that $\TT^s_{\text{Cl}}$ contains a geodesic $n$-dimensional
subtorus $\TT^n$ with all normal curvatures identically equal to $\sqrt{3\cdot\frac n{n+2}}$;
here we consider $\TT^n$ as a submanifold in $\RR^{2\cdot s}$.
(For example, $s(2)=3$; in this case, $\TT^2$ is a subtorus perpendicular to the main diagonal in $\TT^3_{\text{Cl}}$.) 

Now, choose a closed smooth manifold $M$.
By the Whitney embedding theorem, there is a smooth embedding $M\hookrightarrow\RR^n$ for $n>2\cdot \dim M$.
Applying rescaling we can assume that normal curvatures of this embedding are arbitrarily small;
we need it to be smaller than $\sqrt{3}-\sqrt{3\cdot\frac n{n+2}}$.
Composing this embedding with the natural length-preserving covering map $\RR^n\to \TT^n$, we get the needed immersion of $\iota\:M\to \RR^{2\cdot s}$.
\qeds

This construction was discovered by Mikhael Gromov \cite[2.A]{gromov2022}.
In fact, the bound $\sqrt{3\cdot\frac n{n+2}}$ is optimal \cite{petrunin2023}.
The immersion $\iota$ can be easily upgraded to embedding.
Applying the Nash embedding theorem instead of the Whitney embedding theorem,
one gets that the induced metric on $M\hookrightarrow\RR^n$ is proportional to any given Riemannian metric $g$ on $M$.

Here is a closely related problem.

\begin{pr}
Show that any $n$-dimensional submanifold in $\mathbb{S}^q$ with normal curvatures less than $\tfrac1{\sqrt{3}}$ is diffeomorphic to $\mathbb{S}^n$. 
\end{pr}

Note that Veronese embedding $\RP^2\hookrightarrow\mathbb{S}^5$ has normal curvatures exactly $\tfrac1{\sqrt{3}}$ in all directions.
Therefore the $\tfrac1{\sqrt{3}}$ bound is optimal;
see also \cite{petrunin-RPn}.

\csname @openrightfalse\endcsname
\chapter{Curvature-free differential geometry}
\chaptermark{Curvature-free}

The reader should be familiar 
with the notions of
smooth manifolds, 
Riemannian metrics,
and symplectic forms.

\subsection*{Distant involution}
\label{Distant involution}

\begin{pr}
Construct a Riemannian metric $g$ on $\mathbb{S}^3$ and an involution $\iota\:\mathbb{S}^3\z\to\mathbb{S}^3$ such that $\vol (\mathbb{S}^3,g)$ is arbitrarily small and 
\[|x\z-\iota(x)|_g>1\]
 for any $x\in\mathbb{S}^3$.
\end{pr}


\parit{Semisolution.}
Given $\eps>0$, construct a disk $\Delta$ in the plane with 
\[\length\partial \Delta<10\ \ \text{and}\ \ \area \Delta<\eps\]
that admits a continuous involution $\iota$ such that 
\[|\iota(x)-x|\ge 1\]
for any $x\in\partial \Delta$.

\begin{wrapfigure}{o}{20 mm}
\vskip-0mm
\centering
\includegraphics{mppics/pic-402}
\end{wrapfigure}

An example of $\Delta$ can be guessed from the picture;
the involution $\iota$ makes a length preserving half turn of its boundary $\partial \Delta$.

Take the product $\Delta\times \Delta\subset \RR^4$;
it is homeomorphic to the 4-ball.
Note that 
$$\vol_3[\partial(\Delta\times \Delta)]=2\cdot\area \Delta\cdot\length \partial \Delta<20\cdot\eps.$$
The boundary $\partial(\Delta\times \Delta)$ is homeomorphic to $\mathbb{S}^3$
and the restriction of the involution $(x,y)\z\mapsto (\iota(x),\iota(y))$ has the needed property.

All we have to do now is to smooth $\partial(\Delta\times \Delta)$ a little bit.
\qeds

This example is given by Christopher Croke \cite{croke}.
Note that according to Gromov's systolic inequality \cite{gromov-filling}, 
the involution $\iota$ above cannot be made isometric.

The following problem states that a similar construction is not possible for $\mathbb{S}^2$.

\subsection*{Another distant involution}
\label{Another distant involution}

\begin{pr}
Given $x\in \mathbb{S}^2$, denote by $x'$ its antipodal point.
Suppose that $g$ is a Riemannian metric on $\mathbb{S}^2$ such that
\[|x-x'|_g\ge1\]
for any $x\in\mathbb{S}^2$.
Show that the area of $(\mathbb{S}^2,g)$ is bounded below by a fixed positive constant. 
\end{pr}

The expected solution uses Besicovitch inequality described in the next problem.


\subsection*{Besicovitch inequality}
\label{Besicovitch inequality}

\begin{pr}
Let $g$ be a Riemannian metric on an $m$-dimensional cube $Q$ such that any curve connecting opposite faces has length at least $1$. 
Prove that 
\[\vol(Q,g)\ge 1,\] 
and the equality holds if and only if $(Q,g)$ is isometric to the unit cube.
\end{pr}

\subsection*{Minimal foliation\thm}
\label{gromomorphic-curves}

Minimal surfaces in Riemannian manifolds are defined on page \pageref{minimal surface}.

\begin{pr}
Consider the product of spheres $\mathbb{S}^2\times \mathbb{S}^2$ equipped with a Riemannian metric $g$ that is $C^\infty$-close to the product metric. 
Prove that there is a conformally equivalent metric $\lambda\cdot g$ and a re-parametrization of $\mathbb{S}^2\times \mathbb{S}^2$
such that for any $x,y\in \mathbb{S}^2$, the spheres $\{x\}\times\mathbb{S}^2$ and $\mathbb{S}^2\times \{y\}$ are minimal surfaces 
in $(\mathbb{S}^2\times \mathbb{S}^2,\lambda\cdot g)$.
\end{pr}

The expected solution requires pseudo-holomorphic curves introduced by Mikhael Gromov \cite{gromov-pseudoholomorphic}.

\subsection*{Volume and convexity\thm}
\label{Volume and convexity} 

A function $f$ defined on a Riemannian manifold is called convex if, for any geodesic $\gamma$, the composition $f\circ\gamma$ is a convex real-to-real function.

\begin{pr}
Let $M$ be a complete Riemannian manifold that admits a non-constant convex function. 
Prove that $M$ has infinite volume.
\end{pr}

The expected solution uses Liouville's theorem about phase volume.
It implies in particular, that the geodesic flow on the unit tangent bundle of a Riemannian manifold preserves the volume.

\subsection*{Sasaki metric}
\label{pr:Sasaki metric}

Let $(M,g)$ be a Riemannian manifold.
The Sasaki metric is a natural choice of Riemannian metric $\hat g$ on the total space of the tangent bundle $\tau\:\T M\to M$.
It is uniquely defined by the following properties:
\begin{itemize}
\item The map $\tau\:(\T M,\hat g)\to (M,g)$ is a Riemannian submersion.
\item The metric on each tangent space $\T_p\subset \T M$ is the Euclidean metric induced by $g$.
\item Assume that $\gamma(t)$ is a curve in $M$ and $v(t)\in\T_{\gamma(t)}$ is a parallel vector field along $\gamma$. 
Note that $v(t)$ forms a curve in $\T M$.
For the Sasaki metric, we have $v'(t)\perp \T_{\gamma(t)}$ for any $t$;
that is, the curve $v(t)$ normally crosses the tangent spaces $\T_{\gamma(t)}\subset \T M$.
\end{itemize}

In other words, we identify the tangent space 
$\T_u[\T M]$ for any $u\z\in \T_p M$ with the direct sum of vertical and horizontal subspaces $\T_p M\z\oplus \T_p M$.
The projection of this splitting is defined by the differential $d\tau\:\T\T M\to \T M$
and we assume that the velocity of a curve in $\T M$ formed by a parallel field along a curve in $M$ is horizontal.
Then $\T_u[\T M]$ is equipped with the metric $\hat g$ defined by
\[\hat g(X,Y)=g(X^V,Y^V)+g(X^H,Y^H),\]
where $X^V$ and $X^H\in\T_pM$ denote the vertical and horizontal components of $X\in\T_u[\T M]$.

\begin{pr}
Let $g$ be a Riemannian metric on the sphere $\mathbb{S}^2$.
Consider the tangent bundle $\T \mathbb{S}^2$ 
equipped with the induced Sasaki metric $\hat g$.
Show that
the space $(\T \mathbb{S}^2, \hat g)$ lies at a bounded distance to the ray $\RR_{\ge 0}=[0,\infty)$ in the sense of Gromov--Hausdorff.
\end{pr}

\subsection*{Two-systole}

\begin{pr} Given a large real number $L$,
construct a Riemannian metric $g$ on the 3-dimensional torus $\TT^3$ such that $\vol(\TT^3,g)=1$
and \[\area S\ge L\]
for any closed surface $S$ that does not bound in $\TT^3$.
\end{pr}

According to Gromov's systolic inequality \cite{gromov-filling}, the volume of $(\TT^3,g)$ can be bounded below in terms of its \emph{1-systole} defined to be the shortest length of a noncontractible closed curve in $(\TT^3,g)$.
The lower bound on the area of $S$ in the problem is called the 2-systole of $(\TT^3,g)$.

The problem implies that Gromov's systolic inequality does not have a direct 2-dimensional analog.

\subsection*{Normal exponential map\easy}
\label{Normal exponential map}
\label{page:Normal exponential map}

Let $(M,g)$ be a Riemannian manifold;
denote by $\T M$ the tangent bundle over $M$ and by $\T_p=\T_pM$ the tangent space at the point~$p$.

Given a vector $v\in\T_pM$, denote by $\gamma_v$ the geodesic in $(M,g)$
such that $\gamma(0)=p$ and $\gamma'(0)=v$.
The map $\exp\:\T M\to M$ defined by $v\mapsto \gamma_v(1)$ is called the exponential map.

The restriction of $\exp|_{\T_p}$ is called the \index{exponential map}\emph{exponential map at} $p$ and is denoted by $\exp_p$.

Given a smooth immersion $L\to M$,
denote by $\mathrm{N} L$ the normal bundle over $L$.
The restriction $\exp|_{\mathrm{N} L}$ is called the {}\emph{normal exponential map} of $L$ and is denoted by $\exp_L$.

\begin{pr}
Let $M$ be a complete connected Riemannian manifold
with an immersed complete connected Riemannian manifold $L$.
Show that the image  of the 
normal exponential map of $L$ is dense in $M$.
\end{pr}

\subsection*{Symplectic squeezing in the torus}
\label{Symplectic squeezing in the torus}

\begin{pr}
Let 
$\omega=dx_1\wedge dy_1+ dx_2\wedge dy_2$
be the standard symplectic form on $\RR^4$, and $\ZZ^2$ the integral lattice in the $(x_1,y_1)$ coordinate plane of $\RR^4$.

Show that an arbitrary bounded domain $\Omega\subset (\RR^4,\omega)$
admits a symplectic embedding into the quotient space $(\RR^4,\omega)/\ZZ^2$. 
\end{pr}

\subsection*{Diffeomorphism test\easy}
\label{Diffeomorphism test}

\begin{pr}
Let $M$ and $N$ be complete $m$-dimensional simply-connected Riemannian manifolds, and $f\:M\to N$ a smooth map such that 
$$|df(v)|\ge |v|$$
for any tangent vector $v$ of $M$.
Show that $f$ is a diffeomorphism.
\end{pr}

\subsection*{Volume of tubular neighborhoods\thm}
\label{Volume of tubular neighborhoods}

\begin{pr}
Let $M$ and $M'$ be isometric closed smooth submanifolds in a Euclidean space.
Show that for all small $r>0$ we have
$$\vol B(M,r)=\vol B(M',r),$$
where $B(M,r)$ denotes the $r$-neighborhood of $M$.
\end{pr}

\subsection*{Disk\hard}
\label{disk}

\begin{pr}
Given a large real number $L$,
construct a Riemannian metric $g$ on the disk $\mathbb D$ 
with 
\[\diam(\mathbb D,g)\le 1
\ \ 
\text{and}
\ \ 
\length_g \partial\mathbb D\le 1  \]
such that the boundary curve in $\mathbb D$ is not contractible in the class of closed curves with $g$-length less than $L$.
\end{pr}

\subsection*{Shortening homotopy}
\label{short-homotopy}

\begin{pr}
Let $M$ be a compact Riemannian manifold with diameter $D$ and $p\in M$.
Assume that for some $L>D$,
there are no geodesic loops based at $p$ in $M$
with length in the interval $(L-D,L+ D]$.
Show that for any path $\gamma_0$ in $(M,g)$ starting at $p$, 
there is a homotopy $\gamma_t$ relative to its endpoints
such that 
\begin{enumerate}[a)]
\item $\length \gamma_1<L$;
\item $\length \gamma_t\le \length \gamma_0+2\cdot D$ for any $t\in[0,1]$.
 
\end{enumerate}
\end{pr}

Examples of manifolds satisfying the above condition for some $L$ have been found among the Zoll spheres
by Florent Balachev, Christopher Croke, and Mikhail Katz \cite{balacheff-croke-katz}.


\subsection*{Convex hypersurface}
\label{Convex hypersurface}

Recall that a subset $K$ of a Riemannian manifold is called \index{convex set}\emph{convex} if every minimizing geodesic connecting two  points in $K$ lies completely in $K$. 

\begin{pr}
Let $M$ be a totally geodesic hypersurface 
in a closed Riemannian $m$-dimensional manifold $W$.
Assume that the injectivity radius of $M$ is at least $1$
and $M$ forms a convex set in $W$.

Show that the maximal distance from $M$ to the points of $W$ can be bounded below by a positive constant $\eps_m$ that depends only on the dimension $m$ (in fact, $\eps_m=\tfrac2{m+3}$ will~do).
\end{pr}

Note that we did not make any assumption on the injectivity radius of $W$.

\subsection*{Almost constant function}
\label{Almost constant function}

The unit tangent bundle $\UU M$ over a closed Riemannian manifold $M$
admits a natural choice of volume.
Let us equip $\UU M$ with the probability measure that is proportional to this volume.

We say that a unit-speed geodesic $\gamma\:\RR\to M$ is \emph{random}
if $\gamma'(0)$ takes a random value in $\UU M$.

\begin{pr}
Given $\eps>0$,
show that there is a positive integer $m$ such that
for any closed $m$-dimensional Riemannian manifold $M$
and any smooth $1$-Lipschitz function $f\:M\to\RR$ the following holds.

For a random unit-speed geodesic $\gamma$ in $M$ 
the event 
\[|f\circ\gamma(0)-f\circ\gamma(1)|>\eps\]
has probability at most $\eps$.
\end{pr}

\section*{Semisolutions}

\parbf{Another distant involution.}
Let $x\in \mathbb{S}^2$ be a point that minimizes the distance $|x-x'|_g$.
Consider a minimizing geodesic $\gamma$ from $x$ to $x'$.
We can assume that 
\[|x-x'|_g=\length \gamma=1.\]

Let $\gamma'$ be the antipodal arc to $\gamma$.
Note that $\gamma'$ intersects $\gamma$ only at the common endpoints $x$ and $x'$.
Indeed, if $p'=q$ for $p,q\in\gamma$, then $|p-q|\ge 1$.
Since $\length \gamma=1$, the points $p$ and $q$ must be the ends of $\gamma$.

It follows that $\gamma$ together with $\gamma'$ forms a closed simple curve in $\mathbb{S}^2$
that divides the sphere into two disks $D$ and $D'$.

Let us divide $\gamma$ into two equal arcs $\gamma_1$ and $\gamma_2$; each of length $\tfrac12$.
Suppose that $p,q\in\gamma_1$, then 
\begin{align*}
|p-q'|_g&\ge |q-q'|_g-|p-q|_g\ge
\\
&\ge 1-\tfrac12=\tfrac12.
\end{align*}
That is, the minimal distance from $\gamma_1$ to $\gamma_1'$ is at least~$\tfrac12$.
The same way we get that the minimal distance from $\gamma_2$ to $\gamma_2'$ is at least~$\tfrac12$.
By Besicovitch inequality, we get that 
\[\area(D,g)\ge \tfrac14\quad\text{and}\quad \area(D',g)\ge \tfrac14.\]
Therefore 
\[\area(\mathbb{S}^2,g)\ge\tfrac12.\]
\qedsf

This inequality was proved by Marcel Berger \cite{berger}. 
Christopher Croke conjectured that the optimal bound is $\tfrac4\pi$ and the round sphere is the only space that achieves this \cite[see Conjecture 0.3 in][]{croke}.

Let us indicate how to improve the obtained bound to
\[\area(\mathbb{S}^2,g)\ge1.\]

Suppose $x$, $x'$, $\gamma$, and $\gamma'$ are as above.
Consider the function
\[f(z)=\min_t \{\,|\gamma'(t)-z|_g+t\,\}.\]
Observe that $f$ is 1-Lipschitz.

Show that two points $\gamma'(c)$ and $\gamma(1-c)$ lie on one connected component of the level set $L_c=\set{z\in\mathbb{S}^2}{f(z)=c}$;
in particular 
\[\length L_c\ge 2\cdot|\gamma'(c)-\gamma(1-c)|_g.\]
By the triangle inequality, we have that
\begin{align*}
|\gamma'(c)-\gamma(1-c)|_g&\ge 1-|\gamma(c)-\gamma(1-c)|_g=
\\
&=1-|1-2\cdot c|.
\end{align*}

It remains to apply the coarea formula
\[\area(\mathbb{S}^2,g)\ge \int\limits_0^1\length L_c\cdot dc.\]

\parbf{Besicovitch inequality.}
Without loss of generality, we may assume that $Q=[0,1]^m$.
Set 
\[A_i=\set{(x_1,\dots,x_m)\in Q}{x_i=0}.\]

Consider the functions $f_i\:Q\to\RR$ defined by
$$f_i(x)=\min \{1,\dist_{A_i}(x)\}.$$
Note that each $f_i$ is $1$-Lipschitz; 
in particular, $|\nabla f_i|\le 1$ almost everywhere.

Consider the map
\[\bm{f}\:x\mapsto(f_1(x),\dots,f_m(x)).\]
Note that it maps $Q$ to itself
and, moreover, it maps each face of $Q$ to itself.
It follows that the restriction $\bm{f}|_{\partial Q}\:\partial Q\to \partial Q$ has degree one and therefore 
$\bm{f}\:Q\to Q$ is onto.

Let $h$ be the canonical metric on the cube $Q$.
Denote by $\mathrm{J}$ the Jacobian of the map $\bm{f}\:(Q,g)\to (Q,h)$.
Note that 
\[|\mathrm{J}(x)|=|\nabla_x f_1\wedge\dots\wedge\nabla_xf_m|\le 1.\]

By the area formula, we get 
\begin{align*}
\vol(Q,g)
&\ge \int\limits_{x\in Q}|\mathrm{J}(x)|\ge
\\
&\ge\vol(Q,h)=\\
&=1.
\end{align*}

In the case of equality, we have that $\<\nabla_x f_i,\nabla_x f_j\>=0$ for $i\ne j$ 
and $|\nabla_x f_i|=1$ for almost all $x$.
It follows then that the map 
\[\bm{f}\:(Q,g)\z\to (Q,h)\] 
is an isometry.
\qeds

This inequality was proved by Abram Besicovitch \cite{besicovitch}.
It has a number of applications in Riemannian geometry.
For example, using this inequality it is easy to solve the following problem.

\begin{pr}
Assume a metric $g$ on $\RR^m$ coincides with the Euclidean metric outside of a bounded set $K$;
assume further that any geodesic that enters $K$ exits $K$ the same way the Euclidean geodesic would have done. 
Show that $g$ is flat.
\end{pr}

There is a weaker version of the Besicovitch inequality that works for the Hausdorff measure for any metric on the cube; nearly the same proof works.
Here is one of its applications suggested by Stephan Stadler:

\begin{pr}
Let $X$ be a contractible metric space with zero $(n+1)$-dimensional Hausdorff measure.
Assume that $\Delta_1,\Delta_2\subset X$ are two embedded $n$-disks having the same boundary.
Show that $\Delta_1=\Delta_2$.
\end{pr}

\parbf{Minimal foliation.}
The proof is based on the observation that a self-dual harmonic 2-form on $(\mathbb{S}^2\times\mathbb{S}^2,g)$
without zeros defines a symplectic structure.

\medskip

Note that there is a self-dual harmonic 2-form on $(\mathbb{S}^2\times\mathbb{S}^2,g)$;
that is, a 2-form $\omega$ such that $d\omega=0$ and $\mathop{\star}\omega=\omega$,
where $\mathop{\star}$ is the Hodge star operator.
Indeed, take a generic harmonic form $\phi$.
Note that the form $\mathop{\star}\phi$ is also harmonic.
Since $\mathop{\star}(\mathop{\star}\phi)=\phi$,
the form $\omega=\phi+\mathop{\star}\phi$ does the job.

Choose $p\in \mathbb{S}^2\times\mathbb{S}^2$.
We can use $g_p$ to identify the tangent space $\T_p$ and the cotangent space $\T^*_p$.
There is a $g_p$-orthonormal basis $e_1, e_2, e_3, e_4$ on $\T_p$ such that 
\[\omega_p=\lambda_p\cdot e_1\wedge e_2+\lambda'_p\cdot  e_3\wedge e_4.\]
Note that 
\[\mathop{\star}\omega_p=\lambda'_p\cdot e_1\wedge e_2+\lambda_p\cdot  e_3\wedge e_4.\]
Since $\mathop{\star}\omega_p=\omega_p$, we have $\lambda_p=\lambda'_p$.

Consider the rotation $\mathrm{J}_p\:\T_p\to\T_p$ defined by 
\[ 
e_1\mapsto -e_2,
\quad 
e_2\mapsto e_1,
\quad 
e_3\mapsto -e_4,
\quad 
e_4\mapsto e_3.\]
Note that
\[\mathrm{J}_p\circ\mathrm{J}_p =-\id
\quad\text{and}\quad
\omega(X,Y)=\lambda_p\cdot g(X,\mathrm{J}_pY)\] 
for any two tangent vectors $X,Y\in \T_p$.

Consider the canonical symplectic form $\omega_0$ on $\mathbb{S}^2\times\mathbb{S}^2$ which is the sum of the pullbacks of the volume form on $\mathbb{S}^2$  
by the two coordinate projections $\mathbb{S}^2\times\mathbb{S}^2\to \mathbb{S}^2$.
Note that for the canonical metric on $\mathbb{S}^2\times\mathbb{S}^2$,
the form $\omega_0$ is harmonic and self-dual. 
Since $g$ is close to the standard metric,
we can assume that $\omega$ is close to $\omega_0$.
In particular, $\lambda_p\ne0$ for any $p\in \mathbb{S}^2\times\mathbb{S}^2$.

It follows that $\mathrm{J}$ is a pseudo-complex structure for the symplectic form $\omega$ on $\mathbb{S}^2\times\mathbb{S}^2$.
The Riemannian metric $g'=\lambda\cdot g$ is conformal to $g$ and $\omega(X,Y)=g'(X,\mathrm{J} Y)$ 
for any two tangent vectors $X,Y$ at one point.
In this case, the $\mathrm{J}$-holomorphic curves are minimal with respect to $g'$;
in fact, each of them is area-minimizing in its homology class. 

It remains to reparametrize $\mathbb{S}^2\times \mathbb{S}^2$
so that vertical and horizontal spheres would form pseudo-holomorphic curves in the homology classes of $x\times \mathbb{S}^2$ and $\mathbb{S}^2\times y$.
\qeds

For general metrics, the form $\omega$ might vanish at some points.
If the metric is generic, then it happens on disjoint circles \cite{honda}.

\parbf{Volume and convexity.}
We use the idea from the proof of the Poincaré recurrence theorem.

\medskip

Let $M$ be a complete Riemannian manifold that admits a convex function $f$.
Denote by $\tau\:\UU M\to M$ the unit tangent bundle over $M$. 
Consider the function $F\:\UU M\to \RR$ defined by $F(u)=f\circ\tau(u)$.

Note that 
there is a nonempty bounded open set $\Omega\subset \UU M$
such that $df(u)>\eps$ for any $u\in \Omega$ and some fixed $\eps>0$.

Denote by $\phi^t$ the geodesic flow for time $t$ on $\UU M$.
By Liouville's theorem about phase volume, we have
\[\vol[\phi^t(\Omega)]=\vol\Omega\leqno{({*})}\] 
for any $t$.

Given $u\in \UU M$,
consider the function 
$h_u(t)=F\circ\phi^t(u)$.
Since $f$ is convex, so is $h_u$.
Therefore $h'_u(t)>\eps$ for any $t\ge 0$ and $u\in\Omega$.

Since $\Omega$ is a bounded set, the set of values $F(\Omega)$ is bounded as well.
It follows that there is an infinite sequence of time moments 
\[0=t_0<t_1\z<t_2<\dots\]
such that 
\[h_v(t_{i-1})<h_u(t_{i})\]
for any $u,v\in \Omega$ and $i$.
In particular, we have
$$\phi^{t_i}(\Omega)\cap\phi^{t_j}(\Omega)=\emptyset$$ 
for $i\ne j$.
By $({*})$, the latter implies that $\vol (\UU M)=\infty$.
Hence 
\[\vol M=\infty.\qedsin\]
\medskip

The problem is due 
to Richard Bishop and Barrett O'Neill \cite{bishop-oneill};
it was generalized by
Shing-Tung Yau \cite{yau}.

\parbf{Sasaki metric.}
Choose a point $p\in\mathbb{S}^2$.
Note that any rotation of the tangent space $\T_p=\T_p(\mathbb{S}^2,g)$
appears as a holonomy of a loop at $p$;
moreover, the length of such loop can be bounded by a constant, say~$\ell$.

Indeed, fix a smooth homotopy $\gamma_t\:[0,1]\to \mathbb{S}^2$, $t\in[0,1]$ of loops based at $p$ that sweeps out $\mathbb{S}^2$.
By the Gauss--Bonnet formula, the total curvature of $(\mathbb{S}^2,g)$ is $4\cdot\pi$.
Therefore any rotation of $\T_p$ appears as the holonomy of $\gamma_t$ for some $t$ and we can take 
\[\ell=\max\set{\length\gamma_t}{t\in[0,1]}.\]

Denote by $d$ the diameter of $(\mathbb{S}^2,g)$.
From the above, it follows that for any two unit tangent vectors $v\in\T_p$ 
and $w\in T_q$
there is a path 
$\gamma\:[0,1]\to\mathbb{S}^2$ from $p$ to $q$
such that 
\[\length \gamma\le \ell+d\] 
and
$w$ is the parallel transport of $v$ along $\gamma$.

In particular, the diameter of the set of all vectors of fixed magnitude in $(\T \mathbb{S}^2,\hat g)$ has diameter at most $\ell+d$.
Therefore the map $\T\mathbb{S}^2\to[0,\infty)$ defined by $v\mapsto |v|$ 
preserves the distance up to an error of $\ell+d$.
Hence the result follows.
\qeds


\parbf{Two-systole.}
Consider the unit cube with three not intersecting cylindrical tunnels  
between the pairs of opposite faces.
In each tunnel, shrink the metric long-wise and expand it  cross-wise while keeping the volume the same.

\medskip

\begin{wrapfigure}{o}{33 mm}
\vskip-0 mm
\centering
\includegraphics{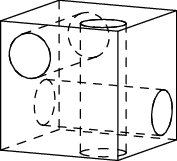}
\vskip2 mm
\end{wrapfigure}

More precisely, assume $(x,y,z)$ is the coordinate system on the cylindrical tunnel $\DD\z\times [0,1]$. 
Then the new metric is defined by
\[g=\phi\cdot [(dx)^2+ (dy)^2]+\tfrac1{\phi^2}\cdot (dz)^2,\]
where $\phi=\phi(x,y)$ is a positive smooth function on $\DD$ taking huge values around the center and equal to 1 near the boundary of $\DD$.

Gluing the opposite faces of the cube, we obtain a 3-dimensional torus with a smooth Riemannian metric.

Since the surface $S$ does not bound in $\TT^3=\mathbb{S}^1\times\mathbb{S}^1\times\mathbb{S}^1$,
one of the three coordinate projections $\TT^3\to\TT^2=\mathbb{S}^1\times\mathbb{S}^1$
induces a map of a non-zero degree $S\to\TT^2$.
It follows that 
\[\area S\ge  \area(\DD,\phi\cdot [(dx)^2+ (dy)^2]).\]
For the right choice of the function $\phi$, the right-hand side can be made larger than the given number $L$.
Hence the statement follows.
\qeds

I learned this problem from Dmitri Burago.
The three-tunnel construction was invented by Gustav Hedlund in a different context \cite{hedlund,bangert1990}.
The following problem of Larry Guth \cite{guth} is closely related:

\begin{pr}
Given a small $\eps>0$, construct a bi-Lipschitz area-nonincreasing degree-one map 
\[[0,1]\times[0,1]\times [0,\eps]\to [0,\tfrac\eps7]\times [0,\tfrac\eps7]\times [0,\tfrac1{7\cdot\eps}].\]

\end{pr}

\parbf{Normal exponential map.}
Assume there are $p\in M$ and $\eps>0$ 
such that the image of the normal exponential map to $L$
 does not intersect the ball $B(p,\eps)_M$; that is, no geodesic normal to $L$ crosses the ball.

Choose a positive real number $R$ such that $B(p,R)_M\cap L\ne \emptyset$.
The sectional curvature of $M$ in the ball $B(p,R)$
is bounded below by a constant, say $K$.

Given $q\in L$, denote by $v_q \in \T_qM$ the direction of a minimizing geodesic $[qp]$.
Note that $v_q\notin \mathrm{N}_qL$.
Moreover, there is $\delta=\delta(\eps,K,R)\z>0$ 
such that for any point $q\in B(p,R)_M\cap L$,
and any normal vector $n\in \mathrm{N}_qL$,
we have 
\[\measuredangle (v_q,n)>\delta.\]
Otherwise, the geodesic in the direction of $n$ would cross $B(p,\eps)_M$.

It follows that starting at any point $q\in B(p,R)_M\cap L$ 
one can construct a unit-speed curve $\gamma$ in $L$ such that 
\[|p-\gamma(t)|\le |p-q|-t\cdot\sin \delta.\]
Following $\gamma$ for some time brings us to $p$;
that is, $p\in L$ --- a contradiction.
\qeds

{

\begin{wrapfigure}{o}{28 mm}
\vskip-0 mm
\centering
\includegraphics{mppics/pic-404}
\end{wrapfigure}

The problem was suggested by Alexander Lytchak.

On the diagram, you see an example of an immersion 
such that one point does not lie in the image of the corresponding normal exponential map.
It might be interesting to understand what type of subsets can be avoided by such images.

}
\parbf{Symplectic squeezing in the torus.}
The embedding will be given as a composition of a linear symplectomorphism $\lambda$ 
with the quotient map $\phi\:\RR^4\to \TT^2\times\RR^2$ by the integral $(x_1,y_1)$-lattice.

\medskip

The composition $\phi\circ\lambda$ will preserve the symplectic structure;
it remains to find $\lambda$ such that the restriction $\phi\circ\lambda|_\Omega$
is injective.

Without loss of generality,
we can assume that $\Omega$ is a ball centered at the origin.
Choose an oriented 2-dimensional subspace $V$ of $\RR^4$ 
such that the integral of $\omega$ over 
$\Omega\cap V$ is a  positive number smaller than $\tfrac\pi4$. 

Note that there is a linear symplectomorphism $\lambda$ that maps planes parallel to $V$ to planes parallel to the $(x_1,y_1)$-plane, and that maps the disk $V\cap\Omega$ to a round disk.
It follows that the intersection of $\lambda(\Omega)$ 
with any plane parallel to the $(x_1,y_1)$-plane is a disk of radius at most $\tfrac 12$.
In particular $\phi\circ\lambda|_\Omega$
is injective.\qeds

This construction was given 
by Larry Guth \cite{guth-symplectic}
and attributed to Leonid Polterovich.

Note that according to Gromov's non-squeezing theorem \cite{gromov-pseudoholomorphic}, 
an analogous statement with $\CC\times \DD$ as the target space does not hold;
here $\DD\subset \CC$ is the open unit disk with the induced symplectic structure.
In particular, it shows that
the projection of $\lambda(\Omega)$ as above 
to the $(x_1,y_1)$-plane
cannot be made arbitrarily small.

\parbf{Diffeomorphism test.}
Note that the map $f$ is an open immersion.

Let $h$ be the pullback metric on $M$ for $f\:M\to N$.
Clearly, $h\ge g$.
In particular, $(M,h)$ is complete and the map $f\:(M,h)\to N$ is a local isometry. 

Note that any local isometry between complete connected Riemannian manifolds of the same dimension is a covering map.
Since $N$ is simply-connected, the result follows.
\qeds 


\parbf{Volume of tubular neighborhoods.}
This problem is a direct corollary of the so-called \emph{tube formula} given by Hermann Weyl \cite{weyl}.
It expresses the volume of the $r$-neighborhood of $M$ as a polynomial $p(r)$;
the coefficients of $p$, up to a multiplicative constant, are integrals over~$M$ of some quantities called the \emph{Lipschitz--Killing curvatures} --- these are certain scalars that can be expressed in terms of the curvature tensor at the given point.
The proof is done by straightforward calculations.


\parbf{Disk.}
The following claim is the key step in the proof.
\begin{cl}{$({*})$}
Given a positive integer $n$, there is a binary tree $T_n$ embedded into the disk $\DD$ such that any null-homotopy of $\partial \DD$ passes a curve that intersects $n$ different edges.
\end{cl}

The proof of the claim can be done by induction on $n$; the base is trivial.
Assuming we constructed $T_{n-1}$, the tree $T_n$ can be obtained by identifying three endpoints of three copies of $T_{n-1}$.

\begin{figure}[h!]
\vskip0mm
\centering
\includegraphics{mppics/pic-406}
\end{figure}

Take $\eps=\tfrac1{10}$ and fix a large integer $n$.
Let us construct a metric on $\DD$ with the embedded tree $T_n$ as in $({*})$ such that
its diameter and the length of its boundary are less than $1$
and  
the distance between any two edges of $T_n$ without a common vertex 
is at least $\eps$.

Choose a Riemannian metric $g$ on the cylinder $\mathbb S^1\times [0,1]$ such that
\begin{itemize}
\item The $\eps$-neighborhoods of the boundary components 
have product metrics.
\item Any vertical segment $x\times[0,1]$ has length $\tfrac 12$.
\item One of the boundary components has length $\eps$.
\item The other boundary component has length $2\cdot m\cdot \eps$, 
where $m$ is the number of edges in the tree $T_n$.
\end{itemize}
Equip $T_n$ with a length-metric so that each edge has length $\eps$.
Glue the cylinder $(\mathbb S^1\times [0,1],g)$ along its long boundary component to the tree $T_n$ by a piecewise isometry 
in such a way that the resulting space is homeomorphic to a disk and the obtained embedding of $T_n$ in $\DD$ is the same as in the claim $({*})$.

By $({*})$, any null-homotopy of the boundary passes a curve that intersects $n$ different edges of $T_n$.
By construction this curve is longer than $\tfrac{\eps}{10}\cdot n=\tfrac{1}{100}\cdot n$.

The obtained metric is not Riemannian, but it is easy to smooth while keeping this property.
Since $n$ is large the result follows.
\qeds
 
This example was constructed by Sidney Frankel and Mikhail Katz \cite{frankel-katz}.

\parbf{Shortening homotopy.}
Set 
\[p=\gamma_0(0)\ \ \text{and}\ \  \ell_0=\length\gamma_0.\]

By a compactness argument,
there exists $\delta>0$ 
such that no geodesic loop based at $p$ has length in the interval $(L-D, L+D+\delta]$. 

Assume that $\ell_0\ge L+\delta$.
Choose $t_0\in [0,1]$ such that
\[\length\left(\gamma_0|_{[0,t_0]}\right)=L+\delta.\]
Let $\sigma$ be a minimizing geodesic from $\gamma(t_0)$
to $p$.
Note that $\gamma_0$ is homotopic to the concatenation 
\[\gamma_0'=\gamma_0|_{[0,t_0]}*\sigma*\bar\sigma*\gamma|_{[t_0,1]},\]
where $\bar\sigma$ denotes the backward parametrization of $\sigma$.

Applying a curve shortening process to the loop $\lambda_0=\gamma|_{[0,t_0]}*\sigma$, 
we get a  homotopy $\lambda_t$
relative to its endpoints 
from the loop $\lambda_0$ to a geodesic loop $\lambda_1$ at $p$.
From the above, 
\[\length\lambda_1\le L-D.\]

\begin{wrapfigure}{o}{31 mm}
\vskip-6mm
\centering
\includegraphics{mppics/pic-408}
\end{wrapfigure}

The concatenation $\gamma_t=\lambda_t*\bar\sigma*\gamma|_{[t_0,1]}$
is a homotopy
from $\gamma_0'$ to another curve $\gamma_1$.
From the construction, it is clear that 
\begin{align*}
 \length \gamma_t&\le \length \gamma_0+2\cdot \length\sigma\le
 \\
 &\le \length \gamma_0+2\cdot D
\end{align*}
for any $t\in[0,1]$
and 
\begin{align*}
 \length \gamma_1&=\length\lambda_1+\length\sigma+\length\gamma|_{[t_0,1]}\le
\\ &\le L-D+D+\length\gamma-(L+\delta)=
\\ &=\ell_0 -\delta.
\end{align*}

Repeating the procedure a sufficient number of times, we get curves $\gamma_2,\dots,\gamma_n$
connected by the needed homotopies so that 
$\ell_{i+1}\le\ell_i-\delta$ and $\ell_n\z< L+\delta$,
where $\ell_i=\length\gamma_i$.

If $\ell_n\le L$, we are done.
Otherwise, repeat the argument again for $\delta'=\ell_n-L$.
\qeds

The problem is due to 
Alexander Nabutovsky 
and Regina Rotman \cite{nabutovsky-rotman}.

\parbf{Convex hypersurface.}
First, let us define the {}\emph{cone construction} of maps into $M$.

Let $\Delta'$ be a simplex 
with a vertex $v$.
Denote by $\Delta$ the facet in $\Delta'$ opposite to $v$.
Let $f\:\Delta\to M$ be a map such that $f(\Delta)\subset B(x,1)_M$ for some $x \in M$.
Given $w\in \Delta$, let $\gamma_w\:[0,1]\to M$ be the minimizing geodesic path from $x$ to  $f(w)$.
Since the injectivity radius of $M$ is at least $1$, the path $\gamma_w$ is uniquely defined.
The map $f'\:\Delta'\to M$ defined as 
\[f'\:(1-t)\cdot v+t\cdot w\mapsto \gamma_w(t)\] 
is called the {}\emph{cone over $f$} with vertex $x$. 

One may start with a map $f_0\:\Delta_0\to M$ and iterate the cone construction for the vertices $x_1,\dots x_k$,
to get a sequence of maps $f_i\:\Delta_i\to M$
as long as $f_{i-1}(\Delta_{i-1})\subset B(x_i,1)$.
A straightforward application of the triangle inequality 
shows that the latter conditions hold if 
$f_0(\Delta_0)\subset B(x_i,s)$ for each $i$ and $s<\tfrac2{2+k}$.

\medskip

Now we go back to the solution.

Choose a fine triangulation of $W$ so that $M$ becomes a sub-complex of $W$.
We can assume that the diameter of each simplex in $\tau$ is less than any given
$\eps>0$.
Furthermore, we can assume that all the vertices of $\tau$ can be colored with $m+2$ colors $(0,\dots, m+1)$
in such a way that the vertices of each simplex 
have different colors;
the latter can be achieved by passing to the barycentric subdivision of $\tau$.
Denote by $\tau_i$ the maximal $i$-dimensional sub-complex of $\tau$ 
with all the vertices colored by $0,\dots, i$.

Let $h$ be the maximal distance from points in $W$ to $M$.
For each vertex $v$ in $\tau$ 
choose a point $v'\in M$ at distance $\le h$.
Note that 
if $v$ and $w$ are vertices of one simplex,
then
\[|v'-w'|_M<2\cdot h+\eps.\]

Assume that $\tfrac{2}{m+3}>h$.
Choose positive $\eps<\tfrac{2}{m+3}-h$ and use it in the construction of the triangulation $\tau$ above.
Applying the iterated cone construction for each simplex of $\tau$
we get an extension of the map $v\mapsto v'$ defined on $\tau_0$ to $\tau_1,\dots\tau_{m+1}$.
According to the above estimates, the cone constructions are defined at each of the needed $m+1$ iterations.

This way we get 
to a retraction $W\to M$.
It follows that the fundamental class of $M$ vanishes in the homology ring of $M$ --- 
a contradiction. 
\qeds

This problem is a stripped version of the bound on filling radius given by Mikhael Gromov \cite{gromov-filling}.  

\parbf{Almost constant function.}
Given a positive integer $m$,
denote by $\delta_m$ 
the expected value of $|x_1|$ for the random unit vector 
$\bm{x}\z=(x_1,\dots,x_m)\z\in\mathbb{S}^{m-1}$ 
with respect to the uniform distribution.

Observe that $\delta_m\to 0$ as $m\to\infty$.
Indeed, from symmetry and Bunyakovsky inequality we get
\[
\tfrac1m=\tfrac1m\cdot\mathrm{E}(|\bm{x}|^2)
=\mathrm{E}(x_1^2)\ge \mathrm{E}(|x_1|)^2=\delta_m^2.
\]

Since $f$ is $1$-Lipschitz,
\[\mathrm{E}(|df(w)|)\le\delta_m\]
for a random vector $w$ in $\UU M$.

Note that 
\begin{align*}
|f\circ \gamma(1)-f\circ\gamma(0)|
&=
\left|\int\limits_0^1df(\gamma'(t))\cdot dt\right|\le \\
&\le \int\limits_0^1\left|df(\gamma'(t))\right|\cdot dt.
\end{align*}

Assume that $\gamma'(0)$ takes a random value in $\UU M$.
By Liouville's theorem about phase volume, the same holds for $\gamma'(t)$
for any fixed $t$.
Therefore
\begin{align*}
\mathrm{E}(|f\circ \gamma(1)-f\circ\gamma(0)|)\le \mathrm{E}\left(\int\limits_0^1|df(\gamma'(t))|\cdot dt\right)\le\delta_m.
\end{align*}

By Markov's inequality,
the probability of the event 
\[|f\circ \gamma(1)-f\circ\gamma(0)|>\eps\]
is at most $\tfrac{\delta_m}{\eps}$.
Hence the result follows.
\qeds

I learned this problem from Mikhael Gromov.
It gives an example in the Riemannian world
of the so-called 
\index{concentration of measure}\emph{concentration of measure phenomenon}
\cite{milman-schechtman,ledoux}.

\csname @openrightfalse\endcsname
\chapter{Metric geometry}

In this chapter, we consider metric spaces.
The relevant background material can be found in \cite{bbi} or \cite{petrunin2021pure}. 

Let us fix a few standard notations.
\begin{itemize}
\item The distance between two points $x$ and $y$ in a metric space $X$
will be denoted by 
\[\dist_x(y),
\quad
|x-y|
\quad\text{or}\quad
|x-y|_X,
\]
the latter notation is used to emphasize that $x$ and $y$ belong to the space $X$.
\item A metric space $X$ is called {}\emph{length-metric space} if, for any two points $x,y\in X$ and any $\eps>0$, the points $x$ and $y$ can be connected by a curve $\alpha$
with
\[\length\alpha<|x-y|_X+\eps.\]
In this case, we say the metric on $X$ is a \index{length-metric}\emph{length-metric}.
\end{itemize}

\subsection*{Embedding of a compact}
\label{compact} 

\begin{pr}
Prove that any compact metric space 
is isometric to 
a subset of a compact length-metric space.
\end{pr}

\parit{Semisolution.}
Let $K$ be a compact metric space.
Denote by $\mathcal{B}(K,\RR)$ the space of real-valued bounded functions on $K$
equipped with sup-norm; 
that is, 
\[|f|=\sup\set{|f(x)|}{x\in K}.\]

Note that the map $K\to \mathcal{B}(K,\RR)$, defined by $x\mapsto \dist_x$
is a distance-preserving embedding.
Indeed, by the triangle inequality we have
\[|\dist_x(z)-\dist_y(z)|\le |x-y|_K\]
for any $z\in K$
and the equality holds for $z=x$.

In other words, we can and will consider $K$ as a subspace of $\mathcal{B}(K,\RR)$.

Denote by $W$ the linear convex hull of $K$ in $\mathcal{B}(K,\RR)$;
that is, $W$ is the intersection of all closed convex subsets containing $K$. 
Clearly, $W$ is a complete subspace of $\mathcal{B}(K,\RR)$.

Since $K$ is compact we can choose a finite $\eps$-net $K_\eps$ in $K$.
The set $K_\eps$ lies in a finite-dimensional subspace;
therefore its convex hull $W_\eps$ is compact.
Note that $W$ lies in the $\eps$-neighborhood of $W_\eps$.
Therefore, $W$ admits a compact $\eps$-net for any $\eps>0$.
That is, $W$ is totally bounded and complete and therefore compact.

Note that line segments in $W$ are geodesics for the metric induced by the sup-norm. 
In particular, $W$ is a compact length-metric space as required.
\qeds

The map $x\mapsto \dist_x$ is called the \index{Kuratowski embedding}\emph{Kuratowski embedding},
it was constructed in \cite{kuratowski}.
Essentially the same map 
was described by Maurice Fr\'echet \cite[][this is the paper where metric spaces were introduced]{frechet}.

The problem also follows directly from a theorem of John Isbell, stating that \emph{injective envelopes} of compact metric spaces are compact;
the injective envelope is an analog of convex hull in the category of metric spaces
\cite[see 2.11 in][]{isbell}.

The following related problem is open even for three-point sets.
This problem was mentioned by Mikhael Gromov \cite[6.B$_1$(f)]{gromov-asymptotic}.

\begin{pr}
Is it true that any compact subset of a complete $\CAT(0)$ length-space lies in a compact convex set?
\end{pr}

\subsection*{Non-contracting map\easy}
\label{Noncontracting map}

A map $f\: X\to Y$ between metric spaces is called \index{non-contracting map}\emph{distance non-contracting} if
\[|f(x)-f(x')|_Y\ge |x-x'|_X\]
for any two points $x,x'\in X$.

\begin{pr}
Let $K$  be a compact metric space and
\[f\:K\z\to K\] 
a distance non-contracting map.
Prove that $f$ is an isometry.
\end{pr}

\subsection*{Finite-whole extension}
\label{Finite-whole extension}

A map $f\: X\to Y$ between metric spaces is called \index{non-expanding map}\emph{non-expanding} if
\[|f(x)-f(x')|_Y\le |x-x'|_X\]
for any two points $x,x'\in X$.

\begin{pr}
Let $X$ and $Y$ be metric spaces, 
$Y$ compact,
$A\subset X$,
and $f\:A\to Y$ a non-expanding map.
Assume that for any finite set $F\subset X$ there is a non-expanding map $F \to Y$
that agrees with $f$ in $F\cap A$.
Show that there is a non-expanding map $X\to  Y$ that agrees with $f$ on $A$.
\end{pr}

\subsection*{Horo-compactification\easy}
\label{Horocompactification}

Let $X$ be a metric space.
Denote by $C(X,\RR)$ the space of continuous functions $X\to \RR$
equipped with the \emph{compact-open topology};
that is, for any compact set $K\subset X$ and any open set $U\subset \RR$
the set of all continuous functions $f\: X\to \RR$ such that $f(K)\subset U$
is declared to be open.

Choose a point $x_0\in X$.
Given a point $z\in X$, let $f_z\in C(X,\RR)$ be the function defined by 
\[f_z(x)=\dist_z(x)-\dist_z(x_0).\]
Let $F_X\:X\to C(X,\RR)$ be the map 
sending $z$ to $f_z$.

Denote by $\bar X$ 
the closure of $F_X(X)$ in $C(X,\RR)$;
note that $\bar X$ is compact.
That is, 
if $F_X$ is an embedding, 
then $\bar X$ is a compactification of $X$,
which is called the \index{horo-compactification}\emph{horo-compactification}.
In this case, the complement 
$\partial_\infty X\z=\bar X\setminus F_X(X)$ 
is called the {}\emph{horo-absolute} of $X$.

The construction above is due to Mikhael Gromov \cite{gromov-hyperbolic}.

\begin{pr}
Construct a proper metric space $X$
such that 
\[F_X\: X \to C(X,\RR)\] 
is not an embedding.
Show that there are no such examples among proper length-metric spaces.
\end{pr}

\subsection*{Approximation of the ball by a sphere}
\label{3-sphere is close to a ball}

\begin{pr}
Construct a sequence of Riemannian metrics on $\mathbb{S}^3$ converging in the sense of Gromov--Hausdorff 
to the unit ball in $\RR^3$.
\end{pr}

\subsection*{Macroscopic dimension\easy}
\label{macroscopic dimension} 

Let $X$ be a locally compact metric space
and $a>0$.

Following Mikhael Gromov \cite{gromov:macroscopic-dimension},
we say that the \index{macroscopic dimension}\emph{macroscopic dimension}  of $X$ at scale $a$ is the smallest integer $m$ such that there is a continuous map $f$ from $X$ to an $m$-dimensional simplicial complex $K$
with 
\[\diam[f^{-1}\{k\}]<a\]
for any point $k\in K$.

Equivalently, the macroscopic dimension of $X$ on the scale $a$ can be defined as 
the smallest integer $m$ such that $X$ admits an open cover with diameter of each set less than $a$ 
and such that each point in $X$ is covered by at most $m+1$ sets in the cover.

\begin{pr}
Let $M$ be a simply-connected Riemannian manifold with the following property: 
every closed curve is null-homotopic 
in its own  1-neighborhood. 
Prove that the macroscopic dimension of $M$ at scale $100$ is at most $1$.
\end{pr}

\subsection*{No Lipschitz embedding\hard}
\label{weird-metric} 

\begin{pr}
Construct a length-metric $d$ on $\RR^3$ such that the space $(\RR^3,d)$ does not admit a locally Lipschitz embedding into the 3-dimensional Euclidean space.
\end{pr}

\subsection*{Sub-Riemannian sphere\thm}
\label{sub-Riemannian} 

Let us explain what is a sub-Riemannian metric.

Let $(M,g)$ be a Riemannian manifold.
Assume that in the tangent bundle $\T M$ 
a choice of sub-bundle $H$ is given.

Let us call the sub-bundle $H$  \index{horizontal distribution}\emph{horizontal distribution}.
The tangent vectors in $H$ will be called {}\emph{horizontal}.
A piecewise smooth curve will be called {}\emph{horizontal}
if all its tangent vectors are horizontal.

The sub-Riemannian distance between any two points $x$ and $y$ is defined as the infimum of lengths of horizontal curves connecting $x$ to~$y$.

Alternatively, the distance can be defined as the limit of Riemannian distances 
for the metrics 
\[g_\lambda(X,Y)=g(X^H,Y^H)+\lambda\cdot g(X^V,Y^V)\] 
as $\lambda\to \infty$,
where $X^H$ denotes the horizontal part of $X$;
that is, the orthogonal projection of $X$ to $H$
and $X^V$ denotes the vertical part of~$X$;
so, $X^V+X^H=X$.

We also need an additional condition to ensure the following properties 
\begin{itemize}
\item The sub-Riemannian metric induces the original topology on the manifold. 
In particular, if $M$ is connected, then the distance cannot take infinite values.
\item Any curve in $M$ can be arbitrarily well approximated by a horizontal curve with the same endpoints.
\end{itemize}
The most common condition of this type is the so-called {}\emph{complete non-integrability};
it means that for any $x\in M$, 
one can choose a basis in its tangent space $\T_xM$
from the vectors of the following type
\[A(x),\quad  [A,B](x),\quad [A,[B,C]](x),\quad [A,[B,[C,D]]](x),\dots\] 
where $[{*},{*}]$ denotes the Lie bracket 
and the vector fields $A,B,C,D, \dots$ are horizontal.

\begin{pr}
Prove that any sub-Riemannian metric 
on $\mathbb{S}^m$ is isometric to the intrinsic metric of a hypersurface in $\RR^{m+1}$.
\end{pr}

It will be difficult to solve the problem without knowing a proof of the Nash--Kuiper theorem about length preserving $C^1$-embeddings.
The original papers of John Nash 
and Nicolaas Kuiper \cite{nash,kuiper} are very readable.

\subsection*{Length-preserving map\thm}
\label{two2one} 

A continuous map $f\:X\to Y$ between metric spaces is called \index{length-preserving}\emph{length-preserving} if it preserves the length of curves; 
that is, for any curve $\alpha$ in $X$ we have
\[\length(f\circ\alpha)=\length\alpha.\]

\begin{pr}
Show that there is no length-preserving map $\RR^2\to \RR$.
\end{pr}

The expected solution uses Rademacher's theorem \cite{rademacher} about differentiability of Lipschitz functions.

\subsection*{Fixed segment}
\label{Fixed segment}

\begin{pr}
Let $\rho(x,y)=\|x-y\|$ be a metric on $\RR^m$ induced by a norm $\|{*}\|$.
Assume that $f\:(\RR^m,\rho)\to(\RR^m,\rho)$ is an isometry that fixes two distinct points $a$ and $b$.
Show that $f$ fixes the line segment between $a$ and~$b$.
\end{pr}

Evidently, $f$ maps the line segment $[ab]$ to a minimizing geodesic connecting $a$ to $b$ in $(\RR^m,\rho)$.
However, in general, there might be many minimizing geodesics connecting $a$ to $b$ in $(\RR^m,\rho)$.
The problem states that $[ab]$ is mapped to itself.

\subsection*{Pogorelov's construction\easy}
\label{Pogorelov's construction}

\begin{pr}
Let $\mu$ be a centrally symmetric Radon measure on $\mathbb{S}^2$ which is positive on every open set and vanishes on every great circle.
Given two points $x,y\in \mathbb{S}^2$,
set 
\[\rho(x,y)=\mu[\,B(x,\tfrac \pi2)\setminus B(y,\tfrac\pi2)\,].\]

Show that $\rho$ is a length-metric on $\mathbb{S}^2$,
and moreover, the geodesics in $(\mathbb{S}^2,\rho)$ run along great circles of $\mathbb{S}^2$.
\end{pr}

\subsection*{Straight geodesics}
\label{Straight geodesics}

Recall that a map $f\:X\to Y$ between metric spaces is called bi-Lipschitz if there is a constant $\eps>0$
such that 
\[\eps\cdot|x-y|_X\le|f(x)-f(y)|_Y\le\tfrac1\eps\cdot|x-y|_X.\]
for any $x,y\in X$.

\begin{pr}
Let $\rho$ be a length-metric on $\RR^m$ that is bi-Lipschitz equivalent to the canonical metric.
Assume that every geodesic $\gamma$ in $(\RR^d,\rho)$ is \emph{affine};
that is, $\gamma(t)=v+w\cdot t$ for constant vectors $v,w\in\RR^m$.

Show that $\rho$ is induced by a norm on $\RR^m$.
\end{pr}

\subsection*{Hyperbolic space}
\label{Hyperbolic space}

\begin{pr}
Construct a bi-Lipschitz map
from the hyperbolic $3$-space 
to the product of two hyperbolic planes.
\end{pr}

\subsection*{Quasi-isometry of a Euclidean space\thm}
\label{hom-near-QI} 

A map $f\:X\to Y$ between metric spaces is called a \index{quasi-isometry}\emph{quasi-isometry} if there is a  real constant $C>1$ such that 
$$\tfrac{1}{C}\cdot|x-x'|_X-C
\le 
|f(x)-f(x')|_Y\le C\cdot|x-x'|_X+C$$
for any $x,x'\in X$ and $f(X)$ is a \index{net}\emph{$C$-net} in $Y$;
that is, for any $y\in Y$ there is $x\in X$ such that $|f(x)-y|_Y\le C$.

Note that a quasi-isometry is not assumed to be continuous;
for example, any map between compact metric spaces is a quasi-isometry.

\begin{pr}
Let $f\:\RR^m\to\RR^m$ be a quasi-isometry.
Show that there is a (bi-Lipschitz) homeomorphism 
$h\:\RR^m\to\RR^m$ at a bounded distance from $f$;
that is, there is a real constant $C$ such that
$$|f(x)-h(x)|\le C$$
for any $x\in\RR^m$.
\end{pr}

The expected solution requires the so-called \emph{gluing theorem},
a corollary of the theorem proved by Laurence Siebenmann \cite{siebenmann}.
It states that 
if $V_1, V_2\subset\RR^m$ are open
and the two embedding $f_1\:V_1\to\RR^m$ and $f_2\:V_2\z\to\RR^m$ 
are sufficiently close to each other 
on the overlap $U=V_1\cap V_2$, 
then
there is an embedding $f$ defined on an open set $W'$
which is slightly smaller than $W=V_1\cup V_2$
and such that $f$ is sufficiently close to each $f_1$ and $f_2$ at the points where they are defined.

The  bi-Lipschitz version requires 
an analogous statement in the category of bi-Lipschitz embeddings;
it was proved by
Dennis Sullivan \cite{sullivan};
a detailed proof is given by Pekka Tukia and Jussi Väisälä \cite[][5.10]{tukia-vaisala}.

\subsection*{Family of sets with no section\easy}
\label{hausdorff-section} 

\begin{pr}
Construct a family of closed sets $C_t\subset\mathbb{S}^1$, $t\z\in [0,1]$ that is continuous in the Hausdorff topology but does not admit a {}\emph{section}.
That is, there is no path $c\:[0,1]\to \mathbb{S}^1$ such that $c(t)\in C_t$ for all $t$.
\end{pr}

\subsection*{Spaces with isometric balls}

\begin{pr}
Construct a pair of locally compact length-metric spaces $X$ and $Y$ 
that are not isometric,
but for some points $x_0\in X$,  $y_0\in Y$ and any radius $R$
the ball $B(x_0,R)_X$ is 
isometric to the ball $B(y_0,R)_Y$.
\end{pr}

\subsection*{Average distance\easy}

\begin{pr}
Show that for any compact length-metric space $X$ there is a number $\ell$ such that for any finite collection of points there is a point $z$ that lies of average distance $\ell$ from the collection;
that is, for any $x_1,\dots,x_n\z\in X$ there is $z\in X$ such that
\[\tfrac1n\cdot\sum_i|x_i-z|_{X}=\ell.\]

\end{pr}

\section*{Semisolutions}

\parbf{Non-contracting map.}
Given any pair of points $x_0,y_0\in K$, 
consider two sequences $x_0,x_1,\dots$ and $y_0,y_1,\dots$
such that $x_{n+1}=f(x_n)$ and $y_{n+1}=f(y_n)$ for each $n$.

Since $K$ is compact, 
we can choose an increasing sequence of integers $n_k$
such that both sequences $(x_{n_i})_{i=1}^\infty$ and $(y_{n_i})_{i=1}^\infty$
converge.
In particular, both are Cauchy sequences;
that is,
\[
|x_{n_i}-x_{n_j}|_K, |y_{n_i}-y_{n_j}|_K\to 0
\ \ 
\text{as}
\ \ \min\{i,j\}\to\infty.
\]

Since $f$ is non-contracting, we get
\[
|x_0-x_{|n_i-n_j|}|
\le 
|x_{n_i}-x_{n_j}|.
\]

It follows that  
there is a sequence $m_i\to\infty$ such that
\[
x_{m_i}\to x_0\ \ \text{and}\ \ y_{m_i}\to y_0\ \ \text{as}\ \ i\to\infty.
\leqno({*})\]

Set \[\ell_n=|x_n-y_n|_K.\]
Since $f$ is non-contracting, the sequence $(\ell_n)$ is non-decreasing.

By $({*})$,  $\ell_{m_i}\to\ell_0$ as $m_i\to\infty$.
It follows that $(\ell_n)$ is a constant sequence.

In particular 
\[|x_0-y_0|_K=\ell_0=\ell_1=|f(x_0)-f(y_0)|_K\]
for any pair of points $(x_0,y_0)$ in $K$.
That is, $f$ is distance-preserving, in particular injective.

From $({*})$, we also get that $f(K)$ is everywhere dense.
Since $K$ is compact $f\:K\to K$ is surjective. Hence the result follows.\qeds

This is a basic lemma in the introduction to Gromov--Hausdorff distance \cite[see 7.3.30 in][]{bbi}.
I learned this proof from Travis Morrison, 
a student in my MASS class at Penn State, Fall 2011.

As an easy corollary, one can see that any surjective non-expanding map from a compact metric space to itself is an isometry.
The following problem due to Aleksander Ca{\l}ka \cite{calka:loc-isom}
is closely related but more involved. 

\begin{pr}
Show that any local isometry from a connected compact metric space to itself is a homeomorphism. 
\end{pr}

\parbf{Finite-whole extension.}
Consider the space $Y^X$ of all maps $X\to Y$ equipped with the product topology.

Given a finite set $F\in X$, denote by $\mathfrak{C}_F$ the set of maps $h\in Y^X$ such that the restriction $h|_F$ is short and the restriction $h|_{A\cap F}$ agrees with $f\:A\to Y$.
By assumption, the sets $\mathfrak{C}_F\subset Y^X$ are closed and nonempty.

Note that for any finite collection of finite sets $F_1,\dots,F_n\subset X$, we have
\[\mathfrak{C}_{F_1}\cap\dots\cap\mathfrak{C}_{F_n}\supset \mathfrak{C}_{F_1\cup\dots\cup F_n}.\]
In particular, the intersection is nonempty.

{\sloppy
According to Tikhonov's theorem \cite[see][and the references therein]{wright}, $Y^X$ is compact.
By the finite intersection property, the intersection $\bigcap_F\mathfrak{C}_F$ with $F$ ranging along all finite subsets of $X$ is nonempty.
It remains to note that any map $h\in \bigcap_F\mathfrak{C}_F$ solves the problem.
\qeds

}

This observation was used by Stephan Stadler and me \cite{petrunin-stadler:revisited}.

\parbf{Horo-compactification.}
For the first part of the problem, take $X$ to be the set of non-negative integers with the metric $\rho$ defined by
\[\rho(m,n)=m+n\] 
for $m\ne n$.

\medskip

The second part is proved by contradiction.
Assume that $X$ is a proper length space and $F_X$ is not an embedding.
That is, there is a sequence of points $z_1,z_2,\dots$ 
and a point $z_\infty$ such that $f_{z_n}\to f_{z_\infty}$ in $C(X,\RR)$
as $n\to \infty$, 
while $|z_n-z_\infty|_X>\eps$ 
for fixed $\eps>0$ and all~$n$.

Note that any pair of points $x,y\in X$ can be connected by a minimizing geodesic $[xy]$.
Choose $\bar z_n$ on a geodesic $[z_\infty z_n]$ such that $|z_\infty-\bar z_n|=\eps$.
Note that 
\begin{align*}
f_{z_n}(z_\infty)-f_{z_n}(\bar z_n)&=\eps
\intertext{and}
f_{z_\infty}(z_\infty)-f_{z_n}(\bar z_n)&=-\eps
\end{align*}
for all $n$.

Since $X$ is proper, we can pass to a subsequence of $z_n$ so that the sequence  $\bar z_n$ converges;
denote its limit by $\bar z_\infty$.
The above identities imply that
\[f_{z_n}(\bar z_\infty)\not\to f_{z_\infty}(\bar z_\infty)
\quad
\text{or}
\quad 
f_{z_n}(z_\infty)\not\to f_{z_\infty}( z_\infty)\]
--- a contradiction.\qeds

I learned this problem from Linus Kramer and Alexander Lytchak;
the example was also mentioned in the lectures of Anders Karlsson
and attributed to Uri Bader \cite[see 2.3 in][]{karlsson}.

\parbf{Approximation of the ball by a sphere.}
Make fine burrows in the standard 3-ball without changing its topology,
but at the same time come sufficiently close to any point in the ball.

Consider the doubling of the obtained ball along  its boundary.
The obtained space is homeomorphic to $\mathbb{S}^3$.
Note that the burrows can be made 
so that the obtained space is sufficiently close to the original ball 
in the Gromov--Hausdorff metric.

It remains to smooth the obtained space slightly 
to get a genuine Riemannian metric with the needed property.\qeds

This problem appeared as an exercise in the textbook of Dmitri Burago, Yuri Burago, and Segei Ivanov \cite[Ex. 7.5.17]{bbi}.

If $M$ is a compact manifold of dimension at least $3$
and $X$ is a reasonable compact length space,
then the existence of a map $M\to X$ that induces a surjective homomorphism on their fundamental groups is a necessary and sufficient condition for the existence of Gromov--Hausdorff approximation of $X$ by Riemannian metrics on $M$.
This statement was proved by Steven Ferry and Boris Okun \cite{ferry-okun}.

\begin{wrapfigure}{o}{41 mm}
\vskip-0mm
\centering
\includegraphics{mppics/pic-501}
\end{wrapfigure}

(A doubled cone over Hawaiian earring gives an example of \emph{unreasonable} space. 
It has a nontrivial fundamental group, but admits an approximation by Riemannian metrics on $\mathbb{S}^3$.)

The two-dimensional case is quite different.
There is no sequence of Riemannian metrics on
$\mathbb{S}^2$ converging to the unit disk in the sense of Gromov--Hausdorff.
In fact, 
if $X$ is a limit of $(\mathbb{S}^2,g_n)$,
then any point $x_0\in X$ either admits a neighborhood homeomorphic to $\RR^2$ or is a cut point;
that is, $X\setminus\{x_0\}$ is disconnected \cite[see 3.32 in][]{gromov-MetStr}.

\parbf{Macroscopic dimension.}
The following claim resembles Besicovitch inequality;
it is key to the proof.
\begin{cl}{$({*})$} Let $a$ be a positive real number.
 Assume that a closed curve $\gamma$ in a metric space $X$ can be subdivided into 4 arcs $\alpha$, $\beta$, $\alpha'$, and $\beta'$ in such a way that 
 \begin{itemize}
 \item $|x-x'|>a$ for any $x\in\alpha$ and $x'\in \alpha'$
 and
 \item $|y-y'|>a$ for any $y\in\beta$ and $y'\in \beta'$.
 \end{itemize}
 Then $\gamma$ is not contractible in its $\tfrac a2$-neighborhood.
\end{cl}

To prove $({*})$, consider two functions defined on $X$ as follows:
\begin{align*}
w_1(x)&=\min \{\,a,\dist_{\alpha}(x)\,\},
\\
w_2(x)&=\min \{\,a,\dist_{\beta}(x)\,\},
\end{align*}
and the map $\bm{w}\:X\to [0,a]\times[0,a]$, defined by
\[\bm{w}\:x\mapsto(w_1(x),w_2(x)).\]

Note that 
\begin{align*}
\bm{w}(\alpha)&=0\times [0,a],
&
\bm{w}(\beta)&=[0,a]\times 0,
\\
\bm{w}(\alpha')&=a\times [0,a],
&
\bm{w}(\beta')&=[0,a]\times a.
\end{align*} 
Therefore, the composition $\bm{w}\circ\gamma$ is a degree 1 map 
\[\mathbb{S}^1\to \partial([0,a]\times[0,a]).\] 
It follows that if $h\:\DD\to X$ shrinks $\gamma$, then there is a point $z\in\DD$ such that 
$\bm{w}\circ h(z)=(\tfrac a2,\tfrac a2)$.
Therefore $h(z)$ lies at distance at least $\tfrac a2$ from $\alpha$, $\beta$, $\alpha'$, $\beta'$
and therefore from $\gamma$.
Hence the claim $({*})$ follows.

\medskip

Choose a point $p\in M$.
Let us cover $M$ by the connected components of the inverse images 
$\dist_p^{-1}((n-1,n+1))$ for all integers $n$.
Clearly, any point in $M$ is covered by at most two of these components.
It remains to show that each of these components has diameter less than $100$.

\begin{wrapfigure}{o}{31 mm}
\vskip-2mm
\centering
\includegraphics{mppics/pic-502}
\end{wrapfigure}

Assume the contrary; let $x$ and $y$ be two points in one connected component 
and $|x-y|_M\ge 100$.
Connect $x$ to $y$ with a curve $\tau$ in this component.
Consider the closed curve $\sigma$ formed by $\tau$ and two geodesics $[px]$, $[py]$.

Note that $|p-x|>40$.
Therefore there is a point $m$ on $[px]$ such that $|m-x|=20$.

By the triangle inequality, the subdivision of $\sigma$ into the arcs $[pm]$, $[mx]$, $\tau$, and $[yp]$ satisfy the conditions of the claim $({*})$ for $a=10$.
Hence the statement follows.\qeds

The problem is due to Mikhael Gromov \cite[Appendix 1(E$_{2}$)]{gromov-filling}.

\parbf{No Lipschitz embedding.}
Consider a chain of circles $c_0,\dots,c_n$ in $\RR^3$;
that is, $c_i$ and $c_{i-1}$ are linked for each $i$. 
\begin{figure}[h!]
\vskip0mm
\centering
\includegraphics{mppics/pic-504}
\end{figure}

Assume that $\RR^3$ is equipped with a length-metric $\rho$ such that the total length of the circles is $\ell$
and $U$ is an open bounded set containing all the circles $c_i$.
Note that for any $L$-Lipschitz embedding $f\:(U,\rho)\z\to\RR^3$, the distance from $f(c_0)$ to $f(c_n)$ is less than $L\cdot\ell$.

The $\rho$-distance from $c_0$ to $c_n$ might be much larger than $L\cdot\ell$.
Indeed, fix a line segment $[ab]$ in $\RR^3$.
Modify 
the length-metric on $\RR^3$ in a small neighborhood of $[ab]$
so that there is a chain $(c_i)$ of circles as above,
that goes from $a$ to $b$ 
such that
(1) the total length, say $\ell$, 
of all the circles $c_i$ is arbitrarily small,
but 
(2) the obtained metric $\rho$ 
is arbitrarily close to the canonical one, say
\[\bigl|\rho(x,y)-|x-y|\bigr|<\eps\]
for any two points $x,y\in\RR^3$
and fixed in advanced small $\eps>0$.
The construction of $\rho$ 
is done by shrinking the length of each circle
and expanding the length in the normal directions  
to the circles in a small neighborhood.
The latter makes it impossible to use the circles $c_i$ as a shortcut; that is, the time needed to go from one circle to another is larger than the time one could save by going along the circle.

Set $a_n=(0,\tfrac1n,0)$ and $b_n=(1,\tfrac1n,0)$.
Note that the line segments $[a_nb_n]$ are disjoint and converging
to $[a_\infty b_\infty]$,
where $a_\infty=(0,0,0)$ and $b_\infty=(1,0,0)$.

Apply the above construction in non-overlapping convex neighborhoods of $[a_nb_n]$ 
for sequences 
$\eps_n$ and $\ell_n$ 
converging to zero very fast.

The obtained length-metric $\rho$ is still close to the canonical metric on $\RR^3$,
but it does not admit 
a locally Lipschitz homeomorphism to $\RR^3$.
Indeed, 
assume that such homeomorphism $h$ exists.
Choose a bounded open set $U$ containing $[a_\infty b_\infty]$;
note that the restriction $h|_U$ is $L$-Lipschitz for some $L$.
From the above construction,
we get 
\begin{align*}
|h(a_\infty)-h(b_\infty)|
&\le 
|h(a_n)-h(b_n)|
+
\\
&\ \ \ \ \ +
|h(a_\infty)-h(a_n)|
+
|h(b_n)-h(b_\infty)|
\le
\\
&\le
L\cdot\ell_n+\tfrac2n+100\cdot\eps_n
\end{align*}
for any positive integer $n$.
The right-hand side converges to $0$ as $n\to\infty$.
Therefore 
\[h(a_\infty)=h(b_\infty)\] 
--- a contradiction.\qeds

The problem is due to
Dmitri Burago, 
Sergei Ivanov,
and David Shoenthal \cite{BIS}.

It is expected that any metric on $\RR^2$ admits locally Lipschitz embeddings into the Euclidean plane.
Also, it seems feasible that any metric on $\RR^3$ admits a locally Lipschitz embedding into $\RR^4$.

Note that any metric on the cube in $\RR^3$ admits a proper locally Lipschitz map to the unit cube with the canonical metric of degree 1.
Moreover, one can make this map injective on any finite set of points.
It is instructive to visualize this map for the metric of the solution.

\parbf{Sub-Riemannian sphere.}
If $d$ is a sub-Riemannian metric on $\mathbb{S}^m$,
then there is a non-decreasing sequence of Riemannian metric tensors
$g_0\z< g_1<\dots$ such that their induced metrics $d_1<d_2<\dots$ converge to $d$.
The metric $g_0$ can be assumed to be the metric of a round sphere,
so it is induced by an embedding $h_0\:\mathbb{S}^m\to \RR^{m+1}$.

Applying the construction from the Nash--Kuiper theorem,
one can produce a sequence of smooth embeddings $h_n\:\mathbb{S}^m\to \RR^{m+1}$ with the induced metrics $g_n'$
such that $|g_n'-g_n|\to 0$.
In particular, if we denote by $d'_n$ the metric corresponding to $g'_n$, then $d'_n\to d$ an $n\to\infty$.

It follows from the same construction that
if one chooses $\eps_n>0$, depending on $h_n$,
then we can assume that 
\[|h_{n+1}(x)-h_n(x)|<\eps_n\] for any $x\in \mathbb{S}^m$.

Let us introduce two conditions on the values $\eps_n$, called \emph{weak} and \emph{strong}.

The weak condition states that $\eps_{n}< \tfrac1{2}\cdot \eps_{n-1}$ for any $n$.
This ensures that the sequence of maps $h_n$ converges pointwise;
denote its limit by $h_\infty$.

Denote by $\bar d$ the length-metric induced by $h_\infty$.
Note that $\bar d\le d$.
The strong condition on $\eps_n$ will ensure that actually $\bar d=d$.

Fix $n$ and assume that $h_n$ and therefore $\eps_{n-1}$ are constructed already.
Set $\Sigma=h_n(\mathbb{S}^m)$, and let $\Sigma_r$ be the tubular $r$-neighborhood of $\Sigma$.
Equip $\Sigma$ and $\Sigma_r$ with the induced length-metrics.
Since $\Sigma$ is a smooth hypersurface, we can choose $r_n\in(0,\eps_{n-1}]$ 
so that the inclusion $\Sigma\hookrightarrow \Sigma_{r_n}$ preserves the distance up to error $\tfrac1{2^n}$.
Then the strong condition states that $\eps_n< \tfrac12\cdot r_n$, 
which is evidently stronger than the weak condition  $\eps_{n}< \tfrac1{2}\cdot \eps_{n-1}$ above.

Note that if the sequence $h_n$ is constructed with the described choice of $\eps_n$,
then $|h_\infty(x)-h_n(x)|<r_n$ for any $x\in\mathbb{S}^m$.
Therefore 
\[\bar d(x,y)+2\cdot r_n+\tfrac1{2^n}\ge d_n'(x,y)\] 
for any $n$ and $x,y\in \mathbb{S}^m$;
hence $\bar d\ge d$ as required. 
\qeds

The problem
on this list was first discovered by Enrico Le Donne \cite{le-donne}.
A similar construction is described in the lecture notes by Allan Yashinski and the author \cite{petrunin-yashinsky} 
which are aimed for undergraduate students. 
Yet the results in \cite{petrunin-paths} are closely relevant.

The construction in the Nash--Kuiper embedding theorem
can be used to prove strange statements.
Here is one example based on the observation that the Weyl curvature tensor 
vanishes on hypersurfaces in the Euclidean space.

\begin{pr}
Let $M$ be a Riemannian manifold diffeomorphic to the $m$-sphere. 
Show that there is a Riemannian manifold $M'$ arbitrarily close to $M$ in the Lipschitz metric whose Weyl curvature tensor is identically 0.
\end{pr}

\parbf{Length-preserving map.}
Assume the contrary;
let $f\:\RR^2\to \RR$ be a length-preserving map.

Note that $f$ is Lipschitz.
Therefore by Rademacher's theorem \cite{rademacher}, the differential $d_xf$ is defined for  almost all $x$.

Choose a unit vector $u$.
Given $x\in\RR^2$,
consider the path $\alpha_x(t)\z=x+t\cdot u$ defined for $t\in [0,1]$.
Note that  
\[\alpha'_x(t)=(d_{\alpha_x(t)}f)(u)\]
holds for almost all $x$ and $t$.
It follows that 
\[\length(f\circ\alpha_x)=\int\limits_0^1 |(d_{\alpha_x(t)}f)(u)| \cdot dt\]
for almost all $x$.

Therefore $|d_xf(v)|=|v|$ for almost all $x,v\in\RR^2$.
In particular, there is $x\in\RR^2$ such that the differential $d_xf$ is defined 
and 
\[|d_xf(e_1)|=|e_1|,
\quad
|d_xf(e_2)|=|e_2|,
\quad
|d_xf(e_1+e_2)|=|e_1+e_2|\]
for a basis $(e_1,e_2)$ of $\RR^2$.
It follows that $d_xf$ has rank 2 --- a contradiction. \qeds

The idea above can also be used to solve the following problem.

\begin{pr} Let $\rho$ be a metric on $\RR^2$ that is induced by a norm.
Show that $(\RR^2,\rho)$ admits 
a length-preserving map
to $\RR^3$ 
if and only if 
$(\RR^2,\rho)$ is isometric to the Euclidean plane.
\end{pr}

\parbf{Fixed segment.}
Note that it is sufficient to show that 
if $f$ is an isometry such that
\[f(a)=a\ \ \text{and}\ \ f(b)=b\]
for some $a,b\in\RR^m$,
then 
\[f(\tfrac{a+b}2)=\tfrac12\cdot(f(a)+f(b)).\]

Without loss of generality, we can assume that $b+a=0$.

Set $f_0=f$.
Consider the sequence of isometries $f_0$, $f_1,\dots$ recursively defined by
\[f_{n+1}(x)= -f_n^{-1}(-f_n(x))\]
for all $n$.

Note that for all $n$ we have $f_n(a)=a$, $f_n(b)=b$ and 
$$|f_{n+1}(0)|=2\cdot|f_n(0)|.$$
Therefore  
if $f(0)\ne 0$,
then $|f_n(0)|\to\infty$ as $n\to\infty$.

On the other hand, since $f_n$ is isometry and $f(a)=a$,
we also have $|f_n(0)|\le 2\cdot |a|$ --- a contradiction.
\qeds

The idea of the proof is due to Jussi Väisälä \cite{vaisala}.
The problem is the main step in the proof of the Mazur--Ulam theorem \cite{mazur-ulam},
which states that any isometry of $(\RR^m,\rho)$ is an affine map.


\parbf{Pogorelov's construction.}
The positivity and symmetry of $\rho$ are evident.

The triangle inequality follows since
\[[B(x,\tfrac \pi2)\setminus B(y,\tfrac\pi2)]
\cup 
[B(y,\tfrac\pi2)\setminus B(z,\tfrac\pi2)]
\supseteq
B(x,\tfrac \pi2) \setminus B(z,\tfrac\pi2).
\leqno(*)\]

\begin{wrapfigure}[8]{o}{31 mm}
\vskip-2mm
\centering
\includegraphics{mppics/pic-505}
\end{wrapfigure}

Observe that
$B(x,\tfrac \pi2)\setminus B(y,\tfrac\pi2)$
does not overlap
$B(y,\tfrac\pi2)\setminus B(z,\tfrac\pi2)$ and  we get equality in $(*)$ if and only if $y$ lies on the great circle arc from $x$ to $z$.
Therefore the second statement follows.\qeds

This construction was given by 
Aleksei Pogorelov \cite{pogorelov}.
It is closely related to the construction given 
by David Hilbert in \cite{hilbert}
which was the motivating example for his 4th problem.

The following problem was suggested by Jairo Bochi \cite{bochi};
it looks similar, but actually quite different.

\begin{pr}
Consider the set $E$ of all ellipses with unit area centered at the origin of the plane.
Equip $E$ with the metric defined by 
\[\rho(A,B)=\area(A\triangle B),\]
where $\triangle$ denotes symmetric difference; that is, $A\triangle B=(A\setminus B)\cup (B\setminus A)$.
Let $\hat\rho$ be the length-metric induced by $\rho$.
Show that $(E,\hat\rho)$ is isometric to a Lobachevsky plane.
\end{pr}

The metrics of the form $\rho(A,B)=\mu(A\bigtriangleup B)$ are very special;
evidently, they can be embedded into the corresponding $L^1$ space.
Also, they satisfy the so-called \index{hypermetric inequality}\emph{hypermetric inequalities}; that is, 
\[\sum_{i,j}b_i\cdot b_j\cdot \rho(A_i,A_j)\le 0\]
for any sequence of sets $A_1,\dots A_n$ and any sequence of integers $b_1,\dots,b_n$ such that $\sum_ib_i\z=1$.
Note that for $n=3$ and $b_1=b_2=-b_3=1$ we get the usual triangle inequality.
For more on the subject, see \cite{deza-laurent}.

\parbf{Straight geodesics.}
From the uniqueness of the straight segment between two given points in $\RR^m$,
it follows that any straight line in $\RR^m$ is a geodesic in $(\RR^m,\rho)$.

Set 
\[\|\bm{v}\|_{\bm{x}}=\rho(\bm{x},(\bm{x}+\bm{v})).\]
Note that 
\[ \|\lambda\cdot\bm{v}\|_{\bm{x}}
=
|\lambda|\cdot\|\bm{v}\|_{\bm{x}}\]
for any $\bm{x}$, $\bm{v}\in\RR^m$, and $\lambda\in\RR$.

Denote by $|x-y|$ the Euclidean distance between the points $x$ and~$y$.
Since $\rho$ and $|{*}-{*}|$ are bi-Lipschitz equivalent,
applying the triangle inequality twice to the points $\bm{x}$, $\bm{x}+\lambda\cdot\bm{v}$, $\bm{x}'$ and $\bm{x}'+\lambda\cdot\bm{v}$, we get
\[
\bigl|\|\lambda\cdot\bm{v}\|_{\bm{x}}
-
\|\lambda\cdot\bm{v}\|_{\bm{x}'}\bigr|
\le 
C\cdot |\bm{x}-\bm{x'}|\]
for any $\bm{x},\bm{x'},\bm{v}\in\RR^m$, 
$\lambda\in\RR$
and a fixed real constant $C$.

Passing to the limit as $\lambda\to\infty$, 
we obtain that
$\|\bm{v}\|_{\bm{x}}$ does not depend on $\bm{x}$;
hence the result follows.\qeds

This idea is due to Thomas Foertsch
and Viktor Schroeder \cite{foertsch-schroeder}.
A more general statement was proved by Petra Hitzelberger and Alexander Lytchak \cite{hitzelberger-lytchak}.
Namely, they showed that 
if any pair of points in a geodesic metric space $X$ can be separated by an \emph{affine function},
then $X$ is isometric to a convex subset of a normed vector space.
(A function $f\:X\to\RR$ is called affine if, for any geodesic $\gamma$ in $X$, the composition $f\circ\gamma$ is affine.)

\parbf{Hyperbolic space.}
The hyperbolic plane $\HH^2$ is isometric to $(\RR^2,g)$, where 
\[g(x,y)=\left(\begin{matrix}
     1&0
     \\
     0&e^{x}
    \end{matrix}\right).\]
The same way, the hyperbolic space $\HH^3$
can be viewed as $(\RR^3,h)$, where 
\[h(x,y,z)=\left(\begin{matrix}
     1&0&0
     \\
     0&e^{x}&0
     \\
     0&0&e^{x}
\end{matrix}\right).\]
    
In the described coordinates, consider the projections $\phi,\psi\:\HH^3\to\HH^2$ defined by 
$\phi\:(x,y,z)\mapsto (x,y)$ and $\psi\:(x,y,z)\mapsto (x,z)$.
Note that 
\begin{align*}
\max&\{\,|\phi(p)-\phi(q)|_{\HH^2},|\psi(p)-\psi(q)|_{\HH^2}\,\}
\le
\\
&\le
|p-q|_{\HH^3}
\le
\\
&\le
|\phi(p)-\phi(q)|_{\HH^2}+ |\psi(p)-\psi(q)|_{\HH^2}
\end{align*}
for any two points $p,q\in \HH^3$.
In particular, the map $\HH^3\to\HH^2\times\HH^2$ defined by $p\mapsto (\phi(p),\psi(p))$
is $2^{\mp1}$-bi-Lipschitz.\qeds

We used that horo-spheres in the hyperbolic space are isometric to the Euclidean plane.
This observation was made by Nikolai Lobachevsky \cite[see 34 in][]{lobachevsky}.
The same observation is used in the following construction discovered by 
K\'{a}roly Böröczky [see \ncite{boroczky} and also \ncite{radin}]. 

\begin{pr}
Construct a tessellation of the hyperbolic plane with one polygonal tile of arbitrarily small area and/or diameter.  
\end{pr}

\parbf{Quasi-isometry of a Euclidean space.}
Choose two constants $M\ge 1$ and $A\ge 0$.
A map $f\:X\z\to Y$ between metric spaces $X$ and $Y$ such that for any $x,y\in X$,
 we have
\[\tfrac1M\cdot |x-y|-A\le |f(x)-f(y)|\le M\cdot |x-y|+A\]
and any point in $Y$ lies on the distance at most $A$ from a point in the image $f(X)$
will be called $(M,A)$-quasi-isometry.

{\sloppy
Note that $(M,0)$-quasi-isometry is a $[\tfrac1M,M]$-bi-Lipschitz map.
Moreover,
if $f_n\:\RR^m\to\RR^m$ is a  $(M,\frac1n)$-quasi-isometry 
for each $n$, 
then any subsequential limit of $f_n$ as $n\to\infty$
is a $[\tfrac1M,M]$-bi-Lipschitz map.

}

Therefore given $M\ge 1$ and $\eps>0$ there is $\delta>0$ such that 
for any $(M,\delta)$-quasi-isometry $f\:\RR^m\to\RR^m$ and any $p\in \RR^m$
there is an $[\tfrac1M,M]$-bi-Lipschitz map $h\:B(p,1)\to \RR^m$
such that
\[|f(x)-h(x)|<\eps\]
for any $x\in B(p,1)$.

Using rescaling, we can get the following equivalent formulation. 
Given $M\ge 1$, $A\ge 0$, and $\eps>0$,
there is sufficiently large $R>0$ such that for any $(M,A)$-quasi-isometry 
$f\:\RR^m\to\RR^m$ and any $p\in\RR^m$ there is a $[\tfrac1M,M]$-bi-Lipschitz map $h\:B(p,R)\to \RR^m$
such that 
\[|f(x)-h(x)|<\eps\cdot R\]
for any $x\in B(p,R)$.

Cover $\RR^m$ by balls $B(p_n,R)$ and construct a $[\tfrac1M,M]$-bi-Lipschitz map $h_n\:B(p_n,R)\to \RR^m$ close to the restrictions $f|_{B(p_n,R)}$ for each $n$.

The maps $h_n$ are $2\cdot \eps\cdot R$ close to each other on the overlaps of their domains of definition.
This makes it possible to deform slightly each $h_n$ so that they agree on the overlaps.
This can be done by Siebenmann's theorem \cite{siebenmann}.
If instead you apply Sullivan's theorem [see \ncite{sullivan} and also 5.10 in \ncite{tukia-vaisala}], you get a bi-Lipschitz homeomorphism $h\:\RR^m\z\to\RR^m$.\qeds

The problem was suggested by Dmitri Burago.

\parbf{Family of sets with no section.} 
Given $t\in (0,1]$, consider the real interval $\tilde C_t=[\tfrac 1t+t, \tfrac 1t+1]$.
Denote by $C_t$ the image of $\tilde C_t$ under the covering map $\pi\:\RR\to \mathbb{S}^1=\RR/\ZZ$.

Set $C_0=\mathbb{S}^1$.
Note that the Hausdorff distance from $C_0$ to $C_t$ is $\tfrac t2$.
Therefore $\{C_t\}_{t\in[0,1]}$ is a family of compact subsets in $\mathbb{S}^1$ that is continuous in the sense of Hausdorff.

\medskip

Assume there is a continuous section $c(t)\in C_t$, for $t\in [0,1]$.
Since $\pi$ is a covering map,
we can lift the path $c$ to a path $\tilde c\:[0,1]\to \RR$ such that $\tilde c(t)\in \tilde C_t$ for all $t$.
In particular, $\tilde c(t)\to\infty$ as $t\to0$
--- a contradiction.\qeds

The problem was suggested by Stephan Stadler.
Here is a simpler, closely related problem.

{

\begin{wrapfigure}{o}{31 mm}
\vskip-0mm
\centering
\includegraphics{mppics/pic-506}
\end{wrapfigure}

\begin{pr}
Show that any Hausdorff continuous family of compact sets in $\RR$ admits a continuous section.
\end{pr}

The existence of sections for a family of sets parametrized by a topological space was considered by Ernest Michael \cite{michael-1,michael-2,michael-3}.


\parbf{Spaces with isometric balls.} 
The needed examples can be constructed by cutting the upper half-plane along a ``dyadic comb'' shown on the diagram;
the obtained space should be equipped with the intrinsic metric induced from the $\ell_\infty$-norm on the plane. 

}

\medskip

First, let us describe the comb precisely.
Choose an infinite sequence $a_0,a_1,\dots$ of zeros and ones.
Given an integer $k$, cut the upper half-plane along the line segment between $(k,0)$ and $(k,2^{m+1})$ 
if $m$ is the maximal number such that 
\[k\equiv a_0+2\cdot a_1+\dots+2^{m-1}\cdot a_{m-1}\pmod{2^{m}};\]
If the equality holds for all $m$, cut the half-plane along the vertical half-line starting at $(k,0)$.

Note that all the obtained spaces, independently from the sequence $(a_n)$, meet the conditions of the problem for the point $x_0=(\tfrac12,0)$.

Further, note that the resulting spaces for two sequences $(a_n)$ and $(a'_n)$ are isometric only in the following two cases 
\begin{itemize}
\item if $a_n=a_n'$ for all large $n$, or
\item if $a_n=1-a_n'$  for all large $n$.
\end{itemize}

It remains to produce two sequences that do not have these identities for all large $n$; 
two random sequences of zeros and ones will do the job with probability one.\qeds

\parbf{Average distance.}
If such a number does not exist then the ranges of average distance functions have empty intersection.
Since $X$ is a compact length-metric space, the range of any continuous function on $X$ is a closed interval.
By 1-dimensional Helly's theorem, there is a pair of such range intervals that do not intersect.
That is, for two point-arrays $(x_1,\dots,x_n)$ and $(y_1,\dots,y_m)$
and their average distance functions 
\[f(z)=\tfrac1n\cdot\sum_i|x_i-z|_X\quad\text{and}\quad h(z)=\tfrac1m\cdot\sum_j|y_j-z|_X,\] we have 
$$\min\set{f(z)}{z\in X}>\max\set{h(z)}{z\in X}.\leqno({*})$$

Note that 
$$\tfrac1m\cdot\sum_j f(y_j)=\tfrac1{m\cdot n}\cdot\sum_{i,j}|x_i-y_j|_X=\tfrac1n\cdot\sum_i h(x_i);$$
that is, the average value of $f(y_j)$ coincides with the average value of $h(x_i)$, 
which contradicts $({*})$.
\qeds

This is a result of Oliver Gross \cite{gross}. 
The value $\ell$ is called the \emph{rendezvous value} of $X$;
in fact, it is uniquely defined.

\csname @openrightfalse\endcsname
\chapter{Actions and coverings}

\subsection*{Bounded orbit}
\label{Bounded orbit}

Recall that a metric space is called \index{proper metric space}\emph{proper} if all its bounded closed sets are compact.

\begin{pr} Let $X$ be a 
proper metric space 
and $\iota\:X\to X$ an isometry.
Assume that for some $x\in X$, the sequence $x_n\z=\iota^n(x)$, $n\in\ZZ$ has a convergent subsequence.
Prove that the sequence $x_n$ is bounded.
\end{pr}

\parit{Semisolution.}
Note that we can assume that the orbit $\{x_n\}$ is dense in $X$;
otherwise, we can pass to the closure of the orbit.
In particular, we can choose a finite number of positive integers $n_1,\dots,n_k$
such that the set of points $\{x_{n_1},\dots,x_{n_k}\}$ is a $1$-net for the ball $B(x_0,10)$;
that is, for any $x\in B(x_0,10)$ there is $x_{n_i}$ such that
\[|x-x_{n_i}|<1.\]

Assume that $x_m\in B(x_0,1)$ for some $m$.
Then 
\[B(x_m,10)=f^m( B(x_0,10))\supset B(x_0,1).\] 
In particular, $\{x_{m+n_1},\dots,x_{m+n_k}\}$ is a $1$-net for the ball $B(x_0,1)$
Therefore $x_{m+n_i}\in B(x_0,1)$ for some $i\z\in\{1,\dots,k\}$.

Set $N=\max_i\{n_i\}$.
Applying the above observation inductively, we get that at least one point from any string $x_{i+1},\dots x_{i+N}$ lies in $B(x_0,1)$.
In particular, the $N$ balls
\[B(x_1,10),\dots,B(x_N,10)\]
cover whole $X$.
Hence the result follows.\qeds

The problem is due to Aleksander Ca{\l}ka \cite{calka}.

\subsection*{Finite action}\label{Finite action}

\begin{pr}
Show that for any nontrivial continuous action of a finite group on the unit sphere
there is an orbit that does not lie in the interior of a hemisphere.
\end{pr}


\subsection*{Covers of the figure eight}\label{figure-eight-1}

Given a covering 
\[f\:\tilde X \to X\]
of the length-metric space $X$,
one can consider the induced length-metric on $\tilde X$,
defining the length of curve $\alpha$ in $\tilde X$ as the length of the composition $f\circ\alpha$; the obtained metric space $\tilde X$ is called the  \index{metric covering}\emph{metric covering} of $X$.

{

\begin{wrapfigure}{o}{23 mm}
\vskip-4mm
\centering
\includegraphics{mppics/pic-602}
\end{wrapfigure}

Let us define the \index{figure eight}\emph{figure eight} as the length-metric space obtained by gluing together all four ends of two unit segments.

}

\begin{pr}
Show that any compact length-metric space 
is a Gromov--Hausdorff limit of a sequence $(\widetilde \Phi_n, \tilde d/n)$
where 
\[(\widetilde \Phi_n, \tilde d)\to(\Phi,d),\]
are metric coverings of the figure eight $(\Phi,d)$.
\end{pr}

\subsection*{Diameter of \textit{m}-fold covering\hard}\label{m-fold-cover}

The metric covering is defined in the previous problem.

\begin{pr}
Let $X$ be a length-metric space,
and let $\tilde X$ be an $m$-fold metric covering of $X$.
Show that
$$\diam\tilde X\le m\cdot \diam X.$$
\end{pr}

The figure below shows a 5-fold covering with the diameter of the total space being exactly 5 times the diameter of the target.

\begin{figure}[h!]
\vskip0mm
\centering
\includegraphics{mppics/pic-604}
\end{figure}

\subsection*{Symmetric square\easy}\label{Symmetric square}

Let $X$ be a topological space.
Note that $X{\times} X$ admits a natural $\ZZ_2$-action generated by the involution $(x,y)\mapsto (y,x)$.
The quotient  space $X{\times} X/\ZZ_2$ is called the \index{symmetric square}\emph{symmetric square} of $X$.

\begin{pr} 
Show that the symmetric square 
of any path-connected topological space 
has a commutative fundamental group.
\end{pr}

{

\begin{wrapfigure}{r}{23 mm}
\vskip-0mm
\centering
\includegraphics{mppics/pic-606}
\end{wrapfigure}

\subsection*{Sierpi\'nski gasket\easy}\label{Sierpinski triangle}

To construct the Sierpi\'nski gasket, start with a solid  equilateral triangle, subdivide it into four smaller congruent equilateral triangles and remove the interior of the central one.
Repeat this procedure recursively for each of the remaining solid triangles.

}

\begin{pr} 
Find the homeomorphism group of the Sierpi\'nski gasket.
\end{pr}

\subsection*{Lattices in a Lie group}\label{Boys and girls in a Lie group}

\begin{pr}
Let $L$ and $M$ be two discrete subgroups of a connected Lie group $G$,
and let $h$ be a left-invariant metric on $G$.
Equip the groups $L$ and $M$ 
with the metric induced from $G$.
Assume that $L\setminus G$ and $M\setminus G$ are compact and
$$\vol(L\setminus (G,h))
=
\vol(M\setminus (G,h)).$$
Prove that there is a bi-Lipschitz one-to-one mapping $f\:L
\to
M$, not necessarily a homomorphism.
\end{pr}

\subsection*{Piecewise Euclidean quotient}\label{Piecewise Euclidean quotient}

Note that the quotient of the Euclidean space by a finite subgroup of $\SO(m)$ is a {}\emph{polyhedral space} as is defined on page \pageref{piecewise linear map};
on the same page, you can find the definition of piecewise linear homeomorphism.

\begin{pr}
Let $\Gamma$ be a finite subgroup of $\SO(m)$.
Denote by $P$ the quotient $\RR^m/\Gamma$ equipped with the induced
polyhedral metric.
Assume that $P$ admits a piecewise linear homeomorphism to $\RR^m$.
Show that $\Gamma$ is generated by rotations  around subspaces of codimension $2$.
\end{pr}

The action of the symmetric group $S_m$ on $\CC^m=\RR^{2\cdot m}$ by permutation of complex coordinates provides a remarkable example.
The homeomorphism $\CC^m/S_m\to \CC^m$ can be given by symmetric polynomials on $\CC^m$;
that is, $(z_1,\dots,z_m)\mapsto (a_0,\dots,a_{m-1})$, where
\[(z+z_1)\cdots(z+z_m)=a_0+a_1\cdot z+\dots+a_{m-1}\cdot z^{m-1}+z^m.\]
This homeomorphism is isotopic to a piecewise linear homeomorphism.

\subsection*{Subgroups of a free group}\label{Subgroups of free group} 

\begin{pr}
Show that every finitely generated subgroup of a free group 
is an intersection of subgroups of finite index.
\end{pr}

\subsection*{Short generators\easy}\label{Lengths of generators of the fundamental group}

\begin{pr}
Let $M$ be a compact Riemannian manifold and $p\in M$.
Show that the fundamental group $\pi_1(M,p)$
is generated by the homotopy classes of the loops with length at most $2\cdot\diam M$.
\end{pr}

\subsection*{Number of generators}\label{Number of generators}

\begin{pr}
Let $M$ be a complete connected Riemannian manifold with non-negative sectional curvature.
Show that the minimal number of generators of the fundamental group $\pi_1 M$
can be bounded above in terms of the dimension of $M$.
\end{pr}

\subsection*{An equation in a Lie group\easy}\label{Equations in the group}

\begin{pr}
Let $G$ be a compact connected Lie group and $n$ a positive integer.
Show that given a collection of elements $g_1,\dots,g_n\in G$
the equation 
\[x\cdot g_1\cdot x\cdot g_2\cdots x\cdot g_n=e\]
has a solution $x\in G$;
here $e$ is the identity element in $G$.
\end{pr}

\subsection*{Quotient of the Hilbert space\hard}\label{Quotient of Hilbert space}

\begin{pr}
Construct an isometric action on the Hilbert space with the quotient space isometric to the sphere $\mathbb{S}^3$.
\end{pr}

\section*{Semisolutions}
\parbf{Finite action.}
Without loss of generality, we may assume that the action is generated by a nontrivial homeomorphism 
\[a\:\mathbb{S}^m\to\mathbb{S}^m\] 
of prime order $p$.

Assume the contrary; that is, assume that any $a$-orbit lies in an open hemisphere.
Then 
\[h(x)=\sum_{n=1}^p a^n\cdot x\ne0\]
for any $x\in\mathbb{S}^m$; here we consider $\mathbb{S}^m$ as the unit sphere in $\mathbb{R}^{m+1}$.

Consider the map $f\:\mathbb{S}^m\to\mathbb{S}^m$ 
defined by $f(x)=\tfrac{h(x)}{|h(x)|}$.
Note that 
\begin{enumerate}[(a)]
\item if $a(x)=x$, then $f(x)=x$;
\item\label{f(x)=f(a(x))} $f(x)=f\circ a(x)$ for any $x\in\mathbb{S}^m$.
\end{enumerate}

Note further that $f$ is homotopic to the identity; 
in particular 
\[\deg f=1.
\leqno({*})\]
The homotopy can be defined by $(x,t)\mapsto \gamma_x(t)$,
where $\gamma_x$ is the minimizing geodesic path in $\mathbb{S}^m$ from $x$ to $f(x)$.
By construction, $|x-f(x)|_{\mathbb{S}^m}<\tfrac\pi2$; 
therefore $\gamma_x$ is uniquely defined.

Choose $x\in \mathbb{S}^m$ such that $a(x)\ne x$.
Note that the group acts without fixed points 
on the inverse image $W=f^{-1}(V)$ 
of a small open neighborhood $V\ni x$.
Therefore the quotient map $\theta\:W\z\to W'\z=W/\ZZ_p$ is a $p$-fold covering.
By (\ref{f(x)=f(a(x))}),
the restriction $f|_W$ factors thru $\theta$;
that is,
there is $f'\:W'\to V$ such that
$f|_W=f'\circ\theta$.

Assume that $p\ne 2$.
Note that $f'$ and $\theta$ have well-defined degrees and 
\[\deg f\equiv\deg \theta\cdot\deg f'\pmod p.\]
Since $\theta$ is a $p$-fold covering, we have $\deg \theta\equiv0\pmod p$.
Therefore
\[\deg f\equiv 0\pmod p.
\leqno({*}{*})\]

Finally, observe that $({*})$ contradicts $({*}{*})$.

In the case $p=2$ the same proof works, 
but the degrees have to be considered modulo $2$.\qeds

Along the same lines, one can get a lower bound for the maximal diameter of the orbits for any nontrivial action of a finite group on a Riemannian manifold.

Applying the problem to the conjugate actions, 
one gets that if a fixed point set of a finite group acting on a sphere
has nonempty interior, 
then the action is trivial.
The same holds for any connected manifold.
All this was proved by Max Newman \cite{newman}.

The following problem from \cite{montgomery} can be solved using Newman's theorem. 

\begin{pr}
Let $h$ be a homeomorphism of a connected manifold $M$ 
such that each $h$-orbit is finite.
Show that $h$ has finite order.
\end{pr}

\parbf{Covers of the figure eight.}
First note that any compact length-metric space $K$ can be approximated by finite metric graphs.

Indeed, fix a finite $\eps$-net $F$ in $K$.
For each pair $x,y\in F$ choose a chain of points $x=x_0,x_1,\dots, x_n=y$ such that
$|x_i-x_{i-1}|_K<\eps$ for each $i$ and 
\[|x-y|_K=|x_0-x_1|_K+\dots+|x_{n-1}-x_n|_K.\]
Denote by $F'$ the union of all these chains with $F$.
Connect a pair of vertices $v,w\in F'$ by an edge of length $|v-w|_K$ if $|v-w|_K<\eps$.
Note that the obtained metric graph is $\eps$ close to $K$ in the Gromov--Hausdorff metric.

\begin{wrapfigure}{o}{25 mm}
\vskip-3mm
\centering
\includegraphics{mppics/pic-610}
\vskip-3mm
\end{wrapfigure}

Further, any finite metric graph can be approximated by a graph made from the fragments shown on the diagram
(we have to attach each pair of free ends of one fragment to a pair of ends in another fragment).

It remains to observe that metric graphs obtained from these fragments are finite coverings of $(\Phi,d/n)$.
\qeds

The same idea works if instead of the figure eight, we have a compact length-metric space $X$ that admits a map $X\to\Phi$ inducing an epimorphism of fundamental groups.
Such spaces $X$ can be found among compact hyperbolic manifolds of any dimension $\ge 2$.
All this is due to Vedrin Šahović \cite{sahovic}.

A similar idea was used later to show that any finitely presented group can appear as a fundamental group of the underlying space of a 3-dimensional hyperbolic orbifold \cite{panov-petrunin-telescopic}.

\parbf{Diameter of \textit{m}-fold covering.}
Choose points $\tilde p,\tilde q\in\tilde M$.
Let  
$\tilde\gamma\:[0,1]\z\to \tilde M$ be a minimizing geodesic path from $\tilde p$ to $\tilde q$. 

We need to show that 
\[\length \tilde\gamma\le m\cdot \diam M.\]
Suppose the contrary.

Denote by $p$, $q$, and $\gamma$ the projections to $M$ of $\tilde p$, $\tilde q$, and $\tilde \gamma$ respectively. 
Represent $\gamma$
as the concatenation of $m$ paths of equal length,
\[\gamma=\gamma_1{*}\dots{*}\gamma_m,\] 
so
\[\length\gamma_i=\tfrac{1}m\cdot\length\gamma>\diam M.\] 

Let $\sigma_i$ be a minimizing geodesic in $M$ connecting the endpoints of $\gamma_i$. 
Note that 
\[\length\sigma_i\le \diam M< \length\gamma_i.\] 

Consider $m+1$ paths $\alpha_0,\dots,\alpha_m$ defined as the concatenations 
\[\alpha_i=\sigma_1{*}\dots{*}\sigma_i{*}\gamma_{i+1}{*}\dots{*}\gamma_m.\]

Let $\tilde\alpha_0,\dots,\tilde\alpha_m$ be their liftings
with $\tilde q$ as an endpoint.
The starting point of each curve $\tilde\alpha_i$ is one of $m$ inverse images of $p$. 
Therefore two curves, $\tilde\alpha_i$ and $\tilde\alpha_j$ for $i<j$, 
have the same starting point in $\tilde M$.

Note that the concatenation
\[\beta=\gamma_1{*}\dots{*}\gamma_i{*}\sigma_{i+1}{*}\dots{*}\sigma_j{*}\gamma_{j+1}{*}\dots{*}\gamma_m.\]
admits a lift $\tilde\beta$ that connects $\tilde p$ to $\tilde q$ in $\tilde M$.
Clearly, $\length \tilde\beta<\length \gamma$ --- a contradiction.
\qeds

The question was asked by Alexander  Nabutovsky
and answered by Sergei Ivanov \cite{ivanov}.
A closely related problem for universal coverings is discussed by Sergio Zamora in \cite{zamora}.

\parbf{Symmetric square.}
Let $\Gamma=\pi_1 X$ and $\Delta=\pi_1((X\times X)/\ZZ_2)$.
Consider the homomorphism $\phi\:\Gamma\times \Gamma\to \Delta$
induced by the quotient map $X\times X\z\to (X\times X)/\ZZ_2$.

Note that $\phi(\alpha,1)=\phi(1,\alpha)$ for any $\alpha\in \Gamma$ and the restrictions $\phi|_{\Gamma\times \{1\}}$ and $\phi|_{\{1\}\times\Gamma}$
are onto.

It remains to note that 
$$\phi(\alpha,1)\cdot\phi(1,\beta)=\phi(1,\beta)\cdot\phi(\alpha,1)$$
for any $\alpha$ and $\beta$ in $\Gamma$.
\qeds

The problem was suggested by Rostislav Matveyev.


\parbf{Sierpi\'nski gasket.}
Denote the Sierpi\'nski gasket by~$\triangle$.

Let us show that any homeomorphism of $\triangle$ is also an isometry.
Therefore its homeomorphism group is the symmetric group $S_3$. 

{

\begin{wrapfigure}{o}{23 mm}
\vskip-4mm
\centering
\includegraphics{mppics/pic-607}
\end{wrapfigure}
Let $\{x,y,z\}$ be a 3-point set in $\triangle$ such that its complement has 3 connected components.
Show that there is a unique choice for the set $\{x,y,z\}$ and 
it is formed by the midpoints of the long sides.

It follows that any homeomorphism of $\triangle$ permutes the set $\{x,y,z\}$.

}

Applying a similar argument recursively to the smaller triangles,
we get that this permutation uniquely describes the homeomorphism.
\qeds

The problem was suggested by Bruce Kleiner.
The homeomorphism group of the Sierpi\'nski carpet is much more interesting \cite{kapovich-kleiner}.

\parbf{Lattices in a Lie group.}
Denote by $V_\ell$ and $W_m$
the Voronoi domains for each $\ell\in L$ and $m\in M$ respectively;
that is,
\begin{align*}
V_\ell&=\set{g\in G}{|g-\ell|_G\le|g-\ell'|_G\ \text{for any}\ \ell'\in L},
\\
W_m&=\set{g\in G}{|g-m|_G\le|g-m'|_G\ \text{for any}\ m'\in M}.
\end{align*}

Note that for any $\ell\in L$ and $m \in M$, we have
\[
\vol V_\ell=\vol(L\setminus (G,h))=\vol(M\setminus (G,h))=\vol W_m.
\leqno({*})
\]

Consider the bipartite graph $\Gamma$ with the parts $L$ and $M$
such that $\ell\in L$ is adjacent  to $m \in M$ if and only if $V_\ell\cap W_m\ne\emptyset$.

By $({*})$ the graph $\Gamma$ satisfies the condition of the marriage theorem \cite{hall-marriage}  ---
any subset $S$ in $L$ has at least $|S|$ neighbors in $M$ and the other way around;
here $|S|$ denotes the number of elements in $S$.
Therefore there is a bijection $f\: L\to M$ such that 
\[V_\ell\cap W_{f(\ell)}\ne\emptyset\] for any $\ell\in L$. 

It remains to observe that $f$ is bi-Lipschitz.
\qeds

The problem is due to 
Dmitri Burago 
and Bruce Kleiner \cite{burago-kleiner}. 
For a finitely generated group $G$, it is not known if $G$ and $G\times \ZZ_2$ can fail to be bi-Lipschitz.
(The groups are assumed to be equipped with the word metric.)


\begin{wrapfigure}{r}{35 mm}
\vskip-4mm
\centering
\includegraphics{mppics/pic-612}
\end{wrapfigure}

\parbf{Piecewise Euclidean quotient.}
Note that the group $\Gamma$ is the holonomy group of the quotient space $P=\RR^m/\Gamma$.
More precisely, one can identify $\RR^m$ with the tangent space to a regular point $x_0$ of $P$ in such a way that
for any $\gamma\in\Gamma$ there is a loop $\ell$ based at $x_0$ that runs in the regular locus of $P$ and has the holonomy~$\gamma$.

Choose $\gamma$ and $\ell$ as above.
Since $P$ is simply-connected, we can shrink $\ell$ by a disk.
By the general position argument, we can assume that the disk 
only passes thru simplices of codimension $0$, $1$ and $2$
and intersects the simplices of codimension $2$ transversely.

In other words, $\ell$ can be presented as a product of 
loops such that each loop goes around a single simplex of codimension $2$ and comes back.
The holonomy for each of these loops is a rotation around a hyperplane.
Hence the result follows.
\qeds

The converse of the problem also holds;
it was proved by Christian Lange \cite{lange};
his proof is based on earlier results of 
Marina Mikhailova \cite{mikhailova}.

Note that the cone over the spherical suspension over the Poincaré sphere is homeomorphic to $\RR^5$ and it is the quotient of $\RR^5$ by the binary icosahedral group, which is a subgroup of $\SO(5)$ of order 120. 
Therefore, 
if one replaces ``piecewise linear homeomorphism'' with ``homeomorphism'' in the formulation, 
then the answer will be different; 
a complete classification of such actions is given in \cite{lange}.

\parbf{Subgroups of a free group.}
The proof exploits the fact that free groups are fundamental groups of graphs.

\begin{wrapfigure}{o}{38 mm}
\vskip-4mm
\centering
\includegraphics{mppics/pic-614}
\end{wrapfigure}

\medskip

Let $F$ be a free group and $G$ a finitely generated subgroup in $F$.
We need to show that $G$ is an intersection of subgroups of finite index in $F$.
Without loss of generality, we can assume that $F$ has a finite number of generators, denote it by $m$.

Let $W$ be the wedge sum of $m$ circles so that $\pi_1(W,p)\z=F$.
Equip $W$ with the length-metric such that each circle has unit length.

Pass to the metric covering $\tilde W$ of $W$ 
such that  $\pi_1(\tilde W,\tilde p)=G$ 
for a lift $\tilde p$ of $p$.

Choose a sufficiently large integer $n$ and consider the doubling of the closed ball $\bar B(\tilde p,n+\frac12)$ along  its boundary.
Let us denote the obtained doubling by $Z_n$ and set $G_n=\pi_1(Z_n,\tilde p)$.

Note that $Z_n$ is a metric covering of $W$;
this allows us to consider $G_n$ as a subgroup of $F$.
By construction, $Z_n$ is compact;
therefore $G_n$ has a finite index in $F$.

It remains to show that 
\[G=\bigcap_{n>k} G_n,\]
where $k$ is the maximal word length in the generating set of $G$.
\qeds

Originally the problem was solved by Marshall Hall \cite{hall1,hall2,burns}.
Our proof is close to the solution of John Stallings \cite{stallings,wilton}.
Note that the statement does not hold for infinitely generated subgroups. 

The same idea can be used to solve many other problems; here are some examples.

\begin{pr}
Show that a subgroup of a free group is free.
\end{pr}

\begin{pr}
 Show that two elements $u$ and $v$ of a free group commute 
if and only if they are both powers of
an element $w$.
\end{pr}

\parbf{Short generators.}
Choose a length-minimizing loop $\gamma$ that represents a given element $a\in\pi_1M$.

\begin{wrapfigure}{o}{21 mm}
\vskip0mm
\centering
\includegraphics{mppics/pic-616}
\vskip1mm
\end{wrapfigure}

Choose $\eps>0$.
Represent $\gamma$ 
as a concatenation of paths
$\gamma=\gamma_1{*}\dots{*}\gamma_n$
such that
\[\length\gamma_i<\eps\] 
for each $i$.

Denote by $p=p_0,p_1,\dots, p_n=p$ the endpoints of these arcs.
Connect $p$ with $p_i$ by a minimizing geodesic $\sigma_i$.
Note that $\gamma$ is homotopic to a product of loops
\[\alpha_i=\sigma_{i-1}{*}\gamma_i{*}\bar\sigma_{i},\]
where $\bar\sigma_{i}$ denotes the path $\sigma_{i}$ traveled backward.
In particular,
\[\length \alpha_i<2\cdot\diam M+\eps \]
for each $i$.

Note that given $\ell>0$, there are only finitely many elements of the fundamental group that can be realized by loops at $p$ with length shorter than $\ell$.
It follows that for the right choice of $\eps>0$, 
any loop $\alpha_i$ is homotopic to a loop of length at most $2\cdot\diam M$.
Hence the result follows.
\qeds

The statement is due to 
Mikhael Gromov \cite[Proposition 3.22 in][]{gromov-MetStr}.

\parbf{Number of generators.}
Consider the universal Riemannian covering $\tilde M$ of $M$.
Note that $\tilde M$ is non-negatively curved and
$\pi_1M$ acts by isometries on $\tilde M$.

Choose $p\in \tilde M$.
Given  $a\in \pi_1M$,
set 
\[|a|=|p- a\cdot p|_{\tilde M}.\]

Consider the so-called \index{short basis}\emph{short basis} in $\pi_1M$;
that is, a sequence of elements $a_1,a_2,\dots{} \in\pi_1M$ defined in the following way:
\begin{enumerate}[(i)]
\item choose $a_1\in\pi_1M$ so that $|a_1|$ takes the minimal value,
\item choose $a_2\in\pi_1M\setminus\langle a_1 \rangle$ so that $|a_2|$ takes the minimal value,
\item choose $a_3\in\pi_1M\setminus\langle a_1,a_2 \rangle$ so that $|a_3|$ takes the minimal value,
and so on.
\end{enumerate}

Note that the sequence terminates at the $n$-th step 
if 
$a_1,\dots,a_n$  generate $\pi_1M$.
By construction, we have
\begin{align*}
|a_j\cdot a_i^{-1}|&\ge |a_j|\ge |a_i|
\intertext{for any $j>i$. 
Set $p_i=a_i\cdot p$.
Note that}
|p_j-p_i|_{\tilde M}
&=|a_j\cdot a_i^{-1}|\ge
\\
&\ge |a_j|=
\\
&=|p_j-p|_{\tilde M}\ge
\\
&\ge|a_i|=
\\
&=|p_i-p|_{\tilde M}.
\intertext{By the Toponogov comparison theorem we get}
\measuredangle \hinge p{p_i}{p_j}&\ge \tfrac\pi3.
\end{align*}
That is, the directions from $p$ to all $p_i$ make an angle of at least $\tfrac\pi3$ with each other.

Therefore the number of points $p_i$ can be bounded in terms of the dimension of $M$.
Hence the result follows.
\qeds

The \emph{short-basis construction}, as well as the result above are due to Mikhael Gromov \cite{gromov-almost-flat}.

\parbf{An equation in a Lie group.} 
We will assume that $G$ is equipped with a bi-invariant metric.
In particular, geodesics starting at the identity element $e\in G$ are given by homomorphisms $\RR\to G$.

Consider the map $\phi\:G\to G$ defined by
\[\phi(x)=x\cdot g_1\cdot x\cdot g_2\cdots x\cdot g_n.\]
We need to show that $\phi$ is onto.
Note that it is sufficient to show that $\phi$ has a non-zero degree.

The map $\phi$ is homotopic to the map $\psi\:x\mapsto x^n$.
Therefore it is sufficient to show that
\[\deg \psi\ne 0\leqno({*})\]

Note that the claim $({*})$ follows from $({*}{*})$.
\begin{cl}{$({*}{*})$} For any $x\in G$ the differential 
 \[d_x\psi\:\T_xG\to \T_{x^n}G\] 
does not revert orientation.
\end{cl}

Indeed, connect $e$ to a given point $y\in G$ by a geodesic path $\gamma$, so $\gamma(0)=e$ and $\gamma(1)=y$.
Since $\gamma\:\RR\to G$ is a homomorphism,
$\psi(x)=y$ for $x=\gamma(\tfrac1n)$.
In particular, the inverse image $\psi^{-1}\{y\}$ is nonempty for any $y\in G$.

By $({*}{*})$, for a regular value $y$, each point in the  inverse image $\psi^{-1}\{y\}$ contributes $1$ to the degree of $\psi$. 
Hence $({*})$ follows.

It remains to prove $({*}{*})$.
Given an element $g\in G$, denote by $L_g$ and $R_g$ the left and right shift $G\to G$ respectively;
that is, $L_g(x)\z=g\cdot x$ and $R_g(x)=x\cdot g$.

Identify the tangent spaces $\T_xG$ and $\T_{x^n}G$ with the Lie algebra $\mathfrak{g}\z=\T_eG$
using $d{R_x}\:\mathfrak{g}\to \T_xG$ and $d{R_x^n}\:\mathfrak{g}\to \T_{x^n}G$ respectively.
Then for any $V\in \mathfrak{g}$, we have
\[d_x\psi(V)=V+\Ad_x(V)+\dots+\Ad_x^{n-1}(V),\]
where $\Ad_x=d_e(L_x\circ R_{x^{-1}})\:\mathfrak{g}\to \mathfrak{g}$. 
Since the metric on $G$ is bi-invariant, $\Ad_x$ is an isometry of $\mathfrak{g}$.
It remains to note that the linear transformation
\[V\mapsto V+T(V)+\dots+T^{n-1}(V)\]
cannot revert orientation for any isometric linear transformation $T$ of the Euclidean space.
The last statement is an exercise in linear algebra.
\qeds

The idea of this solution is due to Murray Gerstenhaber and Oscar Rothaus 
\cite{gerstenhaber-rothaus}.
In fact, the degree of $g$ is $n^k$, where $k$ is the rank of~$G$ \cite{hopf}.

\parbf{Quotient of Hilbert space.}
We consider $\mathbb{S}^3$ as the set of unit quaternions;
in particular, it has a group structure.

Let $\HH$ be the set of paths of \emph{class $W^{1,2}$} in $\mathbb{S}^3$ starting at the identity element $e$;
that is, the path's velocity is square-integrable.
The pointwise multiplication of paths defines a group structure on $\HH$.
Denote by $\Omega$ the subset of all loops in $\HH$.

It remains to equip $\HH$ with the structure of a Hilbert space so that 
the right action of $\Omega$ on $\HH$ is isometric and the quotient is isometric to~$\mathbb{S}^3$.

\medskip

We will prove the statement for any connected Lie group $G$ with a bi-invariant metric; in particular, for $G=\mathbb{S}^3$.
Denote by $\mathfrak{g}=\T_eG$ the Lie algebra of $G$.
Equip $G$ with a bi-invariant metric, and let $\langle{*},{*}\rangle_{\mathfrak{g}}$ be the corresponding scalar product in $\mathfrak{g}$.

Consider the Hilbert space $\HH$ of all $L^2$-functions $f\:[0,1]\to\mathfrak{g}$ with the scalar product defined by
\[\<f,g\>=\int\limits_0^1\<f(t),g(t)\>_{\mathfrak{g}}\cdot dt.\]

\parit{Construction of the quotient map $\phi\:\HH\to G$.}
Given $v\in \mathfrak{g}$, denote by $\tilde v$ the corresponding right-invariant tangent field on $G$.

Given $f\:[0,1]\to \mathfrak{g}$ in $\HH$,
consider the path 
\[\Gamma_f\:[0,1]\to G\] 
with 
$\Gamma_f(0)=1$ and $\Gamma_f'(t)=\tilde f(t)$ for any $t$.

The map $\phi\:\HH\to G$ is the evaluation map $\phi\:f\mapsto \Gamma_f(1)$.
Since $G$ is connected, $\phi$ is onto.

\parit{Group structure on $\HH$.}
Note that the functional $f\mapsto \Gamma_f$ is an injective map from $\HH$ to the space of paths in $G$ starting at $e$.

Given $\alpha\in G$, we denote by $\Ad_\alpha\:\mathfrak{g}\to \mathfrak{g}$ its the adjoint transformation;
that is, $\Ad_\alpha=d_e\Inn_\alpha$, where $\Inn_\alpha\:x\mapsto \alpha\cdot x\cdot \alpha^{-1}$ is the inner automorphism of $G$.
Note that $\Ad_\alpha$ preserves the scalar product on~$\mathfrak{g}$.

Consider the multiplication $\star$ on $\HH$ defined by
\[(h\star f)(t)=h(t)+\Ad_{\Gamma_h(t)}[f(t)].\leqno({*})\]

Note that 
\[\Gamma_{h\star  f}(t)=\Gamma_h(t)\cdot \Gamma_f(t)\]
for any $t\in[0,1]$.
In particular, $(\HH,\star )$ is a group with neutral element~$0$. 

From $({*})$, we get
\[(h\star f)(t)-(h\star g)(t)=\Ad_{\Gamma_h(t)}(f(t)-g(t))\]
and therefore
\[|(f\star h)(t)-(g\star h)(t)|=|f(t)-g(t)|\]
for any $t$.
It follows that for any fixed $h$,
the transformation $f\mapsto h\star f$ is an affine isometry of $\HH$.

The set $\Omega=\phi^{-1}\{e\}$ is a subgroup of $(\HH,\star)$;
it can be viewed as the group of $W^{1,2}$-loops in $G$.
It remains to note that $\phi\:\HH\to G$ is the quotient map for the right action of $\Omega$ on $\HH$.
\qeds

\parit{Alternative solution.} Again, we will prove the statement for any connected Lie group $G$ with a bi-invariant metric.

Denote by $G^n$ the direct product of $n$ copies of $G$.
Consider the map $\phi_n\:G^n\to G$ defined by
\[\phi_n\:(\alpha_1,\dots,\alpha_n)\mapsto \alpha_1\cdots\alpha_n.\]
Note that $\phi_n$ is the quotient map for the $G^{n-1}$-action on $G^n$ defined by
\[(\beta_1,\dots,\beta_{n-1})\cdot(\alpha_1,\dots,\alpha_n)=(\alpha_1\cdot\beta_1^{-1},\beta_1\cdot\alpha_2\cdot\beta_2^{-1},\dots,\beta_{n-1}\cdot\alpha_n).\]

Denote by $\rho_n$ the product metric on $G^n$ rescaled with factor $\sqrt{n}$.
Note that the quotient $(G^n,\rho_n)/G^{n-1}$ is isometric to $G=(G,\rho_1)$.

As $n\to\infty$ the curvature of $(G^n,\rho_n)$ converges to zero and its injectivity radius goes to infinity.
Therefore passing to the ultra-limit of $G^n$ as $n\to\infty$ we get a Hilbert space.
It remains to observe that the limit action has the required property.
\qeds

This construction is given by Chuu-Lian Terng and Gudlaugur Thorbergsson \cite[see section 4 in][]{terng-thorbergsson};
the alternative solution was suggested by Alexander Lytchak.

Instead of the group $\Omega$, 
one could consider the subgroup $\Omega_H$ of paths $\gamma\:[0,1]\to G$ such that the pair $(\gamma(0),\gamma(1))$ belongs to a given subgroup $H<G\times G$.
In this case, the quotient $\HH/\Omega_H$ is isometric to the \emph{double quotient} $G/\!\!/H$;
that is, the quotient of the action on $G$ defined by $(h_1,h_2)\cdot g=h_1\cdot g\cdot h_2^{-1}$ for $(h_1,h_2)\in H<G\times G$.

The following question is open.

\begin{pr} Suppose $R$ is a compact simply-connected Riemannian manifold that is isometric to a quotient of the Hilbert space by a group of isometries (or more generally $R$ is the target of Riemannian submersion from a Hilbert space).
Is it true that $R$ is isometric to a double quotient? That is, is it true that $R$ is a quotient of compact Lie group $G$ by a group of isometries?
 
\end{pr}

\csname @openrightfalse\endcsname
\chapter{Topology}

In this chapter, we consider geometrical problems with strong topological flavor.
A typical introductory course in topology, say \cite{kosniowski},
contains all the necessary material.

\subsection*{Isotopy}\label{Isotopy}

Recall that an isotopy is a continuous one-parameter family of embeddings.

\begin{pr}
Let $K_1$ and $K_2$ be homeomorphic closed subsets of the coordinate subspace $\RR^m$ in $\RR^{2\cdot m}$.
Show that there is a homeomorphism 
\[h\:\RR^{2\cdot m}\z\to \RR^{2\cdot m}\] 
such that $K_2=h(K_1)$.
Moreover, $h$ can be chosen to be isotopic to the identity map.
\end{pr}

\parit{Semisolution.}
Choose a homeomorphism $\phi\:K_1\to K_2$.

By the Tietze extension theorem,
the homeomorphisms $\phi\:K_1\z\to K_2$ and $\phi^{-1}\:K_2\z\to K_1$ can be extended to continuous maps;
denote these maps by $f\:\RR^m\z\to \RR^m$ and $g\:\RR^m\z\to \RR^m$ respectively.

{

\begin{wrapfigure}{o}{56 mm}
\vskip-4mm
\centering
\includegraphics{mppics/pic-702}
\end{wrapfigure}

Consider homeomorphisms
$h_1$, $h_2$, and $h_3$ of $\RR^m\times\RR^m$
defined in the following way:
\begin{align*}
h_1(x,y)&=(x,y+f(x)),
\\
h_2(x,y)&=(x-g(y),y),
\\ 
h_3(x,y)&=(y,-x).
\end{align*}

}

It remains to observe that each homeomorphism $h_i$ is isotopic to the identity map and
\[K_2=h_3\circ h_2\circ h_1(K_1).\qedsin\]

This construction is due to Victor Klee \cite{klee} and it is called \emph{Klee's trick}.
This trick is used in the five-line proof of the Jordan separation theorem by Patrick Doyle \cite{doyle};
a proof of the separation theorem for embeddings $\mathbb{S}^n\hookrightarrow\mathbb{S}^{n+1}$
can be given using the same idea \cite{cohen}. 

The problem ``Monotonic homotopy'' on page \pageref{mono-homotopy} looks similar.

\subsection*{Immersed disks}\label{Immersed disks}

Two immersions $f_1$ and $f_2$ of the disk $\DD$ into the plane will be called {}\emph{essentially different} 
if there is no diffeomorphism $h\:\DD\z\to \DD$ such that
$f_1=f_2\circ h$.

\begin{pr} 
Construct two essentially different smooth immersions of the disk 
into the plane that coincide near the boundary. 
\end{pr}

\subsection*{Positive Dehn twist}\label{Positive Dehn twist} 

\begin{wrapfigure}{r}{43 mm}
\begin{lpic}[t(-4 mm),b(0 mm),r(0 mm),l(0 mm)]{asy/dehn-twist()}
\lbl[b]{20,10;$\xrightarrow{\ h\ }$}
\end{lpic} 
\end{wrapfigure}

Let $\Sigma$ be a surface and 
\[\gamma\:\RR/\ZZ\z\to\Sigma\] 
be a non-contractible closed simple curve.
Let $U_\gamma$ be a neighborhood of $\gamma$ that admits a parametrization 
\[\iota\:\RR/\ZZ\times (0,1)\to U_\gamma.\]
A \index{Dehn twist}\emph{Dehn twist} along $\gamma$ is a homeomorphism $h\:\Sigma\z\to\Sigma$
that is the identity outside of $U_\gamma$ and 
such that
\[\iota^{-1}\circ h\circ \iota\:(x,y)\mapsto(x+y,y).\]

If $\Sigma$ is oriented 
and $\iota$ is orientation preserving,
then the Dehn twist described above is called {}\emph{positive}.

\begin{pr}
Let $\Sigma$ be a compact oriented surface with nonempty boundary.
Prove that any composition of positive Dehn twists of $\Sigma$ is not homotopic to the identity relative to the boundary.

In other words, any product of positive Dehn twists represents a nontrivial class in the mapping class group of $\Sigma$.
\end{pr}

\subsection*{Conic neighborhood}
\label{Conic neighborhood}

Let $p$ be a point in a topological space $X$.
We say that an open neighborhood $U\ni p$ is \index{conic neighborhood}\emph{conic}
if there is a homeomorphism from a cone
to $U$ that sends the vertex to $p$.

\begin{pr}  
Show that any two conic neighborhoods of one point are homeomorphic to each other.
\end{pr}

Note that two cones $\mathop{\rm Cone}(\Sigma_1)$ and $\mathop{\rm Cone}(\Sigma_2)$ might be homeomorphic while $\Sigma_1$ and $\Sigma_2$ are not;
the existence of such examples follows from the double suspension theorem.

\subsection*{Unknots\easy}\label{No knots}

\begin{pr}
Prove that the set of smooth embeddings $f\:\mathbb{S}^1\z\to\RR^3$ equipped with the $C^0$-topology 
forms a connected space.
\end{pr}

\subsection*{Stabilization}\label{Simple stabilization}

\begin{pr}
Construct two compact subsets $K_1, K_2\subset\RR^2$ such that
$K_1$ is not homeomorphic to $K_2$, but $K_1\times[0,1]$ is homeomorphic to $K_2\z\times[0,1]$.
\end{pr}

\subsection*{Homeomorphism of a cube}\label{Homeomorphism of cube}

\begin{pr}
Let $\square$ be a cube in $\RR^m$
and $h\:\square\z\to\square$ be
a homeomorphism that sends each face of $\square$ to itself.
Extend $h$ to a homeomorphism $f\:\RR^m\z\to\RR^m$ that coincides with the identity map outside of a bounded set.    
\end{pr}

\subsection*{Finite topological space\easy}\label{Finite topological space}

\begin{pr}
Given a finite topological space $F$, 
construct a finite simplicial complex $K$
that admits a weak homotopy equivalence  $K\to F$. 
\end{pr}

\subsection*{Dense homeomorphism\easy}\label{Dense homeomorphism}

\begin{pr}
Denote by $\mathcal{H}$ be the set of all orientation preserving homeomorphisms $\mathbb {S}^2\z\to\mathbb {S}^2$ 
equipped with the $C^0$-metric.
Show that there is a homeomorphism $h\in \mathcal{H}$ such that its conjugations $a\circ h\circ a^{-1}$ for all $a\in\mathcal{H}$ form a dense set in $\mathcal{H}$.
 
\end{pr}

\subsection*{Simple path\easy}
\label{Simple path}

\begin{pr}
Let $p$ and $q$ be distinct points in a Hausdorff topological space $X$.
Assume that $p$ and $q$ are connected by a path.
Show that they can be connected by a simple path;
that is, there is an injective continuous map $\beta\:[0,1]\z\to X$
such that $\beta(0)=p$ and $\beta(1)=q$.
\end{pr}

(This statement might be intuitively obvious, but its proof is not simple.)

\subsection*{Path on a surface\easy}
\label{Path on a surface}

\begin{pr}
Show that any path with distinct ends in a surface is homotopic (relative to the ends) to a simple path.  
\end{pr}

{

\begin{wrapfigure}{r}{29 mm}
\vskip-4mm
\centering
\includegraphics{mppics/pic-704}
\bigskip
\includegraphics{mppics/pic-706}
\end{wrapfigure}

\section*{Semisolutions}

\parbf{Immersed disks.}
Both circles on the picture bound essentially different disks.

\medskip

On the first diagram, the dashed lines and the solid lines together bound three embedded disks;
gluing these disks along the dashed lines gives the first immersion.
The reflection of this immersion across the vertical line of symmetry gives another essentially different immersion.
\qeds

}

It is a good exercise to count the essentially different disks in the second example. 
(The answer is 5.) 

{

\begin{wrapfigure}{r}{30 mm}
\vskip-7mm
\centering
\includegraphics{mppics/pic-708}
\end{wrapfigure}

The existence of examples of that type is generally attributed to John Milnor \cite{bennequin}.

An easier problem would be to construct two essentially different immersions of annuli with the same boundary curves; a solution is shown on the picture \cite[for more details and references see][]{eppstein-mumford}.

}

\parbf{Positive Dehn twist.}
Consider the universal covering 
$f\:\tilde\Sigma\z\to\Sigma$.
The surface $\tilde \Sigma$ has a boundary and it comes with the orientation induced from $\Sigma$.

Choose a point $x_0$ on the boundary $\partial \tilde \Sigma$.
Given two other points $y$ and $z$ in $\partial \tilde \Sigma$, we will write
$z\succ y$ if $y$ lies on the right side from a simple curve from $x_0$ to $z$ in $\tilde\Sigma$.
Note that  $\succ $ defines a linear order on $\partial\tilde\Sigma\setminus\{x_0\}$.
We will write $z \succeq y$ 
if $z\succ y$ or $z=y$.

{

\begin{wrapfigure}{o}{27 mm}
\vskip-0mm
\centering
\includegraphics{mppics/pic-710}
\end{wrapfigure}

Note that any homeomorphism $h\:\Sigma\to\Sigma$ identical on the boundary
lifts to the unique homeomorphism $\tilde h\:\tilde \Sigma\to\tilde\Sigma$ 
such that $\tilde h(x_0)\z=x_0$.
The following claim is the key step in the proof:

}

\begin{cl}{$({*})$} 
If $h$ is a positive Dehn twist along a closed curve $\gamma$,
then $y\succeq \tilde h(y)$ for any $y\in\partial\tilde\Sigma\setminus\{x_0\}$
and $y_0\succ\tilde h(y_0)$ for some $y_0\in\partial\tilde\Sigma\setminus\{x_0\}$.
\end{cl}

Note that the problem follows from~$({*})$.
Indeed, the property in $({*})$ is a homotopy invariant 
and it survives under compositions of maps.

\medskip

If $\Sigma$ is not an annulus,
then by the uniformization theorem we can assume that $\Sigma$ has a  hyperbolic metric with geodesic boundary; 
the lifted metric on $\tilde\Sigma$ has the same properties.
Furthermore, we can assume that (1) $\gamma$ is a closed geodesic,
(2) the parametrization $\iota\:\RR/\ZZ\z\times (0,1)\to U_\gamma$ from the definition of Dehn twist is rotationally symmetric 
and (3) for any $u\in \RR/\ZZ$ the arc $\iota(u\times (0,1))$ is a geodesic perpendicular to~$\gamma$. 

Consider the polar coordinates $(\phi,\rho)$ on $\tilde\Sigma$ with the origin at~$x_0$;
since $x_0$ lies on the boundary, the angle coordinate $\phi$ is defined in $[0,\pi]$. 
By construction of the Dehn twist, we get 
\[\phi(x)\ge \phi\circ\tilde h(x)\]
for any $x\ne x_0$,
and if the geodesic $[x_0x]$ crosses $f^{-1}(U_\gamma)$, then 
\[\phi(x)> \phi\circ\tilde h(x).\]
In particular, if $x$ lies on the boundary then $\tilde h(x)$ lies on the right of the geodesic $[x_0x]$; hence the claim $({*})$ follows. 

\begin{figure}[!ht]
\vskip0mm
\centering
\includegraphics{mppics/pic-712}
\end{figure}

If $\Sigma$ is an annulus, then the same argument works except we have to choose a flat metric on $\Sigma$.
In this case, $\tilde \Sigma$ is a strip between two parallel lines in the plane, see the diagram.
\qeds

Note that if the surface has an empty boundary, then the answer is different:

\begin{pr}
Construct a composition of positive Dehn twists on a compact oriented surface without boundary that is homotopic to the identity. 
\end{pr}

I learned the problem from Rostislav Matveyev.
The described construction was given by Hamish Short and Bert Wiest \cite{short-wiest} and attributed to William Thurston.

It turns out that there is no upper bound on the length of the product.
Namely, there is a compact oriented surface $\Sigma$ with a nonempty boundary and a homeomorphism $f\:\Sigma\to \Sigma$ that does not move boundary points such that $f$ is homotopic relative to the boundary to an arbitrary long product of positive Dehn twists.
Such an example was constructed by Refik İnanç Baykur and Jeremy Van Horn-Morris \cite{baykur-vanhornmorris};
see also \cite{baykur-monden-vanhornmorris}.

\parbf{Conic neighborhood.}
Let $V$ and $W$ be two conic neighborhoods of~$p$.
Without loss of generality, we may assume that $V\Subset W$;
that is, the closure of $V$ lies in $W$.

We will need to construct a sequence of embeddings $f_n\:V\to W$
such that 
\begin{itemize}
\item 
For any compact set $K\subset V$ 
there is a positive integer $n=n_K$ such that 
$f_n(k)=f_m(k)$ for any $k\in K$ and $m, n \ge n_K$.
\item For any point $w\in W$ there is a point $v\in V$ such that $f_n(v)=w$ for all large $n$.
\end{itemize}

Note that once such a sequence is constructed, $f\:V\to W$ defined by $f(v)=f_n(v)$ for all large values of $n$ gives the needed homeomorphism.

The sequence $f_n$ can be constructed recursively
\[f_{n+1}=\Psi_n\circ f_n\circ \Phi_n,\]
where $\Phi_n\:V\to V$ 
and $\Psi_n\:W\to W$ 
are homeomorphisms
of the form 
\[\Phi_n(x)=\phi_n(x)\ast x\quad \text{and}\quad \Phi_n(x)=\psi_n(x)\star x,\]
where $\phi_n\:V\to \RR_{\ge 0}$, $\psi_n\:W\to \RR_{\ge 0}$ are suitable continuous functions;
``$\ast$'' and ``$\star$'' denote the {}\emph{multiplication} in the cone structures of $V$ and $W$ respectively.\qeds

The problem is due to Kyung Whan Kwun \cite{kwun}.

\parbf{Unknots.}

\begin{figure}[!ht]
\vskip0mm
\centering
\includegraphics{mppics/pic-714}
\end{figure}

Observe that it is possible to draw an arbitrary tight knot 
while keeping it smoothly embedded at all times including the last moment.\qeds

This problem was suggested by Greg Kuperberg.


\parbf{Stabilization.}
The example can be guessed from the diagram.

\begin{wrapfigure}[6]{r}{50 mm}
\vskip-0mm
\centering
\includegraphics{mppics/pic-716}
\end{wrapfigure}

The two sets $K_1$ and $K_2$ are subspaces of the plane, 
each one being a closed annulus with two attached line segments.
In $K_1$ one segment is attached from the inside and another from the outside; 
in $K_2$ both segments are attached from the outside.

The product spaces $K_1\times[0,1]$ and $K_2\times[0,1]$ are solid tori with attached rectangles.
A homeomorphism $K_1\times[0,1]\z\to K_2\times[0,1]$ can be constructed by twisting a part of one solid torus.

To prove the nonexistence of a homeomorphism $K_1\to K_2$ consider the sets of cut points $V_i\subset K_i$ and the sets $W_i\subset K_i$ of points that admit a punctured simply-connected neighborhood.
Note that the set $V_i$ is the union of the attached line segments 
and $W_i$ is the boundary of the annulus without points where the segments are attached.
Note that $V_i\cup W_i=\partial K_i$;
in particular, a homeomorphism $K_1\to K_2$ (if it exists) sends $\partial K_1$ to $\partial K_2$.

Finally, note that each $\partial K_i$ has two connected components and 
$V_1$ intersects both components of $\partial K_1$
while $V_2$ lies in one component of $\partial K_2$.
Hence $K_1\ncong K_2$.
\qeds

It should be an old puzzle;
I learned it from Maria Goluzina around 1988.

\parbf{Homeomorphism of a cube.}
Let us extend the homeomorphism $h$ to $\RR^m$ by reflecting the cube across its facets.
We get a homeomorphism $\tilde h\:\RR^m\to\RR^m$ such that $\tilde h(x)=h(x)$ for any $x\in\square$ and 
\[\tilde h\circ\gamma=\gamma\circ \tilde h,\]
where $\gamma$ is any reflection with respect to the facets of the cube.

Without loss of generality, we may assume that the cube $\square$ is inscribed in the unit sphere centered at the origin of $\RR^m$.
In this case, $\tilde h$ has \index{displacement}\emph{displacement} at most $2$;
that is, 
\[|\tilde h(x)-x|\le 2\]
for any $x\in\RR^m$.

Choose a smooth increasing function $\phi\:\RR_{\ge 0} \to\RR$ such that
$\phi(r)\z=r$
for $r\le 1$ and $\phi(r)\to 2$ as $r\to\infty$.

Equip $\RR^m$ with polar coordinates $(u,r)$, 
where $u\in\mathbb{S}^{m-1}$, $r\z\ge 0$.
Consider a homeomorphism $\Phi$ from $\RR^m$ to an open ball of radius $2$
defined by 
\[\Phi(u,r)\z=(u,\phi(r)).\]

\begin{wrapfigure}{o}{41 mm}
\vskip-0mm
\centering
\includegraphics{mppics/pic-717}
\vskip0mm
\end{wrapfigure}

Set 
\[
f(x)=\left[
\begin{aligned}
&x&&\text{if}\ |x|\ge 2,
\\
&\Phi\circ \tilde h \circ \Phi^{-1}(x)&&\text{if}\ |x|< 2.
\end{aligned}
\right.
\]
It remains to observe that $f\:\RR^m\to\RR^m$ is a solution.
\qeds

This problem is stripped from the proof of Robion Kirby \cite{kirby}.
The condition that each face is mapped to itself can be removed and 
instead of homeomorphism one could take any embedding close to the identity.

An interesting twist to this idea was given by Dennis Sullivan [see \ncite{sullivan} and also 5.10 in \ncite{tukia-vaisala}].
Instead of the discrete group of motions of the Euclidean space,
he uses a discrete group of motions of the hyperbolic space in the conformal disk model.
To see the idea, note that the construction of $\tilde h$ can be done for a Coxeter polytope in the hyperbolic space instead of a cube.%
\footnote{By Vinberg's theorem \cite{vinberg, vinberg-strong} hyperbolic space of large dimension has no Coxeter polytopes, but the idea works after some modifications.}
Then the constructed map $\tilde h$
coincides with the identity on the absolute and therefore the last ``shrinking'' step in the proof above is not needed.
Moreover, 
if the original homeomorphism is bi-Lipschitz,
then the Sullivan construction produces a bi-Lipschitz homeomorphism ---
this is its main advantage.

\parbf{Finite topological space.}
Given a point $p\in F$,
denote by $O_p$ the minimal open set in $F$ containing $p$. 
Note that we can assume that $F$ is a connected $T_0$-space;
in particular, $O_p=O_q$ if and only if $p=q$.

Let us write $p\preccurlyeq q$ 
if $O_p\subset O_q$.
Evidently, $\preccurlyeq$ is a partial order on~$F$.

Let us construct a simplicial complex $K$ 
by taking $F$ as the set of vertices
and declaring a collection of vertices to be a simplex 
if it can be linearly ordered with respect to $\preccurlyeq$.

Given $k\in K$,
consider the minimal simplex $(f_0,\dots,f_m)\ni k$;
we can assume that $f_0\preccurlyeq \dots\preccurlyeq f_m$.
Set $h\:k\mapsto f_0$;
it defines a map $K\to F$.

It remains to check that $h$ is continuous 
and induces isomorphisms for all the homotopy groups.
\qeds

In a similar fashion, one can construct a finite topological space $F$ for any given simplicial complex $K$ 
such that 
there is a weak homotopy equivalence $K\to F$.
Both constructions are due to Pavel Alexandrov
\cite{alexandrov-finite,mccord}.

\parbf{Dense homeomorphism.}
Note that there is a countable set of homeomorphisms $h_1,h_2,\dots$ that is dense in $\mathcal{H}$
such that
each $h_n$ fixes all the points outside an open round disk, say $D_n$.

Choose a countable disjoint collection of round disks $D_n'$.
Consider the homeomorphism $h\:\mathbb S^2\to \mathbb S^2$
that fixes all the points outside of $\bigcup_nD'_n$ and
for each $n$,
the restriction $h|_{D_n'}$ is conjugate to $h_n|_{D_n}$.

Note that for large $n$, the homeomorphism $h$ is conjugate to a homeomorphism close to $h_n$.
Therefore $h$ is a solution.
\qeds

The problem was mentioned by Frederic Le Rox \cite{rox} on a problem section at a conference in Oberwolfach, 
where he also conjectured that this is not true for the area-preserving homeomorphisms.
An affirmative answer to this conjecture was given by Daniel Dore, Andrew Hanlon, and Sobhan Seyfaddini 
\cite{dore-hanlon,seyfaddini}.
In particular, it implies the following seemingly evident but nontrivial statement.

\begin{pr}
Given $\eps>0$, there is $\delta>0$ such that 
\[\Omega\cap h(\Omega)\ne\emptyset\]
for any topological disk $\Omega\subset \mathbb{S}^2$ with area at least $\eps$
and 
any area-preserving homeomorphism $h\:\mathbb{S}^2\to\mathbb{S}^2$ with displacement at most $\delta$;
that is, such that $|h(x)-x|_{\mathbb{S}^2}<\delta$ for any $x\in \mathbb{S}^2$. 
\end{pr}

\parbf{Simple path.}
We will give two solutions, the first one is elementary and the second one is involved. 

\parit{First solution.}
Let $\alpha$ be a path connecting $p$ to $q$.

Passing to a subinterval if necessary,
we can assume that $\alpha(t)\ne p,q$ for $t\ne0,1$.

An open set $\Omega$ in $(0,1)$ will be called {}\emph{suitable} if, for any connected component $(a,b)$ of $\Omega$, we have $\alpha(a)=\alpha(b)$.
Since the union of nested suitable sets is suitable, we can find a maximal suitable set $\hat \Omega$.

Define $\beta(t)=\alpha(a)$ for any $t$ in a connected component $(a,b)\subset\Omega$.
Note that $\beta$ is continuous and monotonic;
that is, for any $x\in [0,1]$ the set $\beta^{-1}\{\beta(x)\}$ is connected.

It remains to reparametrize $\beta$ to make it injective.
In other words, we need to construct a non-decreasing surjective function $\tau\:[0,1]\z\to[0,1]$ such that 
$\tau(t_1)=\tau(t_2)$ if and only if there is a connected component $(a,b)$ such that $t_1,t_2\z\in [a,b]$.
The construction is similar to the construction of the devil's staircase.
\qeds

\parit{Second solution.}
Note that one can assume that $X$ coincides with the image of $\alpha$.
In particular, $X$ is a connected locally connected compact Hausdorff space.

Any such space admits a length-metric.
This statement is not at all trivial;
it was conjectured by Karl Menger \cite{menger}
and proved independently 
by R.~H.~Bing  \cite{bing-length-0, bing-length-1} 
and Edwin Moise \cite{moise}.

It remains to consider a geodesic path from $p$ to $q$.
\qeds

The problem was inspired by a lemma 
proved by 
Alexander Lytchak
and Stefan Wenger \cite[see 7.13 in][]{lytchak-wenger}.

\parbf{Path on a surface.}
Denote the surface by $\Sigma$; assume that the path runs from $p$ to $q$.
The following picture suggests an idea for an induction on the number of self-crossings.

\begin{figure}[!ht]
\vskip0mm
\centering
\includegraphics{mppics/pic-718}
\end{figure}

To do the proof formally,
let us present the path as a concatenation $\alpha*\beta$ of two paths  such that $\alpha$ is simple
and $\beta$ does not pass thru~$p$.
We can assume that $\beta\:[0,1]\to \Sigma$ is smooth.

Choose a smooth time-dependent vector field $V_t$ on $\Sigma$ such that
\[V_t(\beta(t))=\beta'(t)\quad\text{and}\quad V_t(p)=0\]
for any $t\in[0,1]$. 

Consider the flow $\Phi^t\:\Sigma\to \Sigma$ along $V_t$;
that is,
\[\Phi^0(x)=x\quad\text{and}\quad \tfrac{d}{dt}(\Phi^t(x))=V_t(\Phi^t(x))\]
for any $t\in[0,1]$ and $x\in \Sigma$.
The map $\Phi^1\:\Sigma\to \Sigma$ is a diffeomorphism;
in particular, $\Phi^1$ sends the simple path $\alpha$ to a simple path $\alpha_1=\Phi^1\circ\alpha$.
By construction $\alpha_1(1)=q$. 
Since $V_t(p)=0$ for any $t$, we have $\alpha_1(0)\z=p$.
That is, the path $\alpha_1$ runs from $p$ to $q$.

It remains to show that $\alpha_1$ is homotopic to $\alpha*\beta$ relative to the ends.
Set $\alpha_\tau=\Phi^\tau\circ\alpha$ and denote by  $\beta_\tau$ the path running along $\beta$ from $\beta(\tau)$ to $q$;
that is, 
\[\beta_\tau(t)=\beta(\tau+\tfrac1{1-\tau}\cdot t).\]
The concatenation $\alpha_\tau*\beta_\tau$ provides a homotopy from $\alpha*\beta$ to $\alpha_1*\beta_1$. 
Since $\beta_1$ is a constant path, $\alpha*\beta$ is homotopic to $\alpha_1$.
Hence the statement follows.
\qeds

This is a stripped version of the problem suggested by Jaros{\l}aw K\k{e}dra \cite{One-step}; 
it was used by Michael Khanevsky \cite[Lemma 3 in][]{khanevsky}.

\csname @openrightfalse\endcsname
\chapter{Piecewise linear geometry}

A \index{polyhedral space}\emph{polyhedral space} is a complete length-metric space that admits a locally finite triangulation 
such that each simplex is isometric to a simplex in a Euclidean space.
By a {}\emph{triangulation} of a polyhedral space, we always mean a triangulation of that type. 

A point in a polyhedral space is called \index{regular point}\emph{regular} if it has a neighborhood isometric to an open set in a Euclidean space;
otherwise, it is called {}\emph{singular}.

If we replace the Euclidean spaces by the unit spheres or the hyperbolic spaces,
we arrive at the definition of {}\emph{spherical} and {}\emph{hyperbolic polyhedral spaces} respectively.

The term \index{piecewise}\emph{piecewise} typically means that there is a triangulation with some property on each triangle.
For example,  if $P$ and $Q$ are polyhedral spaces, then
\begin{itemize}
\item a map $f\:P\to Q$ is called {}\emph{piecewise distance-preserving} if there is a triangulation $\mathcal{T}$ of $P$ such that for any simplex $\Delta\in \mathcal{T}$ the restriction $f|_\Delta$ is distance-preserving;
\item a map $h\:P\z\to Q$  is called {}\emph{piecewise linear} if both spaces $P$ and $Q$ admit triangulations such that each simplex of $P$ is mapped to a simplex of $Q$ by an affine map.
In particular, a {}\emph{piecewise linear homeomorphism} is a piecewise linear map which is a homeomorphism.\label{piecewise linear map}
\end{itemize}

\subsection*{Spherical arm lemma}\label{Spherical arm lemma}

Recall that a polygon without self-intersections is called \index{simple polygon}\emph{simple}.

\begin{pr}
Let $A=[a_1\dots a_n]$ and $B=[b_1\dots b_n]$ be two simple spherical polygons 
with equal corresponding sides.
Assume that $A$ lies in a hemisphere and $\measuredangle a_i\ge\measuredangle b_i$ for each $i$.
Show that $A$ is congruent to $B$.
\end{pr}

\parit{Semisolution.}
Cut out $A$ from the sphere and glue $B$ in its place.
Denote by $\Sigma$ the obtained spherical polyhedral space.
Note that 
\begin{itemize}
\item $\Sigma$ is homeomorphic to $\mathbb S^2$.
\item $\Sigma$ has curvature $\ge 1$ in the sense of Alexandrov; that is, the total angle around each singular point is less than $2\cdot \pi$.
\item All the singular points of $\Sigma$ 
lie outside of an isometric copy of a hemisphere $\mathbb{S}^2_+\subset \Sigma$.
\end{itemize}

Denote by $n$ the number of singular points in $\Sigma$.
It is sufficient to show that $n=0$.

Assume the contrary; that is, $n\ge 1$.
We can assume that $n$ takes the minimal possible value.

Clearly, $n>1$;
that is, $\Sigma$ cannot have a single singular point.
Therefore we can choose two singular points $p,q\in \Sigma$.
Cut $\Sigma$ along a geodesic $[pq]$.
The obtained hole can be patched so that we obtain a new polyhedral space $\Sigma'$ of the same type but with $n-1$ singular points.
Since $n$ is minimal, we arrive at a contradiction.

Namely, if the total angles around $p$ and $q$ are $2\cdot \pi-\alpha$ and $2\cdot \pi-\beta$ respectively,
consider the spherical triangle $\triangle$ with the base $|p\z-q|_\Sigma$ and the adjacent angles $\tfrac\alpha2$, $\tfrac\beta2$. 
The needed patch is obtained by doubling $\triangle$ along its lateral sides.
\qeds

\parit{Alternative end of the proof.}
By the Alexandrov embedding theorem, $\Sigma$ is isometric to the surface of a convex polyhedron $P$ in the unit sphere $\mathbb S^3$. 
The center of the hemisphere has to lie in a facet of $P$, say $F$.
It remains to note that $F$ contains the equator and therefore $P$ has to be a hemisphere in $\mathbb S^3$ or an intersection of two hemispheres.
In both cases, its surface is isometric to $\mathbb S^2$.
\qeds

The problem is due to Victor Zalgaller \cite{zalgaller-shperical-polygon};
the result of Victor Toponogov in \cite{toponogov} gives a smooth analog of this statement.
The patch construction above was introduced by 
Aleksandr Alexandrov
in his proof of convex embeddability of polyhedra
\cite[see VI, \S7 in][]{alexandrov1948}.
The alternative end of the proof is taken from \cite{panov-petrunin}.

\subsection*{Triangulation of 3-sphere}\label{4-poly}

\begin{pr}
Construct a triangulation of $\mathbb{S}^3$ 
with $100$ vertices
such that any two vertices are connected by an edge.
\end{pr}

\subsection*{Folding problem}\label{Folding problem}

\begin{pr}
Let $P$ be a compact $2$-dimensional 
polyhedral space. 
Construct a 
piecewise distance-preserving map
$f\:P\to \RR^2$.
\end{pr}

\subsection*{Piecewise distance-preserving extension}\label{iso-kirzhbraun}

\begin{pr}
Prove that any 1-Lipschitz map from a finite subset $F\subset \RR^2$
to 
$\RR^2$ can be extended to a 
piecewise distance-preserving map
$\RR^2\to\RR^2$.
\end{pr}

\subsection*{Closed polyhedral surface}\label{Closed polyhedral surface}

\begin{pr}
Construct a closed polyhedral surface $\Sigma$ in $\RR^3$ with nonpositive curvature;
that is, the total angle around each vertex of $\Sigma$ is at least~$2\cdot\pi$.
\end{pr}

\subsection*{Minimal polyhedral disk}\label{Minimal polyhedral disk}

By a polyhedral disk in $\RR^3$
we mean a triangulation of a plane polygon $P$ with a map $P\to\RR^3$ that is affine on each triangle.
The area of the polyhedral disk is defined as the sum of areas of the images of the triangles in the triangulation.

\begin{pr}
Consider the  class of polyhedral disks glued from $n$ triangles in $\RR^3$ 
with a fixed broken line as the boundary.
Let $\Sigma_n$ be a disk of minimal area in this class.
Show that $\Sigma_n$ is a \index{saddle surface}\emph{saddle surface};
that is, a plane cannot cut all the edges coming from one of the interior vertices of~$\Sigma_n$.
\end{pr}

\subsection*{Coherent triangulation\easy}\label{Coherent triangulation} 

A triangulation of a convex polygon is called coherent if there is a convex function that is linear on each triangle and changes its gradient on every edge of the triangulation.

\begin{pr}
Find a non-coherent triangulation of a triangle.
\end{pr}

\subsection*{Sphere with one edge\hard}\label{panov-S^3} 

\begin{pr}
Construct a polyhedral space that is homeomorphic to $\mathbb{S}^3$ and such that its singular locus is formed by a circle.
\end{pr}

\subsection*{Triangulation of a torus}\label{Triangulation of a torus}

\begin{pr}
Show that the torus does not admit a triangulation 
such that one vertex has 5 edges,
one has 7 edges and 
all other vertices have 
6 edges. 
\end{pr}

\subsection*{No simple geodesics\easy}\label{No simple geodesics}

\begin{pr}
Construct a convex polyhedron $P$ whose surface 
does not have a closed simple geodesic.
\end{pr}

\section*{Semisolutions}

\parbf{Triangulation of 3-sphere.}
Choose 100 distinct points $p_1\z\dots,p_{100}$
on the {}\emph{moment curve} 
\[\gamma\:t\mapsto (t,t^2,t^3,t^4)\] 
in $\RR^4$.
Denote by $P$ the convex hull of $\{p_1,\z\dots,p_{100}\}$.

The surface of $P$ is homeomorphic to $\mathbb{S}^2$.
Therefore it is sufficient to show that any two vertices of $P$ are connected by an edge.
The latter follows from the following claim.

\begin{cl}{$({*})$}
Given two points $p$ and $q$ on $\gamma$, there is a hyperplane $H$ in $\RR^4$ that intersects $\gamma$ only at $p$ and $q$ and leaves $\gamma$ on one side.
\end{cl}

To prove the claim, assume that $p=\gamma(t_1)$ and $q=\gamma(t_2)$. 
Consider the polynomial
\[f(t)=a+b\cdot t+c\cdot t^2+d\cdot t^3+t^4=(t-t_1)^2\cdot(t-t_2)^2.\]
Clearly, $f(t)\ge 0$, and the equality holds only at $t_1$ and $t_2$.
It follows that the affine function $\ell\:\RR^4\to\RR$ defined by 
\[\ell\:(w,x,y,z)\mapsto a+b\cdot w+c\cdot x+d\cdot y+z\]
is nonnegative at the points of $\gamma$ and vanishes only at $p$ and $q$.
Therefore the zero-set of $\ell$ is the required hyperplane $H$ in $({*})$. 
\qeds

The polyhedron $P$ above is an example of the so-called \index{cyclic polytope}\emph{cyclic polytopes}.

\parbf{Folding problem.}
Given a triangulation of $P$,
consider the Voronoi domain $V_v$ for each vertex $v$;
that is, $V_v$ is the set of all points in $P$ closer to $v$ than to any other vertex.
Note that the triangulation can be subdivided if necessary
so that the Voronoi domain of each vertex is isometric to a convex subset in the cone with the vertex at its tip.

The boundaries of all the Voronoi domains form a graph with straight edges.
Let us triangulate $P$ so that each triangle has one of those edges as the base 
and the opposite vertex is the center of an adjacent Voronoi domain; 
such a vertex will be called the {}\emph{main} vertex of the triangle.

\begin{wrapfigure}[7]{o}{47 mm}
\vskip-0mm
\centering
\includegraphics{mppics/pic-802}
\end{wrapfigure}

Choose a solid triangle $\triangle=[vab]$ in the constructed triangulation; 
let $v$ be its main vertex.
Given a point 
$x\in  \triangle$, set 
\begin{align*}
\rho(x)&=|x-v|
\intertext{and}
\theta(x)&=\min \{\measuredangle \hinge vax,\measuredangle\hinge vbx\}.
\end{align*}
Let us map $x$ to the point with polar coordinates $(\rho(x),\theta(x))$ in the plane.

Note that for each triangle $\triangle$, 
the constructed map $\triangle\to\RR^2$ is piecewise distance-preserving.
It remains to check that these maps agree on the common sides of the triangles.
\qeds

This construction was given by Victor Zalgaller \cite{zalgaller-polyhedra}.
Svetlana Krat generalized the statement to higher dimensions \cite{krat}.

\parbf{Piecewise distance-preserving extension.}
Let $a_1,\dots,a_n$
and $b_1,\z\dots,b_n$
be two collections of points in $\RR^2$
such that 
\[|a_i-a_j|\ge |b_i-b_j|\] 
for all pairs $i$, $j$.
We need to construct a piecewise distance-preserving map $f\:\RR^2\to\RR^2$
such that $f(a_i)=b_i$ for each $i$.

Assume that the problem is already solved for $n<m$;
let us do the case $n=m$.
By assumption, 
there is a piecewise linear length-preserving map $f\:\RR^2\to\RR^2$
such that $f(a_i)=b_i$ for each $i>1$.

Consider the set 
\[\Omega=\set{x\in\RR^2}{|f(x)-b_1|>|x-a_1|}.\]
Since $|a_i-a_1|\ge|b_i-b_1|$, we get $a_i\notin \Omega$ for $i>1$.

Note that we can assume that the map $f$ and therefore the set $\Omega$ are bounded.
Indeed, let $\square$ be a square containing all the points $b_i$.
There is a piecewise isometric map $h\:\RR^2\to\square$ obtained by folding plane along the lines of the grid defined by $\square$.
Then the composition $h\circ f$ is bounded and it satisfies all the properties of $f$ described above.

If $\Omega=\emptyset$,
then $f(a_1)=b_1$; 
that is, $f$ is a solution.
It remains to consider the case $\Omega\ne\emptyset$. 

Note that $\Omega$ is star-shaped with respect to $a_1$.
Indeed, if $x\in \Omega$, then $|a_1-x|<|b_1-f(x)|$.
If $y\in [a_1x]$ then 
$|a_1-y|+|y-x|=|a_1-x|$ and since $f$ is length-preserving we get $|f(x)-f(y)|\le |x-y|$.
By the triangle inequality, 
$|a_1-y|<|b_1-f(y)|$; that is, $y\in\Omega$. 

\begin{wrapfigure}{o}{47 mm}
\vskip-1mm
\centering
\includegraphics{mppics/pic-804}
\vskip-1mm
\end{wrapfigure}

The boundary $\partial\Omega$ can be subdivided into a finite collection of line segments $\{E_i\}$
so that $f$ maps rigidly each $E_i$.
Note that 
\[|f(x)\z-b_1|=|x-a_1|\]
for any $x\in E_i$.
Denote by $T_i$ the triangle with the base $E_i$ and the vertex $a_1$.
From the above, there is a rigid motion $m_i$ of $T_i$ such that $m_i(x)\z=f(x)$ for any $x\in E_i$ and $m_i(a_1)=b_1$.
Let us redefine the map $f$ in $\Omega$ by sending $x$ to $m_i(x)$ for any $x\in T_i$.

The maps $m_i$ agree on the common sides of triangles $T_i$.
Therefore we have produced a new piecewise isometric map $f'\:\RR^2\to \RR^2$ satisfying all the requirements.
\qeds

The same proof works in all dimensions.

The statement was proved by Ulrich Brehm 
and rediscovered by Arseniy Akopyan and Alexey Tarasov [see \ncite{brehm}, \ncite{akopyan-tarasov}, and Section 2 in \ncite{petrunin-yashinsky}].


\parbf{Closed polyhedral surface.}
An example can be constructed by drilling a polyhedral cave in your favorite convex polyhedron.
On the diagram, you see the result of this construction for the octahedron.

{

\begin{wrapfigure}{r}{43 mm}
\vskip-2mm
\centering
\includegraphics{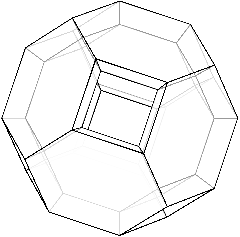}
\bigskip
\includegraphics{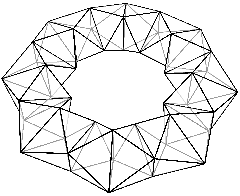}
\end{wrapfigure}

\medskip

Choose a convex polyhedron $K$.
We can assume that the interior of $K$ contains the origin $0\in\RR^3$.
Remove from $K$ the interior of $K'=\tfrac45\cdot K$.

Note that one can drill from each vertex of $K$ a polyhedral tunnel to the corresponding vertex of $K'$
so that the surface $\Sigma$ of the obtained non-convex polytope is a solution.
\qeds

The problem suggested by Jaros{\l}aw K\k{e}dra.

The given construction above produces a surface of genus at least 3.
Another example shown on the diagram is isometric to a flat torus.
It is a bent version of the so-called \index{Schwarz boot}\emph{Schwarz boot} \cite{schwarz1890definition}.
It is made by joining a few identical cylinders, each made from six triangles.

The existence of such a torus also follows from a general result of Yuri Burago and Victor Zalgaller \cite{burago-zalgaller:pl}.
They show in particular that any 1-Lipschitz smooth embedding of the flat torus in $\RR^3$ can be approximated by a piecewise distance-preserving embedding.

}

The following related problem was proposed by Brian Bowditch.

\begin{pr}
Construct a polyhedral metric on the 3-sphere such that the total angle around any edge of its triangulation is at least $2\cdot\pi$.
\end{pr}

A solution can be built using the construction of Joel Hass \cite{hass}.
Another solution was given by Karim Adiprasito \cite{adiprasito};
he proved that an example can be found among spaces that admit a cubulation into unit cubes.
One can also show that it is impossible to find such an example by starting with a doubled cube (as well as other simple polyhedral metrics on the sphere) and passing to coverings branching along unknots а finite number of times.

\parbf{Minimal polyhedral disk.}
Arguing by contradiction, 
assume that a polyhedral disk $\Sigma$ minimizing the area is not saddle;
that is, there is an interior vertex $v$ of $\Sigma$ such that all the edges from $v$ can be cut with a plane.

Note that we can move $v$ in such a way that the lengths of all its edges decrease.

\begin{wrapfigure}{r}{43 mm}
\vskip0mm
\centering
\includegraphics{mppics/pic-806}
\end{wrapfigure}

Since the area is minimal,  this deformation does not decrease the area. 
Taking the derivative of the total area along this deformation implies that $\Sigma$
contains two adjacent non-coplanar triangles $[pvw]$ and $[qvw]$ such that
\[\measuredangle \hinge pvw+\measuredangle \hinge qvw> \pi.\]
In this case, replacing the triangles $[pvw]$ and $[qvw]$
by the triangles $[vpq]$ and $[wpq]$
leads to a polyhedral surface with a smaller area.
That is, $\Sigma$ is not area-minimizing --- a contradiction.
\qeds

{

\begin{wrapfigure}{r}{43 mm}
\vskip-8mm
\centering
\includegraphics{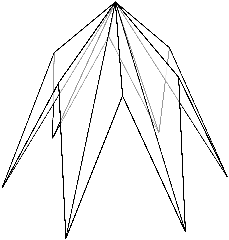}
\end{wrapfigure}

For a general polyhedral surface, a deformation decreasing the lengths of all edges may not decrease the area.
Moreover, the surface that minimizes the area among all surfaces with a fixed  triangulation might not be saddle;
the symmetric tent shown on the diagram provides an example [see \ncite{petrunin-monthly} for more details].


\parbf{Coherent triangulation.} 
An example is shown on the diagram.
The triangulation of the triangle $[x'y'z']$ has a homothetic triangle $[xyz]$ and the edges
$[xx']$, $[yy']$, $[zz']$, 
$[yx']$, $[zy']$, $[xz']$.

}

\medskip

\begin{wrapfigure}{r}{33 mm}
\vskip0mm
\centering
\includegraphics{mppics/pic-808}
\end{wrapfigure}

Assume this triangulation is coherent;
let $f$ be the corresponding piecewise linear convex function.
Without loss of generality, we can assume that $f$ vanishes on the boundary of the big triangle.

From the convexity of $f$ at the edges $[x'y]$,  $[y'z]$, and $[z'x]$, we get 
\[f(x)>f(y)>f(z)>f(x)\]
--- a contradiction.
\qeds

The problem is discussed in the book of 
Israel Gelfand, 
Mikhail Kapranov,
and Andrei Zelevinsky  \cite[see 7C in][]{GKZ}.
The given example is closely related to the so-called \index{Schönhardt polyhedron}\emph{Schönhardt polyhedron}, an example of a non-convex polyhedron that does not admit a triangulation \cite{schoenhardt}.

\parbf{Sphere with one edge.} 
An example can be found among flat orbifolds;
in other words, the required polyhedral space can be chosen to be a quotient of $\mathbb{R}^3$ by a crystallographic action.

\begin{wrapfigure}{r}{23 mm}
\vskip-4mm
\centering
\includegraphics{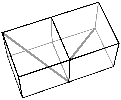}
\end{wrapfigure}

\medskip

Consider the action on $\mathbb{R}^3$ generated by order-3 rotations around two diagonals of the cubes shown on the diagram.
Note that this action is crystallographic;
in fact, it preserves the cubical grid.
Therefore the quotient of $\mathbb{R}^3$ by this action, say $P$, is a compact polyhedral space.

For the described action,
the isotropy group of any point is either trivial or has order 3.
Therefore any point in $P$ admits a neighborhood that is isometric to an open set either in $\RR^3$ or $\RR^3/\ZZ_3$.
Both of these spaces are topological 3-manifolds;
therefore, $P$ is a 3-dimensional manifold as well. 

Note that $P$ is simply-connected;
this follows since the action is generated by rotations. 
By the Poincaré conjecture, $P$ is a topological sphere.

The singular locus of $P$ is the image of the axes of all the order-3 rotations.
Note that the action is transitive on the set of all these axes.
It follows that the singular locus $P_s$ of $P$ is connected; that is, $P_s$ is a circle.

\medskip

This is the so-called P$2_13$ action.
Among 219 crystallographic actions, this is the only one with the quotient space that has the required property; see \cite{dunbar}.

Note that a 3-fold covering of $P$ that is branching in $P_s$ is a flat manifold. 
In particular, this covering is \emph{not} simply-connected.
Therefore $P_s$ is \emph{not} a trivial knot;
it is, in fact, the figure-eight knot.

\begin{wrapfigure}{r}{23 mm}
\vskip-4mm
\centering
\includegraphics{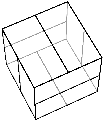}
\vskip-4mm
\end{wrapfigure}

There are a few crystallographic actions with singular locus formed by links. 
The simplest is the Borromean rings;
it is the singular locus of a flat orbifold obtained by gluing each face of a cube to itself
along the reflections with respect to the middle lines shown on the picture; the corresponding action is called I$2_12_12_1$. 

Other examples (nonorbifold, spherical, hyperbolic) are discussed by Michel Boileau and Joan Porti \cite[Chapter 9]{boileau-porti};
more examples are given by Alexander Mednykh \cite{mednykh}.

I get intrigued by this problem because of the following connection:
spaces of directions in 4-dimensional polyhedral Kähler manifolds
are 3-spheres with a spherical polyhedral metric with a singular locus formed by a knot or link.
In addition, these spaces have an isometric $\RR$-action induced by the complex structure.

Polyhedral Kähler manifolds admit a natural stratification into flat manifolds of even codimensions.
So the existence problem of polyhedral Kähler metrics on a given manifold or with a given stratification might be considered as a higher-dimensional analog of the original problem.

A polyhedral metric on a manifold is called \index{polyhedral Kähler manifold}\emph{Kähler} if it comes with a complex structure on each simplex of maximal dimension that agrees on each simplex of codimension 1 (plus a minor condition).
This class of spaces was defined and studied by Dmitri Panov \cite{panov-Kaeler}.
In particular, he constructed many examples of polyhedral Kähler metrics on 
$\CP^2$ ---
try to construct one, it is not that easy.
The following problem is completely open:

\begin{pr}
Is it possible to approximate the canonical metric on $\CP^2$ by polyhedral metrics with nonnegative curvature in the sense of Alexandrov?
\end{pr}

\parbf{Triangulation of a torus.}
Assume the contrary;
let $\tau$ be a triangulation of the torus with the vertex $z_5$ meeting $5$ triangles,
vertex $z_7$ meeting $7$ triangles, 
and every other vertex meeting $6$ triangles.

Let us equip the torus with the flat metric such that each triangle is equilateral.
The metric will have two singular cone points $z_5$ and $z_7$.
The total angle around $z_5$ is $\tfrac53\cdot\pi$
and the total angle around $z_7$ is $\tfrac73\cdot\pi$.
Note that

\begin{cl}{$({*})$}
the holonomy group of the obtained polyhedral metric on the torus is generated by the rotation by $\tfrac\pi3$.
\end{cl}

Indeed, since parallel translation along any loop preserves the directions of the sides of any triangle;
it can only permute it cyclically, which corresponds to rotations by multiples of $\tfrac\pi3$. 
On the other hand, the holonomy of the loop that surrounds $z_5$ is a rotation by $\tfrac\pi3$.

Consider a closed geodesic $\gamma_1$ minimizing the length among all not null-homotopic circles.
Let $\gamma_2$ be another closed geodesic that minimizes the length and is not homotopic to any power of $\gamma_1$.

Note that $\gamma_1$ and $\gamma_2$ intersect at a single point;
otherwise, one could shorten one of them keeping the defining property.

Note that $\gamma_i$ does not contain $z_5$.
In fact, no geodesic can pass thru any singular point with a total angle smaller than $2\cdot\pi$.

Assume that $\gamma_i$ passes thru $z_7$.
Then by $({*})$, one of the two angles cut by $\gamma_i$ at $z_7$ is $\pi$.
It follows that one can push $\gamma_i$ aside so that it does not longer pass thru $z_7$, but remains to be a closed geodesic with the same length.

\begin{wrapfigure}{r}{33 mm}
\vskip0mm
\centering
\includegraphics{mppics/pic-810}
\end{wrapfigure}

Cut the torus along $\gamma_1$ and $\gamma_2$.
In the obtained quadrilateral, connect $z_5$ to $z_7$ by a minimizing geodesic and cut along it.
This way we obtain an annulus $\Omega$ with a flat metric.

Note that a neighborhood of the first boundary component is a parallelogram --- it has equal opposite sides and its angles add up to $2\cdot \pi$.
In particular, $\Omega$ admits an isometric immersion into the plane.

The second boundary component has to be mapped to a diangle with straight sides and angles $\tfrac\pi3$.
Such diangle does not exist in the plane --- a contradiction.
\qeds

The problem was originally discovered and solved by {\fontencoding{T1}\selectfont Stanislav Jendro\v{l}}
and Ernest Jucovi\v{c} \cite{jendrol-jucovich},
their proof is combinatorial.
The solution described above was given by Rostislav Matveyev \cite{matveyev}.
A complex-analytic proof was found by 
Ivan Izmestiev, 
Robert Kusner, 
Günter Rote, 
Boris Springborn, 
and John Sullivan \cite{izmestiev-rote-springborn-kusner}.

\begin{wrapfigure}{r}{20 mm}
\vskip0mm
\centering
\includegraphics{mppics/pic-812}
\end{wrapfigure}

There are flat metrics on the torus with 
only two singular points of total angles $\tfrac53\cdot\pi$ and $\tfrac73\cdot\pi$.
Such an example can be obtained by identifying the hexagon on the picture  according to the arrows.
However, the holonomy group of the obtained torus is generated by the rotation by $\tfrac\pi6$. 
In particular, the observation $({*})$ is essential in the proof.

The same argument shows that the
holonomy group of a flat torus with exactly two singular points of total angle $2\cdot(1\pm \tfrac1n)\cdot\pi$ has more than $n$ elements.
In the solution, we did the case $n=6$.

If one denotes by $v_m$ the number of vertices in a triangulation of the torus with $m$ incoming edges,
then by Euler's formula, we get
\[\sum_m(m-6)\cdot v_m=0.\leqno({*}{*})\]
Note that this equation says nothing about $v_6$.
It turns out that for almost any sequence $v_3,v_4,\dots$ satisfying $({*}{*})$ one can adjust $v_6$ so that it corresponds to a triangulation of the torus --- the sequence 
\[0,0,1,v_6,1,0,0,\dots\] 
discussed in the problem is the only exception. 

The following problem is harder. 
Recall that the curvature of a point $s$ in a polyhedral surface is defined as $2\cdot\pi-\theta$, where $\theta$ denotes the total angle around~$s$.
Note that all regular points in a polyhedral surface have zero curvature.

\begin{pr}
Let $\Sigma$ be a spherical polyhedral space homeomorphic to the 2-sphere
and $\omega_1,\dots,\omega_n$ be the curvatures of its singular points.
Set
\[\delta_i=\min\set{|\tfrac{\omega_i}2-2\cdot k\cdot\pi|}{k\in\ZZ}.\]
Show that there is a closed polygonal line in the unit sphere with sides 
$\delta_1,\dots,\delta_n$.  
\end{pr}

This problem was stated and solved by Gabriele Mondello and Dmitri Panov \cite{mondello-panov}.
The solution requires another holonomy group ---
it assigns an element of the double covering of $\SO(3)$ (which is $\SU(2)\z=\mathbb{S}^3$) to any loop in $\Sigma$ that avoids singularities.

\parbf{No simple geodesics.}
The curvature of a vertex on the surface of a convex polyhedron
is defined as $2\cdot\pi-\theta$, where $\theta$ is the total angle around the vertex.

By the Gauss--Bonnet formula, a simple closed geodesic cuts the surface into two disks each with total curvature $2\cdot\pi$.
Therefore it is sufficient to construct a convex polyhedron with curvatures of the vertices $\omega_1,\dots,\omega_n$ such that
$2\cdot\pi$ cannot be obtained as a sum of some of the $\omega_i$.

An example of that type can be found among the tetrahedrons.
\qeds

The problem is due to Gregory Galperin \cite{galperin};
it was rediscovered by Dmitry Fuchs and Serge Tabachnikov \cite[see 20.8 in][]{fuchs-tabachnikov}.
The following problem is closely related.

\begin{pr}
Assume that the surface of convex polyhedron $P$ contains arbitrary long closed simple geodesics. 
Show that $P$ is an \emph{isosceles tetrahedron};
that is, the opposite edges of the tetrahedron are equal.
\end{pr}

The latter statement was proved by Vladimir Protasov [see \ncite{protasov} and also \ncite{akopyan-petrunin}, \ncite{itoh-vilcu}].

\csname @openrightfalse\endcsname
\chapter{Discrete geometry}

In this chapter, we consider geometrical problems with a strong combinatoric flavor.
No special prerequisite is needed.

\subsection*{Round circles in 3-sphere}\label{Round circles}

\begin{pr}
Suppose that $\mathcal{C}$ is a finite collection of pairwise linked round circles in the unit 3-sphere. 
Prove that there is an isotopy of $\mathcal{C}$ that moves all of them into great circles.
\end{pr}

\parit{Semisolution.}
For each circle in $\mathcal{C}$, consider the plane containing it.
Note that the circles are linked 
if and only if 
the corresponding planes intersect at a single point inside the unit sphere $\mathbb{S}^3\subset \RR^4$.

Consider the collection of circles formed by the intersections of the planes with the sphere of radius $R\ge 1$.
Rescale the sphere and pass to the limit as $R\to\infty$.  
This way we get the needed isotopy.\qeds

This problem was discussed 
by Genevieve Walsh \cite{walsh}.
The same idea was used by Michael Freedman and  Richard Skora to show that any link made from pairwise not linked round circles is trivial; 
in particular, Borromean rings cannot be realized by round circles
\cite[see Lemma 3.2 in ][]{freedman-skora}.

\subsection*{Box in a box}\label{box-in-box}

\begin{pr}
Suppose a rectangular parallelepiped with sides $a,b,c$ 
lies inside another rectangular parallelepiped with sides $a',b',c'$. 
Show that 
\[a'+b'+c'\ge a+b+c.\]

\end{pr}

\subsection*{Piercing the cube}\label{Piercing the cube}

\begin{pr}
Let $\Pi$ be a $k$-plane that passes thru the center of a unit $n$-cube $Q$.
Show that $k$-area of the intersection $\Pi\cap Q$ is at least 1.
\end{pr}

\subsection*{Harnack's circles}\label{Harnack}

\begin{pr}
Prove that a smooth algebraic curve of degree $d$ in $\RP^2$ consists of at most $n=\tfrac12\cdot(d^2-3\cdot d+4)$ connected components.
\end{pr}

\subsection*{Two points on each line}\label{2pts-on-line}

\begin{pr}
Construct a set in the Euclidean plane that intersects each line at exactly 2 points. 
\end{pr}

\subsection*{Balls without gaps}
\label{Balls without gaps}

\begin{pr}
Let $B_1,\dots,B_n$ be balls  
of radii $r_1,\dots,r_n$ 
in a Euclidean space.
Assume that no hyperplane divides the balls into two
non-empty sets without intersecting at least one of the balls. 
Show that the balls
$B_1,\dots,B_n$ can be covered by a ball of radius
$r=r_1+\dots+r_n$.

\end{pr}

\subsection*{Covering lemma}

\begin{pr}
Let $\{B_i\}_{i\in F}$ be any finite collection of  balls in $\RR^m$. 
Show that there is a subcollection of pairwise disjoint balls $\{B_i\}_{i\in G}$, $G\subset F$
such that
\[\vol \left(\bigcup_{i\in F}B_i\right)\le 3^m\cdot \vol \left(\bigcup_{i\in G}B_i\right).\]
\end{pr}

{


\begin{wrapfigure}{r}{27 mm}
\vskip0mm
\centering
\includegraphics{mppics/pic-902}
\end{wrapfigure}

\subsection*{Kissing number\easy}\label{pr:Kissing number}

Let  $W_0$ be a convex body in $\RR^m$.
We say that $k$ is the \index{kissing number}\emph{kissing number} of $W_0$ (briefly $k=\mathop{\rm kiss}W_0$)
if $k$ is the maximal integer such that there are $k$ bodies $W_1,\dots,W_k$ such that 
(1) each $W_i$ is congruent to $W_0$,
(2) $W_i\cap W_0\not=\emptyset$ for each $i$ 
and (3) no pair $W_i,W_j$ has common interior points.

}

As you may have guessed from the diagram, the kissing number of the round disk in a plane is $6$.

\begin{pr}
Show that for any convex body $W_0$ in $\RR^m$ we have that
$$\mathop{\rm kiss}W_0\ge \mathop{\rm kiss}B,$$
where $B$ denotes the unit ball in $\RR^m$.
\end{pr}

\subsection*{Monotonic homotopy}
\label{mono-homotopy}

\begin{pr}
Let $F$ be a finite set and $h_0,h_1\: F\to\RR^m$ be two maps.
Consider $\RR^m$ as a subspace of $\RR^{2\cdot m}$.
Show that there is a homotopy  $h_t\:F\z\to\RR^{2\cdot m}$ from $h_0$ to $h_1$ such that  the function 
$t\z\mapsto |h_t(x)-h_t(y)|$
is monotonic for any pair $x,y\in F$.
\end{pr}

\subsection*{Facet cover}\label{Facet cover}

\begin{pr}
Show that any facet of a convex polyhedron can be covered by the remaining facets.
\end{pr}

\subsection*{Cube}\label{Cube}

\begin{pr}
Half of the vertices 
of an $m$-dimensional cube
are colored in white and the other half in black.
Show that the cube has at least $2^{m-1}$ edges connecting vertices of different colors. 
\end{pr}

\subsection*{Geodesic loop}\label{Geodesic loop}

\begin{pr}
Show that the surface of a cube in $\RR^3$
does not admit a geodesic loop with a vertex as the base point.
\end{pr}

\subsection*{Right and acute triangles}\label{Right and acute triangles}

\begin{pr}
Let $x_1,\dots,x_n\in\RR^m$
be a collection of points such that any triangle $[x_ix_jx_k]$ is right or acute.
Show that $n\le 2^m$.
\end{pr}

\subsection*{Upper approximant}\label{One-sided approximants}

\begin{pr}
Let $\mu$ be a Borel probability measure on the plane.
Show that given $\eps>0$, there is a finite set of points $S$ that intersects every convex figure of measure at least $\eps$.
Moreover, we can assume that $|S|$ --- the number of points in $S$, depends only on $\eps$ (in fact, one can take $|S|=\lceil\tfrac1{\eps^5}\rceil$).
\end{pr}

\subsection*{Right-angled polyhedron\thm}\label{Right-angled polyhedron}

A polyhedron is called {}\emph{right-angled} if all its dihedral angles are right.

\begin{pr}
Show that in all sufficiently large dimensions, there is no compact convex hyperbolic right-angled polyhedron. 
\end{pr}

Here is a summary of the Dehn--Sommerville equations that can help.

Let $P$ be a \index{simple polyhedron}\emph{simple} Euclidean $m$-dimensional polyhedron;
that is, exactly $m$ facets meet at each vertex of $P$.
Denote by $f_k$ the number of $k$-dimensional faces of $P$;
the array of integers $(f_0,\dots f_m)$ is called the $f$-vector of $P$.

Choose a linear function $\ell$ that takes different values on the vertices of~$P$.
The \index{index of vertex}\emph{index} of a vertex $v$ 
is defined as the number of edges $[vw]$ of $P$ such that $\ell(v)>\ell(w)$. 
The number of vertices of index $k$ will be denoted by $h_k$.
The array of integers $(h_0,\dots h_m)$ is called the $h$-vector of $P$.

Each $k$-face of $P$ contains a unique vertex that maximizes $\ell$.
If the vertex has index $i$,
then $i\ge k$, and then it is the maximal vertex of exactly $\tbinom i k$
faces of dimension $k$.
This observation can be packed in the following polynomial identity:
\[\sum_k h_k\cdot (t+1)^k=\sum_k f_k\cdot t^k.\]

This identity implies that the $h$-vector does not depend on the choice of $\ell$.
Changing the sign of $\ell$, we get that, the $h$-vector is the same for the reversed order;
that is,
\[h_k=h_{m-k}\quad\text{for all}\quad k.
\leqno({*})\]

The identities $({*})$ for all $k$ are called the Dehn--Sommerville equations.
They give a complete list of linear equations for $h$-vectors (and therefore $f$-vectors) of simple polyhedrons.
For more on the subject, see \cite[Chapter 9]{gruenbaum}.

\subsection*{Real roots of random polynomial\easy}\label{Real roots of random polynomial}

Consider the \index{moment curv}\emph{moment curve} $\gamma_n\:t\mapsto(1,t,\dots,t^n)$  in $\RR^{n+1}$.
Let 
\[\ell_n=\length \tfrac{\gamma_n}{|\gamma_n|}.\]

\begin{pr}
Show that $\tfrac{\ell_n}{\pi}$ is the expected number of real roots of a polynomial of degree $n$ with independent normally distributed real coefficients.
\end{pr}

\subsection*{Space coloring\easy}\label{Space coloring}

\begin{pr}
Let $A$ be a set of colored points in $\mathbb{R}^d$.
We are allowed to color any line containing at least $k$ already colored points.
Suppose that in a finite number of steps, we can color any point in $\mathbb{R}^d$.
What is the minimal number of points in $A$ as a function of $d$ and $k$?
\end{pr}

\section*{Semisolutions}

\parbf{Box in a box.} 
Let $\Pi$ be a parallelepiped
with dimensions $a$, $b$, and $c$.
Denote by $v(r)$ the volume of the $r$-neighborhood of $\Pi$,
 
Note that for all positive $r$ we have
\[v_{\Pi}(r)=w_3(\Pi)+w_2(\Pi)\cdot r+w_1(\Pi)\cdot r^2+w_0(\Pi)\cdot r^3,\leqno({*})\]
where 
\begin{itemize}
\item $w_0(\Pi)=\tfrac43\cdot \pi$ is the volume of the unit ball,
\item $w_1(\Pi)=\pi\cdot (a+b+c)$,
\item $w_2(\Pi)=2\cdot(a\cdot b+b\cdot c+c\cdot a)$ is the surface area of $\Pi$,
\item $w_3(\Pi)=a\cdot b\cdot c$ is the volume of $\Pi$,
\end{itemize}

Let $\Pi'$ be another parallelepiped
with dimensions $a'$, $b'$ and~$c'$.
If $\Pi\subset \Pi'$,
then $v_{\Pi} (r)\le v_{\Pi'}(r)$ for any $r$.
For $r\to\infty$, these inequalities imply
\[a+b+c\le a'+b'+c'.\qedsin\]

\parit{Alternative proof.}
Note that the average length of the projection of $\Pi$ to a line is
$\Const\cdot(a+b+c)$ for some $\Const>0$.
(In fact, $\Const=\tfrac12$, but we will not need it.)

Since $\Pi\subset \Pi'$,
the average length of the projection of $\Pi$
cannot exceed the average length of the projection of $\Pi'$.
Hence the statement follows.
\qeds

The problem was discussed by Alexander Shen \cite{shen}.

A formula analogous to $({*})$
holds for an arbitrary convex body $B$ of arbitrary dimension $m$.
It was discovered by Jakob Steiner \cite{steiner}.
The coefficient $w_i(B)$ in the polynomial (with different normalization constants)
appears under different names, most commonly
\emph{intrinsic volumes} and
\emph{quermassintegrals}.
Up to a normalization constant,
they can also be defined as the average 
area of the projections of $B$ to the $i$-dimensional planes.
In particular, 
if $B'$ and $B$ are convex bodies such that $B'\subset B$,
then $w_i(B')\le w_i(B)$ for any $i$.
This generalizes our problem quite a bit.
Further generalizations lead to the theory of \index{mixed volumes}\emph{mixed volumes} \cite{burago-zalgaller}.

The equality $w_1(\Pi)=\pi\cdot (a+b+c)$ still holds for all parallelepipeds, not only rectangular ones.
In particular, if one parallelepiped 
lies inside another then the sum of all edges of the first one cannot exceed the sum for the second.

\parbf{Piercing the cube.}
Observe that there is an odd increasing 1-Lipschitz function $\phi\:\RR\to(-\tfrac12,\tfrac12)$ that pushes the measure with the density $e^{-\pi x^2}$ to the Lebesgue measure on $(-\tfrac12,\tfrac12)$.

We can assume that $Q=\set{(x_1,\dots,x_n)\in \RR^n}{|x_i|\le \tfrac12}$.
Applying $\phi$ to each coordinate,
we get a 1-Lipschitz diffeomorphism from $\RR^n$ to the interior of $Q$ that pushes the measure with the Gauss density 
$\rho(x_1,\dots,x_n)=e^{-\pi (x_1^2+\dots+x_n^2)}$ to the Lebesgue measure on $Q$.

Note that the inverse image $\Pi'=\phi^{-1}(\Pi)$ is a symmetric $k$-surface in~$\RR^n$.
Since $\phi$ is $1$-Lipschitz, we get that the area of $\Pi\cap Q$ cannot be smaller than the $\rho$-weighted $k$-area of $\Pi'$.
Moreover, we have equality for coordinate $k$-planes.

Denote by $S_r$ the sphere of radius $r$ in $\RR^n$ centered at the origin.
Note that $\Pi'\cap S_r$ intersects each great $(n-k-1)$-sphere in $S_r$.
By Crofton's formula, we get that $(k-1)$-area of $\Pi'\cap S_r$ cannot be smaller than $(k-1)$-area of the great $(k-1)$-sphere in $S_r$.
By the coarea formula, the $\rho$-weighted $k$-area of $\Pi'$ cannot be smaller than $\rho$-weighted $k$-area of a $k$-plane that passes thru the origin --- hence the result.
\qeds

The problem was posed by Anton Good
and solved by Jeffrey Vaaler \cite{vaaler}.
The presented solution was found by Arseniy Akopyan, Alfredo Hubard, and Roman Karasev \cite{akopyan-hubard-karasev};
note that it proves the following more general statement:
\textit{If $f\:Q \to \RR^{n-k}$ is an odd continuous map, such that $Z = f^{-1}\{0\}$ is smooth $k$-surface, then the $k$-area of $Z$ is at least 1.}
The map $\RR^n\to Q$ was first used by Gilles Pisier \cite[p. 182]{pisier}

\parbf{Harnack's circles.}
Let $\sigma\subset \RP^2$ be an algebraic curve of degree~$d$.
Consider the complexification $\Sigma\subset \CP^2$ of~$\sigma$.
Without loss of generality, we may assume that $\Sigma$ is regular.

Note that all regular complex algebraic curves of degree $d$ in $\CP^2$
are isotopic to each other in the class of regular algebraic curves of degree~$d$.
Indeed, the set of equations of degree $d$ that correspond to singular curves has real codimension 2.
It follows that the set of equations of degree $d$ that correspond to regular curves is connected.
Therefore, one can construct an isotopy from one regular curve to any other by changing continuously the parameters of the equations.

In particular, it follows that all regular complex algebraic curves of degree $d$ in $\CP^2$ have the same genus,
denote it by $g$.
Perturbing a singular curve formed by $d$ lines in $\CP^2$,
we can see that 
\[g=\tfrac12\cdot(d-1)\cdot(d-2).\]

The real curve $\sigma$ forms the fixed point set in $\Sigma$ by the complex conjugation.
In particular $\Sigma\setminus\sigma$ has at most two connected components.
Therefore, the number of components of $\sigma$ cannot exceed $g+1$.
\qeds

This problem is a background for Hilbert's 16th problem.
The inequality was originally proved 
by Axel Harnack using a different method \cite{harnack}.
The idea to use complexification is due to Felix Klein \cite{klein}.
In fact any number of connected components up to $g+1$ is realizable, with one exception:
if $d$ is odd, then $\sigma$ has at least one connected component. 

\parbf{Two points on each line.}
Take any complete ordering of the set of all lines 
so that each beginning interval has cardinality less than continuum.

Assume we have a set of points $X$ of cardinality less than continuum such that each line intersects $X$ in at most $2$ points.

Choose the least line $\ell$ in the ordering that intersects $X$ 
in $0$ or $1$ point.
Note that the set of all lines intersecting $X$ at two points has cardinality less than continuum.
Therefore we can choose a point on $\ell$ and add it to $X$ so that the remaining lines are not overloaded.

It remains to apply the well-ordering principle.
\qeds

This problem has an endless list of variations.
The following problem looks similar but far more involved;
a solution follows from the proof of Paul Monsky that a square cannot be cut into triangles with equal areas \cite{monsky}.

\begin{pr}
Subdivide the plane into three everywhere dense sets $A$, $B$, and $C$ such that each line meets exactly two of these sets.
\end{pr}

\parbf{Balls without gaps.} 
Assume the mass of each ball is proportional to its radius.
Denote by $z$  the center of mass of the balls.
It is sufficient to show the following.
\begin{cl}{$({*})$}
The ball $B(z,r)$ contains all $B_1,\dots,B_n$.
\end{cl}

Assume this is not the case.
Then there is a line $\ell$ thru $z$, 
such that the orthogonal projection of a ball $B_i$ to $\ell$ 
does not lie completely inside the projection of $B$.
(This observation reduces the problem to the one-dimensional case.)

Note that the projection of all balls $B_1,\dots,B_n$ has to be connected and it contains a line segment longer than $r$ on one side from $z$. 
In this case, the center of mass of the balls projects inside of this segment --- a contradiction.
\qeds

The statement was conjectured by Paul Erdős.
The solution was given by Adolph and Ruth Goodmans
[see \ncite{goodman-goodman} and also \ncite{hadwiger}].


\parbf{Covering lemma.} 
The required collection $\{B_i\}_{i\in G}$ is constructed using the \emph{greedy algorithm}. 
We choose the balls one by one;
on each step we take the largest ball that does not intersect those which we have chosen already.

\medskip

Note that each ball in the original collection $\{B_i\}_{i\in F}$ intersects a ball in $\{B_i\}_{i\in G}$ with a larger radius.
Therefore 
\[\bigcup_{i\in F}B_i\subset \bigcup_{i\in G}3\cdot B_i,\leqno({*})\]
where $3\cdot B_i$ denotes the ball concentric to $B_i$ and three times larger radius.
Hence the statement follows.
\qeds

The constant $3^m$ can be improved slightly \cite{domotorp}.
For $m=1$ the optimal constant is $2$.
Possibly, for any $m$, the optimal constant is $2^m$;
it can not be smaller, an example can be found among collections of unit balls that contain a fixed point.

The inclusion $({*})$ is called the \emph{Vitali covering lemma}.
The following statement is called the \emph{Besicovitch covering lemma};
it has a similar proof.

\begin{pr}
For any positive integer $m$, there is a positive integer $M$ such that 
any finite collection of balls $\{B_i\}_{i\in F}$ in $\RR^m$ 
contains a subcollection $\{B_i\}_{i\in G}$
such that (1) center of any ball in $\{B_i\}_{i\in F}$ lies inside one of a ball from $\{B_i\}_{i\in G}$
and (2) the collection $\{B_i\}_{i\in G}$ can be subdivided into $M$ subcollections of pairwise disjoint balls.
\end{pr}

Both lemmas were used to prove the so-called \emph{covering theorems} in measure theory,
which state that ``undesirable sets'' have vanishing measures.
Their applications overlap but aren't identical, the \emph{Vitali covering theorem} works for nice measures in arbitrary metric spaces while the \emph{Besicovitch covering theorem} works in nice metric spaces with arbitrary Borel measures.

More precisely, Vitali works in arbitrary metric spaces with a \index{doubling measure}\emph{doubling measure} $\mu$;
the latter means that 
\[\mu [2\cdot B]\le C\cdot \mu B\] 
for a fixed constant $C$ and any ball $B$ in the metric space.
On the other hand, Besicovitch works for all Borel measures in the so-called \emph{directionally limited} metric spaces \cite[see 2.8.9 in][]{federer};
these include Alexandrov spaces with curvature bounded below.


\parbf{Kissing number.}
Set $n=\mathop{\rm kiss} B$.
Let $B_1,\dots, B_n$ be copies of the ball $B$ that touch $B$ and don't have common interior points.
For each $B_i$ consider the vector $v_i$ from the center of $B$ to the center of $B_i$.
Note that $\measuredangle(v_i,v_j)\ge \tfrac\pi3$ if $i\ne j$.

For each $i$,
consider the supporting hyperplane $\Pi_i$
of $W$
with the outer normal vector $v_i$.
Denote by $W_i$ the reflection of $W$ with respect to $\Pi_i$.

\begin{wrapfigure}{r}{43 mm}
\vskip-0mm
\centering
\includegraphics{mppics/pic-903}
\vskip-2mm
\end{wrapfigure}

Note that $W_i$ and $W_j$ have no common interior points if $i\ne j$;
the latter gives the needed inequality.
\qeds

The proof is given by 
Charles Halberg, 
Eugene Levin, 
and Ernst Straus 
\cite{halberg-levin-straus}.
It is not known if the same inequality holds for the orientation-preserving version of the kissing number.

\parbf{Monotonic homotopy.}
Note that we can assume
that $h_0(F)$ and $h_1(F)$ both lie in the coordinate $m$-spaces of $\RR^{2\cdot m}=\RR^m\times \RR^m$;
that is,
$h_0(F)\z\subset\RR^m\times\{0\}$
and $h_1(F)\subset  \{0\}\times\RR^m$.

Direct calculations show that the following homotopy is monotonic
\[h_t(x)=\bigl(h_0(x)\cdot \cos\tfrac{\pi\cdot t}2
\,,\,
 h_1(x)\cdot\sin\tfrac{\pi\cdot t}{2}\bigr).\qedsin\] 
\medskip

This homotopy was discovered by Ralph Alexander \cite{ralexander}.
It has a  
number of applications, 
one of the most beautiful is given 
by K\'aroly Bezdek 
and Robert Connelly \cite{bezdek-connelly};
they proved the Kneser--Poulsen  
and Klee--Wagon conjectures in the two-dimensional case.

The dimension $2\cdot m$ is optimal;
that is, for any positive integer $m$,
there are two maps $h_0,h_1\:F\to \RR^m$ that cannot be connected by a monotonic homotopy $h_t\:F\to\RR^{2\cdot m-1}$.
The latter was shown by Maria Belk and Robert Connelly \cite{belk-connelly}

\parbf{Facet cover.}
Choose a convex polyhedron $P$ and its facet $F$.
Denote by $\Pi$ the hyperplane of $F$.

For each point $p\in F$ consider the maximal ball $B_p\subset P$ that contains~$p$.
Note that $B_p$ touches another facet of $P$ at some point $q$.

Show that the restriction of the partially defined map $q\mapsto p$ to any other facet $F'$ can be extended to a distance-preserving map $F'\to \Pi$.
By construction union of all such images cover $F$.
\qeds

This problem was considered by Igor Pak and Rom Pinchasi \cite{pak-pinchas};
the presented proof was given by Arseniy Akopyan \cite{akopyan-2012}.
A slightly different version of this problem made it to the All-Russian olympiad of 2012 \cite[№ 116774]{problems}.

In the three-dimensional case,
a more involved, but straightforward solution can be built on the fact that orthogonal projection of a convex figure can be covered by the figure itself.
The latter statement is not as simple as one might think  \cite[see][and the references therein]{kos-toroscik}.

\parbf{Cube.}
Consider the cube $[-1,1]^m\subset \RR^m$.
Any vertex of this cube has the form $\bm{q}=(q_1,\dots,q_m)$,
where  $q_i=\pm1$.

For each vertex $\bm{q}$,
consider the intersection of the corresponding hyperoctant with the unit sphere;
that is, consider the set
\[V_{\bm{q}}=\set{(x_1,\dots,x_m)\in\mathbb{S}^{m-1}}{q_i\cdot x_i\ge 0\ \text{for each}\ i}.\]

Let $\mathcal{A}\subset\mathbb{S}^{m-1}$ be the union of all the sets $V_{\bm{q}}$ for black $\bm{q}$.
Note that 
\[\vol_{m-1}\mathcal{A}=\tfrac12\cdot\vol_{m-1}\mathbb{S}^{m-1}.\]

By the spherical isoperimetric inequality,
\[\vol_{m-2}\partial\mathcal{A}
\ge \vol_{m-2}\mathbb{S}^{m-2}.\] 

It remains to observe that
\[\vol_{m-2}\partial\mathcal{A}
=
\tfrac k{2^{m-1}}\cdot\vol_{m-2}\mathbb{S}^{m-2},\]
where $k$ is the number of edges of the cube with one black end and the other in white.
\qeds

The problem was suggested by Greg Kuperberg.

\parbf{Geodesic loop.}
Let $\gamma$ be a geodesic from vetrex to vertex;
denote by $v$ its midpoint.

Show that there is symmetry of the cube that fixes $v$, reverts $\gamma$, and moves all the vertices of the cube.
Conclude that the endpoints of $\gamma$ are different.
\qeds

I learned this problem from Jaros{\l}aw K\k{e}dra; 
it was rediscovered independently by
Diana Davis,
Victor Dods,
Cynthia Traub,
and Jed Yang \cite{DDTY}.
The presented solution is taken from the paper of Serge Troubetzkoy \cite{troubetzkoy}.
This idea can be used to solve the following harder problems.

\begin{pr}
Show the same for the surface of the $n$-dimensional cube, $n\ge 4$.
\end{pr}

\begin{pr}
 Show the same for the surface of the tetrahedron, octahedron, and icosahedron.
\end{pr}

For the dodecahedron such loops exist;
\begin{figure}[!ht]
\vskip0mm
\centering
\includegraphics{mppics/pic-905}
\end{figure}
a development of one example is on the diagram.
Vertices of an inscribed tetrahedron are circled.
A classification of all such loops is found by Jayadev S. Athreya, David Aulicino, and W. Patrick Hooper \cite{athreya-aulicino-hooper}.

These and related problems are discussed by Dmitry Fuchs \cite{fuchs-2016}.

\parbf{Right and acute triangles.}
Denote by $K$ the convex hull of $\{x_1,\z\dots,x_n\}$.
Without loss of generality, we can assume that $\dim K=m$. 

\begin{wrapfigure}{o}{35 mm}
\vskip-0mm
\centering
\includegraphics{mppics/pic-907}
\end{wrapfigure}

Note that for any distinct points $x_i$, $x_j$,
and any interior point $z$ in $K$,
we have 
\[\measuredangle \hinge{x_i}{x_j}z<\tfrac\pi2.\leqno({*})\]
Indeed, if $({*})$ does not hold, then $\langle x_j-x_i,z-x_i\rangle<0$.
Since $z\in K$ we have $\langle x_j-x_i,x_k-x_i\rangle\z<0$ for some vertex $x_k$.
That is, $\measuredangle \hinge{x_i}{x_j}{x_k} > \tfrac\pi2$ --- a contradiction.

Denote by $h_i$ the homothety with center $x_i$ and coefficient $\tfrac12$.
Set $K_i=h_i(K)$.

Let us show that $K_i$ and $K_j$ have no common interior points.
Assume the contrary; 
that is, \[z=h_i(z_i)=h_j(z_j);\]
for some interior points $z_i$ and $z_j$ in $K$.
Note that 
\[
\measuredangle\hinge{x_i}{x_j}{z_j}
+
\measuredangle\hinge{x_j}{x_i}{z_i}
=
\pi,
\]
which contradicts $({*})$.

Note that $K_i\subset K$ for any $i$;
it follows that 
\begin{align*}
\tfrac n{2^m}\cdot \vol K
&=\sum_{i=1}^n\vol K_i\le\vol K.
\end{align*}
Hence the result follows.
\qeds

The problem was posted by Paul Erdős \cite{erdos}
and solved by Ludwig Danzer and Branko Grünbaum \cite{danzer-gruenbaum}.

Grigori Perelman noticed that the same proof works for a similar problem in Alexandrov spaces \cite{perelman-Erdos}.
The latter led to interesting connections with the crystallographic groups \cite{lebedeva};
in particular, it gives an approach to the following open problem.

\begin{pr}
Let $\Gamma\acts\RR^n$ be an effective properly discontinuous isometric action.
Denote by $m$ the number of maximal finite subgroups $\Gamma$ up to conjugation. 
Is it true that $m\le 2^n$?
\end{pr}

Surprisingly, the maximal number of points that make only acute triangles grows exponentially with $m$ as well.
The latter was shown by Paul Erdős and Zolt\'an Füredi \cite{erdos-fueredi} using the \index{probabilistic method}\emph{probabilistic method}.
Later, an elementary constructive argument was found and improved by Dmitriy Zakharov,
\texttt{grizzly} (an anonymous mathematician),
Bal{\'a}zs Gerencs{\'e}r, and Viktor Harangi
\cite{zakharov,grizzly,gerencser-harangi};
the current lower bound is $2^{m-1}+1$, which is exponentially optimal.

\parbf{Upper approximant.}
We assume that the measure is given by a distribution.
The general case is done by a straightforward modification.

Cut the plane by two lines into 4 angles of equal measure.
Let $p$ be the intersection point of the two lines.
Note that every convex set avoiding $p$ is fully contained in three angles out of the four.
In particular, its measure cannot exceed $\tfrac34$.

Apply this construction recursively for the restriction of the measure to each triple of angles.
After $n$ steps we get a $(1+\z\dots+4^{n-1})$-point set 
that intersects each convex figure $F$ of measure $(\tfrac34)^n$.
\qeds

This is a stripped version of a theorem proved by Boris Bukh and Gabriel Nivasch \cite{bukh-nivasch}.

\parbf{Right-angled polyhedron.}
Let $P$ be a right-angled hyperbolic polyhedron of dimension $m$.
Note that $P$ is simple; 
that is, exactly $m$ facets meet at each vertex of $P$.

From the projective model of the hyperbolic plane, 
one can see that for any simple compact hyperbolic polyhedron, there is a simple Euclidean polyhedron with the same combinatorics. 
In particular, the Dehn--Sommerville equations hold for $P$.

Denote by $(f_0,\dots f_m)$ and $(h_0,\dots h_m)$ the $f$- and $h$-vectors of $P$.
Recall that 
\begin{align*}
h_i&=h_{m-i},
\\
f_1&=h_1+2\cdot h_2+\dots+m\cdot h_m,
\\
f_2&= h_2+3\cdot h_3+\dots+\tbinom m 2\cdot h_m.
\end{align*}
Observe that $h_i\ge 1$; together with the identities above it implies that
\[f_2> \tfrac{m-2}4\cdot f_1.
\leqno({*})\]

Since $P$ is hyperbolic, each 2-dimensional face of $P$ has at least 5 sides.
It follows that
\[f_2\le \tfrac{m-1}5\cdot f_1.\]
The latter contradicts $({*})$ for $m\ge 6$.
\qeds

This is the core of the proof of nonexistence of compact hyperbolic Coxeter's polyhedrons of large dimensions 
given by Ernest Vinberg \cite{vinberg, vinberg-strong}.

Playing a bit more with the same inequalities, 
one gets the nonexistence of right-angled hyperbolic polyhedrons,
in all dimensions starting from~$5$.
In the $4$-dimensional case,
there is a regular right-angled  hyperbolic polyhedron with \index{120-cells}\emph{120-cells} --- a $4$-dimensional uncle of the dodecahedron.

The following related question is open:

\begin{pr}
Let $m$ be a large integer.
Is there a cocompact properly discontinuous isometric action  on the $m$-dimensional Lobachevsky space that is generated by finite order elements (for example, involutions)?
\end{pr}

\parbf{Real roots of random polynomial.}
Choose a polynomial 
\[p(t)=a_0+\dots+a_n\cdot t^n.\]
Consider the hyperplane $\Pi$ in $\RR^{n+1}$ defined by the equation 
\[a_0\cdot x_0+\dots+a_n\cdot x_n=0.\]
Note that the number of real roots of $p$ equals the number of intersections of $\Pi$ with the moment curve.

It remains to apply the spherical Crofton formula.
\qeds

The observation is due to Alan Edelman and Eric Kostlan \cite{edelman-kostlan}.

\parbf{Space coloring.}
The answer is $\binom{d+k-1}{d}$.

Choose $d+k-1$ hyperplanes in $\RR^d$ in general position.
Let $A$ be the set of all $\binom{d+k-1}{d}$ intersection points of every $d$-tuple of the chosen hyperplanes. 
Show that starting with $A$ one can color every point in $\RR^d$.

To prove the lower bound, note that $\binom{d+k-1}{d}$ is the dimension of the space of polynomials $p\:\RR^d\to\RR$ of total degree at most $k-1$.
Therefore, if $A$ has less than $\binom{d+k-1}{d}$ points, then there is a non-zero polynomial $p$ of degree at most $k-1$ such that $A$ lies in its zero-set $Z$.
Observe that if a line does not lie in $Z$, then it has at most $k-1$ common points with $Z$.
Therefore, it is impossible to color a point outside of $Z$.
\qeds

This is a problem of Ivan Mitrofanov and Fedor Petrov \cite[№ 10]{kanel-belov};
a three-dimensional version of this problem appeared at the 28$^\text{th}$ Annual Vojtěch Jarník International Mathematical Competition, Category II \cite{vjimc}.

The problem illustrates the so-called \emph{polynomial method};
for more on the subject check the book of Lary Guth \cite{guth2016}.



\backmatter
\newpage
\phantomsection
{\scriptsize
\input{problems.ind}
}
\sloppy
\printbibliography[heading=bibintoc]
\fussy

\end{document}